\documentclass[letterpaper, 12pt, oneside]{book}

    \usepackage[left=1.5in, right=1in, top=1in, bottom=1in, dvips]{geometry}
    
    \usepackage{fancyhdr} 
    \renewcommand{\bibname}{References}

    \usepackage{amsmath,amsthm, amsfonts,amssymb} 
    \usepackage{setspace}

    \usepackage{graphicx} 
    \usepackage{pstricks, pst-plot}

    \usepackage{lineno} 
    \usepackage{appendix}

    \newtheorem{theorem}{Theorem}[section]
    \newtheorem{proposition}[theorem]{Proposition}
    \newtheorem{lemma}[theorem]{Lemma}
    \newtheorem{sublemma}[theorem]{Sublemma}
    \newtheorem{corollary}[theorem]{Corollary}
    \newtheorem{claim}[theorem]{Claim}
    \newtheorem{conjecture}[theorem]{Conjecture}
    \newtheorem{criterion}[theorem]{Criterion}
    \theoremstyle{definition}
        \newtheorem{definition}[theorem]{Definition}
        \newtheorem{example}[theorem]{Example}
        \newtheorem{algorithm}[theorem]{Algorithm}

    \theoremstyle{remark}
        \newtheorem{remark}[theorem]{Remark}

    \numberwithin{equation}{section}
    \numberwithin{figure}{section}

    \newcommand{\n}{\vspace{12pt}} 
      
    \newcommand{\newchap}[1] 
	{                    
        \chapter{#1}
        \thispagestyle{myheadings}
	}

    \newcommand{\BC}{\mathcal{BC}}

    \newcommand{\brs}{B^{r,s}}		
    \newcommand{\chsig}{\check{\sigma}}
    \newcommand{\muz}{\mu_{\hat{0}}}
    \newcommand{\muo}{\mu_{\hat{1}}}
    \newcommand{\sig}{\sigma}
    \newcommand{\tbrs}{\tilde{B}^{r,s}}

    \newcommand{\bcdual}{*_{\mathcal{BC}}}    
    \newcommand{\bks}[1]{B^{#1,s}}
    \newcommand{\cl}{\mathrm{cl}}
    \newcommand{\cw}{\mathrm{cw}}
    \newcommand{\dual}{*}
    \newcommand{\es}{\varnothing}
    \newcommand{\etil}[1]{e_{#1}}
    \newcommand{\ftil}[1]{f_{#1}}
    \newcommand{\hsig}{\check{\sigma}}
    \newcommand{\inner}[2]{\langle #1\,,\,#2\rangle}
    \newcommand{\iotil}{\check{\iota}}    
    \newcommand{\La}{\Lambda}
    \newcommand{\Le}{\hbox{\rotatedown{$\Gamma$}}}
    \newcommand{\Lee}{\hbox{\rotatedown{L}}}
    \newcommand{\Lt}{\om_{2}}
    \newcommand{\la}{\lambda}
    \newcommand{\lan}{\langle}
    \newcommand\lev{\mbox{\sl lev}\,}
    \newcommand{\Lz}{\Lambda_{0}}
    \newcommand{\Lo}{\Lambda_{1}}
    \newcommand{\Mmin}[1]{\mathcal{M}_{\min}(#1)}
    \newcommand{\om}{\varpi}
    
    \newcommand{\qbin}[2]{\genfrac{[}{]}{0pt}{}{#1}{#2}}
    \newcommand{\qq}{\mathbb{Q}(q)}    
    \newcommand{\ran}{\rangle}
    \newcommand{\stack}[2]{\genfrac{}{}{0pt}{}{#1}{#2}}
    \newcommand{\Sup}{\mathcal{BC}}
    \newcommand{\Tred}{T^{\#}}
    
    \newcommand{\uq}[1]{U_{q}(\mathfrak{#1})}
    \newcommand{\uqp}[1]{U_{q}'(\mathfrak{#1})}
    
    \newcommand{\uqsophat}[1]{U_{q}'(\widehat{\mathfrak{so}}_{#1})}
    
    \newcommand{\vare}{\varepsilon}
    \newcommand{\wt}{\mathrm{wt}}
    \newcommand{\T}{\mathcal{M}}
    \newcommand{\tbt}{\tilde{B}^{2,2}}
    \newcommand{\tbs}{\tilde{B}^{2,s}}
    \newcommand{\Z}{\mathbb{Z}}

\begin{document}
    
        \pagenumbering{roman}
    \pagestyle{plain}

    %
    %

    \singlespacing

    \begin{center}

        \begin{huge}
	    Applications of Crystal Bases to Current Problems in
	    Representation Theory
	\end{huge}\\\n
        By\\\n
        {\sc Philip Max Sternberg}\\
        B.A. (University of California, Berkeley) 2002\\\n
        DISSERTATION\\\n
        Submitted in partial satisfaction of the requirements for the degree of\\\n
        DOCTOR OF PHILOSOPHY\\\n
        in\\\n
        MATHEMATICS\\\n
        in the\\\n
        OFFICE OF GRADUATE STUDIES\\\n
        of the\\\n
        UNIVERSITY OF CALIFORNIA\\\n
        DAVIS\\\n\n
        Approved:\\\n\n
        
        \rule{4in}{1pt}\\
        Anne Schilling (Chair)\\\n\n
        
        \rule{4in}{1pt}\\
        Eric M. Rains\\\n\n
        
        \rule{4in}{1pt}\\
        Monica Vazirani\\\n\n
        
        Committee in Charge\\[8pt]\n\n

    \end{center}

    \newpage
    
    \begin{center}
	\raisebox{-4in}{\copyright\ Philip Sternberg, 2006.  All
	rights reserved.}
    \end{center}
    
    \thispagestyle{empty}
    \addtocounter{page}{-1}

    \newpage

    \begin{center}
	\raisebox{-4in}{Dedicated to the memory of David
	Sternberg,}\\
	\raisebox{-16pt}{father, scientist, and mentor.}
    \end{center}

    \newpage

    \doublespacing

    %
    %
    
    \tableofcontents
    
    \newpage

    %
    %
    

    \begin{flushright}
        \singlespacing
        Philip Sternberg\\
        June 2006\\
        Mathematics
    \end{flushright}
    
    \centerline{\textbf{\underline{Abstract}}}
    
        Over the past sixteen years, Kashiwara's construction of crystal bases
\cite{Kash:1990, Kash:1991} has been an invaluable tool in addressing
questions in representation theory, statistical mechanics, and soliton
cellular automata.  This dissertation addresses two current questions
in the theory of crystal bases: the structure of crystal bases for
Kirillov-Reshetikhin modules over algebras of type $D_{n}^{(1)}$, and
the local structure of doubly laced crystals.  It also includes a user
manual for the CrystalView software package, which was written by the
author.

The first of these questions is addressed in Chapters \ref{sec:KR} and
\ref{sec:B2s}.  Specifically, in Chapter \ref{sec:KR} we describe a
conjecture for the structure of the crystal bases $B^{r,s}$ that
correspond to the Kirillov-Reshetikhin modules $W^{r,s}$ of type
$D_{n}^{(1)}$.  Specifically, we prove that many aspects of the
combinatorial structure of these crystals are determined by
representation theoretical concerns, leaving only a small aspect of
this structure unspecified.  Our conjecture specifies this part of the
structure of these crystals uniformly.  This result is based on the
work and conjectures of the $X=M$ program put forward in
\cite{HKOTT:2001, HKOTY:1999}.  In Chapter \ref{sec:B2s} we prove that
this conjecture is true for the special case of $r=2$.

The second of these questions is the subject of Chapter
\ref{sec:local}.  Building on the work of \cite{Stem}, we confirm a
conjecture of Stembridge that limits the local behavior of crystals
for representations of Lie algebras all of whose rank two subalgebras
are either $A_{1}\times A_{1}$, $A_{2}$, or $C_{2}$.  Specifically, we
show that if $v$ is a vertex in such a crystal with both an
$i$-colored and a $j$-colored edge directed into it, one of four pairs
of sequences $(E_{1},E_{2})$, each composed of the operators $e_{i}$
and $e_{j}$, satisfies $E_{1}v = E_{2}v$.

In Chapter \ref{sec:CV}, we explain the features, usage, and
implementation of CrystalView.

    \newpage

    %
    %

    \chapter*{Acknowledgments}

        I first wish to thank Anne Schilling for her excellent supervision of
the work contained in this dissertation.  The training I received
under her guidance both as a creator and communicator of mathematics
has been truly exceptional.  The results contained here could not have
been obtained without Prof.  Schilling's help as critic and
collaborator.

My thanks are also extended to the other members of my dissertation
committee, Eric Rains and Monica Vazirani, who also served on my
qualifying examination committee.  They have both been a great help to
me in many conversations as I have explored many diverse mathematical
interests.

I could not have completed my degree without the support of my wife,
Mary Orland.  The obstacles to completing a Ph.D. are not solely
intellectual--Mary helped me to meet many of these other challenges.
My gratitude is due to her for keeping me out of the trouble in which
I could have very easily found myself without her assistance.

It is unlikely that I would have ever considered pursuing an advanced
degree had I not had as a role model my father, David Sternberg.  From
my earliest memories, my father encouraged me to pursue new knowledge
and understanding in every way possible.  I am proud to be able to
dedicate this dissertation to the honor of his memory.  I was also
always supported by my mother, Linda Sternberg, throughout my studies.
I am forever grateful to both of my parents for their love and
support, and their belief that my capability was limited only by my
ambition.

Many thanks are also extended to Isaiah Lankham, who has been a
splendid office-mate, travelling companion, and dear friend for the
past four years.  In particular, Isaiah's editorial improvements on
many of my papers and his suggestions pertaining to CrystalView were
indispensible.

The result of Chapter \ref{sec:local} could not have been obtained
without John Stembridge, to whom I also extend my thanks; not only is
he responsible for the seminal paper which it generalizes, but he also
provided very helpful insights into the mechanics of this problem.  My
thanks are also given to Masato Okado, Tomoki Nakanishi, and Mark
Shimozono for their very helpful comments regarding the work of
Chapters \ref{sec:KR} and \ref{sec:B2s}.

I would also like to thank Jes\'us de Loera and Greg Kuperberg for
serving on my qualifying examination committee, and for many helpful
and enjoyable conversations over the past years.  My thanks are
particularly due to Prof.  de Loera for his suggestion that I
investigate working with Prof.  Schilling, and to Prof.  Kuperberg for
agreeing to serve on my examination committee on very short notice
after one of my original committee members found himself suddenly
unavailable.

In the course of publishing and presenting the research contained
here, several anonomous referees have provided comments and
criticisms which have ultimately improved the exposition of these
results considerably.  I wish to thank these referees for volunteering
their time for the betterment of the community.

Parts of this research were carried out while I was a guest at two
excellent institutions, the Max-Plank-Institut f\"ur Mathematik in
Bonn, Germany, and the Research Institute for the Mathematical
Sciences in Kyoto, Japan.  I extend my thanks to these organizations
for their hospitality.

This work was supported in part by NSF grants DMS-0135345,
DMS-0200774, and DMS-0501101.

    \newpage

    %
    %

    \pagestyle{fancy}
    \pagenumbering{arabic}
       
    %
    %

    \newchap{Introduction}
    \label{sec:Intro}
    
        \section{Historical Overview and Motivation}
        \label{sec:Intro.Motivation}
        
            \subsection{Representation Theory}
            \label{sec:Intro:Motivation:RT}
            The study of representation theory has as its goal the understanding
of algebraic objects in terms of their action on other objects.  For
instance, we may study groups by their action on a set of elements (a
permutation representation) or on a vector space (a linear
representation).  This notion is formally described by realizing this
action as a homomorphism from the object being studied into the the
most general set of symmetries of the object being acted upon.  In the
above examples, these would be a homomorphism from a group into a
group of permutations or into the general linear group of a vector
space, respectively.

If we wish to study the representation theory of an algebra $A$ over a
field $\mathbb{K}$, we consider morphisms $\rho:A\to\mathfrak{gl}(V)$,
where $V$ is a vector space over $\mathbb{K}$.  If $A$ is a Lie
algebra, we take $\mathfrak{gl}(V)$ to be the general linear Lie
algebra of transformations of $V$; if $A$ is instead an associative
algebra, we take as $\mathfrak{gl}(V)$ the associative general linear
algebra of $V$.

In the following work, we are concerned with representations of
Quantum Groups, also known as Quantized Universal Enveloping Algebras.
These are associative algebras over the field $\mathbb{C}(q)$ of
rational functions with complex coefficients.  Denoted by
$U_{q}(\mathfrak{g})$, they can be realized as a quantization
deformation of the universal enveloping algebra $U(\mathfrak{g})$ of a
Lie algebra $\mathfrak{g}$.  When $\mathfrak{g}$ is a Kac-Moody
algebra, there is a straightforward presentation of
$U_{q}(\mathfrak{g})$ by generators and relations.

Specifically, we have the following, after \cite{HK}.

\begin{definition}
    \label{def:qgroup}
    Let $I$ be a finite indexing set; let $A = (a_{ij}):i,j\in 
    I$ be a symmetrizable generalized Cartan matrix with symmetrizing 
    matrix $D = \mathrm{diag(s_{i})}:i\in I$, and let 
    $\Pi,\Pi^{\vee},P,$ and $P^{\vee}$ be the set of simple roots, 
    simple coroots, weight lattice, and dual weight lattice of the 
    corresponding Kac-Moody Lie algebra $\mathfrak{g}$.  The 
    associated quantum group $U_{q}(\mathfrak{g})$ is the associative 
    $\mathbb{C}(q)$-algebra with unit generated by $e_{i}$ and 
    $f_{i}$ for $i\in I$ and $q^{h}$ for $h\in P^{\vee}$ with the 
    following relations:
    \begin{itemize}
        \item $q^{0} = 1$, $q^{h}q^{h'} = q^{h+h'}$ for $h,h'\in 
	P^{\vee}$,
    
        \item $q^{h}e_{i}q^{-h} = q^{\alpha_{i}(h)}e_{i}$ for $h\in 
	P^{\vee}$,
    
        \item $q^{h}f_{i}q^{-h} = q^{-\alpha_{i}(h)}f_{i}$ for $h\in 
	P^{\vee}$,
    
        \item $e_{i}f_{j} - f_{j}e_{i} = 
	\delta_{ij}\frac{K_{i}-K_{i}^{-1}}{q_{i}-q_{i}^{-1}}$ for 
	$i,j\in I$,
    
        \item 
	$\sum_{k=0}^{1-a_{ij}}(-1)^{k}\big[\stack{1-a_{ij}}{k}\big]_{q_{i}}e_{i}^{1-a_{ij}-k}e_{j}e_{i}^{k}
	= 0$ for $i\neq j$,
    
        \item 	
	$\sum_{k=0}^{1-a_{ij}}(-1)^{k}\big[\stack{1-a_{ij}}{k}\big]_{q_{i}}f_{i}^{1-a_{ij}-k}f_{j}f_{i}^{k}
	= 0$ for $i\neq j$,

    \end{itemize}
    where $q_{i} = s_{i}$ and $K_{i} = q^{s_{i}h_{i}}$, and the symbol
    $[n]_{q} = \frac{q^{n}-q^{-n}}{q-q^{-1}}$ for
    $n\in\mathbb{Z}_{\geq0}$.
\end{definition}

The algebra $U_{q}(\mathfrak{g})$ has a representation theory that
refines the representation theory of the associated Lie algebras by
means of quantization deformation.  More precisely, given a simple
highest weight module $M$ over $\mathfrak{g}$, we can construct a
simple module $M_{q}$ over $U_{q}(\mathfrak{g})$ whose weight spaces
multiplicities are the same as those of $M$, as proved by Lusztig in
\cite{Lus:1988}.  Furthermore, tensor products of such deformed
modules decompose into simple modules in the same way as their
classical counterparts.  Thus, the representation theory of
$U_{q}(\mathfrak{g})$ retains all algebraic information about
$\mathfrak{g}$, in addition to it modules having a natural grading by
powers of $q$.

            \subsection{Crystal Bases}
            \label{sec:Intro:Motivation:CB}
            The representation theory of Lie algebras has always had close ties to
combinatorics.  Using the technique of crystal bases, Kashiwara
developed a completely discrete language for describing modules over
quantum groups, and thus also modules over Lie algebras
\cite{Kash:1990, Kash:1991}.  Specifically, given any highest weight
integrable $U_{q}(\mathfrak{g})$-module $V$, there is a colored
directed graph $B(V)$ called the crystal graph of $V$ that encodes
nearly all of the information of the action of $U_{q}(\mathfrak{g})$
on $V$.  Abstracting this notion, one can define crystals without
making an explicit reference to a particular representation.  This
definition includes all crystals that come from
$\mathfrak{g}$-modules, as well as crystals that cannot be thus
realized.  Such axiomatic crystals have been used to prove several
significant results in the theory of crystal bases.

Precisely, we have the following definition of crystals, after
\cite{Kash:1993}.

\begin{definition}
    \label{def:crystal}
    Let $I$ be a finite index set, let $A=(a_{ij})$ for $i,j\in I$ be
    a Cartan matrix, and let $(\Pi,\Pi^{\vee},P,P^{\vee})$ be the
    associated Cartan data of a Kac-Moody algebra $\mathfrak{g}$ as in
    Definition \ref{def:qgroup}.  A $\mathfrak{g}$-crystal is a set
    $B$ together with maps $\wt:B\to P, e_{i},f_{i}:B\to B\cup\{0\}$,
    and $\varepsilon_{i},\varphi_{i}:B\to\mathbb{Z}\cup\{-\infty\}$
    for $i\in I$ such that:
    \begin{itemize}
        \item $\varphi_{i}(b) = \varepsilon_{i}(b) + \lan 
	h_{i},\wt(b)\ran$ for all $i\in I$,
    
        \item $\wt(e_{i}b) = \wt b + \alpha_{i}$ if $e_{i}b \in B$,
    
        \item $\wt(f_{i}b) = \wt b - \alpha_{i}$ if $f_{i}b \in B$,
    
        \item $\varepsilon_{i}(e_{i}b) = \varepsilon_{i}(b) - 1, 
	\varphi_{i}(e_{i}b) = \varphi_{i}(b) + 1$ if $e_{i}b\in B$,
    
        \item $\varepsilon_{i}(f_{i}b) = \varepsilon_{i}(b) + 1, 
	\varphi_{i}(f_{i}b) = \varphi_{i}(b) - 1$ if $f_{i}b\in B$,
    
        \item $f_{i}b = b'$ if and only if $b= e_{i}b$ for $b,b'\in 
	B, i\in I$,
    
        \item if $\varphi_{i}(b) = -\infty$ for $b\in B$, then 
	$e_{i}b = f_{i}b = 0$.
    \end{itemize}
    
\end{definition}

The corresponding crystal graph has a vertex for each $b\in B$ and an
$i$-colored edge from $b$ to $b'$ if $f_{i}b = b'$.

            \subsection{Quantum Affine Algebras
	    }
            \label{sec:Intro:Motivation:QAA}
            Many current questions involving crystal bases pertain to crystals for
modules of quantum affine algebras; i.e., quantum groups
$U_{q}(\mathfrak{g})$ where the underlying algebra $\mathfrak{g}$ is
of affine type.  Two important classes of representations of these
algebras are finite dimensional modules and highest weight modules.
Many existing methods for analyzing representations of Lie algebras
and quantum groups, both combinatorial and algebraic, only apply to
highest weight modules.  A great deal of new work was therefore needed
to address questions about the finite-dimensional representation
theory of quantum affine algebras.

Finite-dimensional modules of quantum affine algebras were classified
by Chari and Pressley using Drinfeld polynomials \cite{CP:1995},
taking a major step toward understanding these representation
theoretical concerns.  However, many of the internal aspects of these
representations remain poorly understood.  One approach to
understanding these representations more explicitly is by determining
if they have crystal bases, and if so, describing their structure via
an explicit discrete construction.

In \cite{KKR:1988, KR:1988}, a relationship was observed between the
Bethe ansatz (a method for solving lattice models in statistical
mechanics) and the combinatorics of Young tableaux.  There it was
proposed that a certain family of modules over quantum affine algebras
of type $A_{n}^{(1)}$, which have come to be called
Kirillov-Reshetikhin modules, should be combinatorializable (as these
papers precede \cite{Kash:1990}, we do not use the term
``crystallizable'').  This claim is motivated by the fact that another
method for solving the same lattice models, the corner transfer
matrix, makes use of the combinatorics that govern type $A$
representation theory.  The Bethe ansatz gives rise to another set of
combinatorial objects, called rigged configurations; since the formulas
arising from these two methods are equal, there should be a bijection
between rigged configurations and Young tableaux.  As this notion was
formalized over the next decade, the combinatorial notion of crystals,
type-independent tableaux, and type-independent rigged configurations
allowed the work of \cite{KKR:1988, KR:1988} to be generalized.

In \cite{HKOTT:2001, HKOTY:1999}, the Kyoto school proposed that all
Kirillov-Reshetikhin modules of all types have crystal bases, and
furthermore that the only finite-dimensional modules of quantum affine
algebras with crystal bases are tensor products of
Kirillov-Reshetikhin modules.  Considerable progress has been made
toward understanding these crystals and proving these claims.  For a
thorough review of the current state of this program, see
\cite{S:2005}.  Here, we summarize the known results for
Kirillov-Reshetikhin crystals \cite{O}.

For type $A_{n}^{(1)}$, the crystals are known to exist
\cite{KKMMNN:1992a} and their explicit structure is known
\cite{Sh:2002}.  For non-exceptional algebras, \cite{KKM} gives a
realization of $B^{1,s}$ as the set of lattice points in a convex
polytope.  The crystal $B^{r,1}$ for all algebras was shown to exist
in \cite{Kash:2002}, and a realization using Lakshmibai-Seshadri (L-S)
paths has been constructed in \cite{NS}.  Other cases include those of
the spin represantations $B^{n,s}$ for $C_{n}^{(1)}$, $B^{n,s}$ and
$B^{n-1,s}$ for $D_{n}^{(1)}$, and $B^{n,s}$ for $D_{n+1}^{(2)}$,
which are shown to exist and given explicitly in \cite{KKMMNN:1992a}.
In \cite{BFKL} any case in which the highest classical component of
$B^{r,s}$ is the adjoint representation of the classical subalgebra is
dealt with.  The crystal structures have also been given for
$B^{r,1+\delta_{rn}}$ for $B_{n}^{(1)}$, $B^{r,2}$ with $r\neq n$ for
$C_{n}^{(1)}$, $B^{r,1}$ with $r\neq n$ and $r\neq n-1$ for
$D_{n}^{(1)}$ \cite{Koga:1999}, $B^{r,1}$ for $C_{n}^{(1)}$, $B^{r,1}$
for $A_{2n}^{(2)}$, $B^{r,1}$ for $A_{2n-1}^{(2)\dagger}$ \cite{JMO},
and $B^{1,s}$ for $G_{2}^{(1)}$ \cite{Ya}.

The results of \cite{OSS:2003,OSS:2003a} reduce the problem to simply
laced types using virtual crystals, a technique whereby Dynkin diagram
folding is extended to crystals.  Therefore, these conjectures need
only to be resolved for types $D_{n}^{(1)}$ ($n\geq 4$) and
$E_{n}^{(1)}$ ($n\in\{6,7,8\}$).  Furthermore, it has been shown in
\cite{FSS} that given certain assumptions about the crystals
$B^{r,s}$, their affine structure is uniquely specified by an
isomorphism with Demazure crystals.  Thus, defining a crystal
structure satisfying these assumptions constitutes a significant step
towards understanding their combinatorial structure.

	    \subsection{Local Characterization of Crystals}
	    \label{sec:Intro:local}
	    Crystals were originally developed to describe representations of Lie
algebras; however, many objects satisfying Definition
\ref{def:crystal} do not correspond to any
$U_{q}(\mathfrak{g})$-module.  It is therefore of interest to know
when a colored directed graph is the crystal graph of a
representation.  Such a characterization can be used to show that a
graph that does not {\it a priori} encode representation theoretic
information is in fact a crystal graph.

Many explicit combinatorial models have been developed for crystal
graphs of representations; two of these are paths in the weight
lattice of $\mathfrak{g}$ \cite{L1,L2} and generalized Young tableaux
\cite{KN}.  In all of these cases, the combinatorial properties of the
crystals are defined globally.  In \cite{Stem}, Stembridge introduced
a set of graph theoretic axioms, each of which addresses only local
properties of a colored directed graph, that characterizes highest
weight crystal graphs that come from representations of simply laced
algebras.

        \section{Summary of Main Results}
        \label{sec:Intro:Summary}
        The first main result presented here is a conjectural explicit
description of crystals for Kirillov-Reshetikhin modules over an
algebra of type $D_{n}^{(1)}$.  We investigate the structure of such
modules as modules over the embedded algebras of type $D_{n}$ and
$D_{n-1}$, discovering some very significant symmetries.  Based on a
refinement of the branching properties of $B^{r,s}$ over these
algebras, we formulate a combinatorial conjecture for the structure of
these crystals.

Specifically, given $r\in \{1,2,\ldots,n\}$ and $s\in\mathbb{Z}_{\geq
0}$, we construct a crystal $\tilde{B}^{r,s}$.

\begin{conjecture}
    \label{conj:brs}
    Assuming that $B^{r,s}$ exists with the properties conjectured in
    \cite{HKOTT:2001, HKOTY:1999}, we have
    $$
    B^{r,s}\simeq\tilde{B}^{r,s}.
    $$
\end{conjecture}

The work of \cite{FSS} shows that if we assume $B^{r,s}$ to have
certain properties, its affine structure is specified by the inclusion
of a Demazure module in it.  Thus, if these hypotheses (which are less
demanding than those of \cite{HKOTT:2001, HKOTY:1999}) are satisfied,
this conjecture is true.

As evidence of this conjecture, we prove in Section \ref{sec:B2s:r2}
that in the case of $r = 2$, the affine structure of the
Kirillov-Reshetikhin modules respects the observed symmetries.  We
thus have the following Theorem, as proved in \cite{SS:2006}.

\begin{theorem}
    \label{thm:b2s}
    Assuming that $B^{2,s}$ exists with the properties conjectured in
    \cite{HKOTT:2001, HKOTY:1999}, we have
    $$
    B^{2,s}\simeq\tilde{B}^{2,s}.
    $$    
\end{theorem}

Further elucidating the stucture of $B^{2,s}$ for type $D_{n}^{(1)}$, 
we define in Section \ref{sec:KR:tabs} a new generalization of Young 
tableaux and show that it can be used to index the vertices of 
$B^{2,s}$.  We also give an explicit definition of the affine 
Kashiwara operators on these tableaux.  We encapsulate these results 
as follows.

\begin{theorem}
    \label{thm:b2stabs}
    The tableaux defined in Definition \ref{def:affine_tab} are in
    bijective correspondence with the tableaux used to index the
    classical components of $B^{2,s}$.  Furthermore, Algorithm
    \ref{alg:iota} explicitly specifies the action of the Kashiwara
    operators $e_{0}$ and $f_{0}$ on these tableaux.
\end{theorem}

Finally, we confirm a conjecture of Stembridge on the local structure
of crystals coming from representations of doubly laced algebras
\cite{Ster}.

\begin{theorem}
    \label{proposition:int}
    Let $\mathfrak{g}$ be a doubly laced algebra, i.e., an algebra
    all of its regular rank 2 subalgebras are of type $A_{1}\times
    A_{1}$, $A_{2}$, or $C_{2}$.  Let $B$ be the crystal graph of an
    irreducible highest weight module of $\mathfrak{g}$, and let $v$
    be a vertex of $B$ such that $e_{i}v \neq 0$ and $e_{j}v \neq 0$,
    where $e_{i}$ and $e_{j}$ denote two different Kashiwara raising
    operators.  Then one of the following is true:
    \begin{enumerate}
        \item $e_{i}e_{j}v = e_{j}e_{i}v$,
    
	\item $e_{i}e_{j}^{2}e_{i}v = e_{j}e_{i}^{2}e_{j}v$ and no
	other sequences of the operators $e_{i}, e_{j}$ with
	length less than or equal to four satisfy such an equality,
    
	\item $e_{i}e_{j}^{3}e_{i}v = e_{j}e_{i}e_{j}e_{i}e_{j}v =
	e_{j}^{2}e_{i}^{2}e_{j}v$ and no other sequences of the
	operators $e_{i}, e_{j}$ with length less than or equal to
	five satisfy such an equality,
    
	\item $e_{i}e_{j}^{3}e_{i}^{2}e_{j}v =
	e_{i}e_{j}^{2}e_{i}e_{j}e_{i}e_{j}v =
	e_{j}e_{i}^{2}e_{j}^{3}e_{i}v = 
	e_{j}e_{i}e_{j}e_{i}e_{j}^{2}e_{i}v$ and 
	no other sequences of the operators $e_{i}, e_{j}$ with length
	less than or equal to seven satisfy such an equality.
    \end{enumerate}
     
    The equivalent statement with $f_{i}$ and $f_{j}$ in place of
    $e_{i}$ and $e_{j}$ also holds.

\end{theorem}

\subsection{Outline by chapter}

In Chapter \ref{sec:KR}, we include in the first section all the
necessary prerequisite constructions and well-known theorems on
quantum groups and crystal bases, including perfect crystals.  In
section \ref{sec:KR:tabs} we go into more detail on type $D$ crystals,
specifying the combinatorics of representations of even orthogonal
algebras and specializing the conjectures on finite-dimensional
modules over quantum affine algebras to the type $D_{n}^{(1)}$ case.
Sections \ref{sec:KR:BC} and \ref{sec:KR:sigma} contain a detailed
analysis of the affine crystals that are supposed to exist by these
conjectures, and we present our conjecture on the exact structure of
these crystals.

In Chapter \ref{sec:B2s}, we begin by specializing the work of the
last two sections to the case of $r=2$ in section \ref{sec:B2s:tabs}.
We extend these results in section \ref{sec:B2s:tabs2} by constructing
an explicit bijection between classical tableaux and what we define to
be ``affine tableaux''.  In section \ref{sec:B2s:r2} we prove Theorem 
\ref{thm:b2s}.  This chapter is based on \cite{SS:2006}.

In Chapter \ref{sec:local} we deal with local properties of highest
weight crystals coming from representations of some classes of Lie
algebras.  In section \ref{sec:local:simply}, we recall the previously
known case of simply laced crystals.  Section \ref{sec:local:B} is a
review of the combinatorics of type $C_{2}$ crystals.  Proceding to
section \ref{sec:local:genB}, we analyze a generic $C_{2}$ crystal to 
see what its neighborhood in the crystal graph may look like.  In
section \ref{sec:local:proof}, we prove that only four possible local 
behaviors may appear in doubly laced crystals.  This chapter is based 
on \cite{Ster}.

Finally, Chapter \ref{sec:CV} deals with the software package
CrystalView, detailing requirements, user input options, current
limitations of the software, and a description of the algorithms used 
to generate the output of this program.

    %
    %

    \newchap{Crystals for Kirillov-Reshetikhin Modules over Quantum 
    Affine Algebras of Type $D_{n}^{(1)}$}
    \label{sec:KR}
        
        \section{Review of Quantum Groups and Crystal Bases}
        \label{sec:KR:Review}
        \subsection{Quantum groups}
For $n\in\mathbb{Z}$ and a formal parameter $q$, we use the notation
\begin{displaymath}
    [n]_{q}=\frac{q^{n}-q^{-n}}{q-q^{-1}},\quad
    [n]_{q}!=\prod_{k=1}^{n}[k]_{q},\textrm{ and }
    \qbin{m}{n}_q=\frac{[m]_{q}!}{[n]_{q}![m-n]_{q}!}.
\end{displaymath}
These are all elements of $\qq$, called the $q$-integers, 
$q$-factorials, and $q$-binomial coefficients, respectively.

Let $\mathfrak{g}$ be an arbitrary Kac-Moody Lie algebra with Cartan
datum
$(A,\Pi,\Pi^{\vee},P,P^{\vee})$ and a Dynkin diagram indexed by $I$.
Here $A=(a_{ij})_{i,j\in I}$ is the Cartan matrix, $P$ and $P^\vee$
are the weight lattice and dual weight lattice, respectively,
$\Pi=\{\alpha_i \mid i\in I\}$ is the set of simple roots and
$\Pi^\vee=\{h_i \mid i\in I\}$ is the set of simple coroots.
Furthermore, let $\{s_i \mid i\in I\}$ be the entries of the diagonal 
symmetrizing matrix of $A$ and define $q_{i}=q^{s_{i}}$ and
$K_{i}=q^{s_{i}h_{i}}$.
Then the quantum enveloping algebra $\uq{g}$ is the associative
$\qq$-algebra 
generated by $e_{i}$ and $f_{i}$ for $i\in I$, and $q^{h}$ for $h\in
P^{\vee}$, 
with the following relations (see e.g. \cite[Def. 3.1.1]{HK}):
\begin{enumerate}
    \item $q^{0}=1$, $q^{h}q^{h'}=q^{h+h'}$ for all $h,h'\in
P^{\vee}$,
    \item $q^{h}e_{i}q^{-h}=q^{\alpha_{i}(h)}e_{i}$ for all $h\in
P^{\vee}$,
    \item $q^{h}f_{i}q^{-h}=q^{\alpha_{i}(h)}f_{i}$ for all $h\in
P^{\vee}$,
    \item
$e_{i}f_{j}-f_{j}e_{i}=\delta_{ij}\frac{K_{i}-K_{i}^{-1}}{q_{i}-q_{i}^{-1}}$

          for $i,j\in I$,
    \item
$\sum_{k=0}^{1-a_{ij}}(-1)^{k}\qbin{1-a_{ij}}{k}_{q_{i}}e_{i}^{1-a_{ij}-k}
          e_{j}e_{i}^{k}=0$ for all $i\neq j$,
    \item
$\sum_{k=0}^{1-a_{ij}}(-1)^{k}\qbin{1-a_{ij}}{k}_{q_{i}}f_{i}^{1-a_{ij}-k}
          f_{j}f_{i}^{k}=0$ for all $i\neq j$.
\end{enumerate}

\subsection{Crystal bases}\label{sec:crystal bases}
The quantum algebra $\uq{g}$ can be viewed as a $q$-deformation of the
universal enveloping algebra $U(\mathfrak{g})$ of $\mathfrak{g}$.
Lusztig~\cite{Lus:1988} showed that the integrable highest weight
representations of $U(\mathfrak{g})$ can be deformed to
$U_q(\mathfrak{g})$ representations in such a way that the dimension
of the weight spaces are invariant under the deformation, provided
$q\neq0$ and $q^{k}\neq1$ for all $k\in\mathbb{Z}$ (see
also~\cite{HK}).  Let $M$ be a $U_q(\mathfrak{g})$-module and $R$ the
subset of all elements in $\qq$ which are regular at $q=0$.
Kashiwara~\cite{Kash:1990,Kash:1991} introduced Kashiwara operators
$\tilde{e}_{i}$ and $\tilde{f}_{i}$ as certain linear combinations of
powers of $e_i$ and $f_i$.  In the sequel we will use the notation
$e_{i}$ and $f_{i}$ to refer to the Kashiwara operators, rather than
the standard generators of $U_{q}(\mathfrak{g})$.  A crystal lattice
$\mathcal{L}$ is a free $R$-submodule of $M$ that generates $M$ over
$\qq$, has the same weight decomposition and has the property that
$\etil{i}\mathcal{L}\subset \mathcal{L}$ and
$\ftil{i}\mathcal{L}\subset \mathcal{L}$ for all $i\in I$.  The
passage from $\mathcal{L}$ to the quotient $\mathcal{L}/q\mathcal{L}$
is referred to as taking the crystal limit.  A crystal basis (also
called a crystal base by some authors) is a $\mathbb{Q}$-basis of
$\mathcal{L}/q\mathcal{L}$ with certain properties.

We can think of a $U_q(\mathfrak{g})$-crystal as a nonempty set $B$
equipped with maps $\wt:B\rightarrow P$ and
$\etil{i},\ftil{i}:B\rightarrow B\cup\{\es\}$ for all $i\in I$,
satisfying
\begin{align}
\label{eq:e-f}
\ftil{i}(b)=b' &\Leftrightarrow \etil{i}(b')=b
\text{ if $b,b'\in B$} \\
\wt(\ftil{i}(b))&=\wt(b)-\alpha_i \text{ if $\ftil{i}(b)\in B$} \\
\label{eq:string length}
\inner{h_i}{\wt(b)}&=\varphi_i(b)-\varepsilon_i(b).
\end{align}
In the most general setting, $\varepsilon_{i}$ and $\varphi_{i}$ can
be integer valued functions on $B$ satisfying further axioms; however,
here we consider only normal crystals, in which case we simply define
for all $b \in B$
\begin{equation*}
\begin{split}
\varepsilon_i(b)&= \max\{n\ge0\mid \etil{i}^n(b)\not=\es \} \\
\varphi_i(b) &= \max\{n\ge0\mid \ftil{i}^n(b)\not=\es \}.
\end{split}
\end{equation*}
Here we assume that $\varphi_i(b),\varepsilon_i(b)<\infty$ for all
$i\in I$ and $b\in B$, since we restrict our attention to crystals
that come from integrable $U_{q}(\mathfrak{g})$-modules.  A
$U_q(\mathfrak{g})$-crystal $B$ can be viewed as a directed
edge-colored graph (the crystal graph) whose vertices are the elements
of $B$, with a directed edge from $b$ to $b'$ labeled $i\in I$ if and
only if $\ftil{i}(b)=b'$.

Let $B_1$ and $B_2$ be $U_q(\mathfrak{g})$-crystals. The Cartesian
product
$B_2\times B_1$ can also be endowed with the structure of a
$U_q(\mathfrak{g})$-crystal.
The resulting crystal is denoted by $B_2\otimes B_1$ and its elements
$(b_2,b_1)$ are 
written $b_2\otimes b_1$.
(The reader is warned that in the interests of compatibility with the 
combinatorics of Young tableaux, our convention is opposite to that of 
Kashiwara~\cite{Kash:1995}). For $i\in I$ and $b=b_2\otimes b_1\in
B_2\otimes B_1$,
we have $\wt(b)=\wt(b_1)+\wt(b_2)$,
\begin{equation} \label{eq:f on two factors}
\ftil{i}(b_2\otimes b_1) = \begin{cases} \ftil{i}(b_2)\otimes b_1
& \text{if $\varepsilon_i(b_2)\ge \varphi_i(b_1)$} \\
b_2\otimes \ftil{i}(b_1) & \text{if $\varepsilon_i(b_2)<\varphi_i(b_1)$}
\end{cases}
\end{equation}
and
\begin{equation} \label{eq:e on two factors}
\etil{i}(b_2\otimes b_1) = \begin{cases} \etil{i}(b_2) \otimes b_1 &
\text{if $\varepsilon_i(b_2)>\varphi_i(b_1)$} \\
b_2\otimes \etil{i}(b_1) & \text{if $\varepsilon_i(b_2)\le
\varphi_i(b_1)$.}
\end{cases}
\end{equation}
Combinatorially, this action of $\ftil{i}$ and $\etil{i}$ on tensor
products can be described by the signature rule.
The $i$-signature of $b$ is the word consisting of the symbols $+$
and $-$ 
given by
\begin{equation*}
\underset{\text{$\varphi_i(b_2)$ times}}{\underbrace{-\dotsm-}}
\quad \underset{\text{$\varepsilon_i(b_2)$
times}}{\underbrace{+\dotsm+}} \quad 
\underset{\text{$\varphi_i(b_1)$ times}}{\underbrace{-\dotsm-}}
\quad \underset{\text{$\varepsilon_i(b_1)$
times}}{\underbrace{+\dotsm+}} .
\end{equation*}
The reduced $i$-signature of $b$ is the subword of the
$i$-signature of $b$, given by the repeated removal of adjacent
symbols $+-$ (in that order); it has the form
\begin{equation*}
\underset{\text{$\varphi$ times}}{\underbrace{-\dotsm-}} \quad
\underset{\text{$\varepsilon$ times}}{\underbrace{+\dotsm+}}.
\end{equation*}
If $\varphi=0$ then $\ftil{i}(b)=\es$; otherwise $\ftil{i}$ acts
on the tensor factor corresponding to the rightmost symbol $-$ in the 
reduced $i$-signature of $b$. Similarly, if $\varepsilon=0$ then
$\etil{i}(b)=\es$; otherwise $\etil{i}$ acts on the leftmost symbol
$+$ 
in the reduced $i$-signature of $b$. From this it is clear that
\begin{equation*}
\begin{split}
&\varphi_{i}(b_{2}\otimes 
b_{1})=\varphi_{i}(b_{2})+\mathrm{max}(0,\varphi_{i}(b_{1})-\vare_{i}(b_{2}))
,\\
&\vare_{i}(b_{2}\otimes
b_{1})=\vare_{i}(b_{1})+\mathrm{max}(0,-\varphi_{i}(b_{1})+\vare_{i}(b_{2})).
\end{split}
\end{equation*}

For any dominant weight $\la$, there is an irreducible highest weight 
module $V(\la)$ with highest weight $\la$.  Provided that $V(\la)$ is 
integrable, it has a crystal basis denoted by $B(\la)$.

\subsection{Perfect crystals} \label{subsec:perfect}
Of particular interest is a class of crystals called perfect crystals,
which are crystals for affine algebras satisfying a set of very
special properties.  These properties ensure that perfect crystals can
be used to construct the path realization of crystals for highest
weight modules by taking the semi-infinite tensor product of these
crystals~\cite{KKMMNN:1992a}.  To define them, we need a few
preliminary definitions.

Recall that $P$ denotes the weight lattice of a Kac-Moody algebra
$\mathfrak{g}$; for the remainder of this section, $\mathfrak{g}$ is
of affine type.  The center of $\mathfrak{g}$ is one-dimensional and
is generated by the 
canonical central element $c=\sum_{i\in I}a^\vee_i h_i$, where the
$a^\vee_i$ are the numbers on the nodes of the Dynkin diagram of the
algebra dual to $\mathfrak{g}$ given in Table Aff of~\cite[section
4.8]{Kac:1990}.  Moreover, the imaginary
roots of $\mathfrak{g}$ are nonzero integer multiples of the null
root
$\delta=\sum_{i\in I} a_i \alpha_i$, where the $a_i$ are the numbers
on the nodes of the Dynkin diagram of $\mathfrak{g}$ given in Table
Aff of~\cite{Kac:1990}.  Define $P_{\cl}=P/\mathbb{Z}\delta$,
$P_{\cl}^{+}=\{\la\in P_{\cl}|\lan h_{i},\la\ran\geq 0\textrm{ for
all 
}i\in I\}$, and $\uqp{g}$ to be the quantum enveloping algebra with
the Cartan datum
$(A,\Pi,\Pi^{\vee},P_{\cl},P_{\cl}^{\vee})$.

Define the set of level $\ell$ weights to be
$(P_{\cl}^{+})_{\ell}=\{\la\in P_{\cl}^{+}|\lan c,\la \ran=\ell\}$.
For
a crystal basis element $b\in B$, define
\begin{equation*}
\vare(b)=\sum_{i\in I}\vare_{i}(b)\Lambda_{i} \quad \text{and} \quad
\varphi(b)=\sum_{i\in I}\varphi_{i}(b)\Lambda_{i},
\end{equation*}
where $\Lambda_i$ is the $i$-th fundamental weight of $\mathfrak{g}$.
Finally, for a crystal basis $B$, we define $B_{\min}$ to be the set
of crystal basis
elements $b$ such that $\lan c,\vare(b)\ran$ is minimal over $b\in B$.

Following \cite{KKMMNN:1992},
\begin{definition} \label{def:perfect}
A $\uqp{g}$-crystal $B$ is a perfect crystal of level $\ell$ if:
\begin{enumerate}
    \item $B\otimes B$ is connected;
    \item there exists $\la\in P_{\cl}$ such that $\wt(B)\subset
    \la+\sum_{i\neq0}\mathbb{Z}_{\leq0}\alpha_{i}$ and
    there is a unique $b_{\la}\in B$ such that $\wt(b_{\la}) = \la$;
    \item there is a finite-dimensional irreducible $\uqp{g}$-module
$V$ with a
    crystal base whose crystal graph is isomorphic to $B$;
    \item for any $b\in B$, we have $\langle c,\vare(b)\rangle \geq
    \ell$;
    \item the maps $\vare$ and $\varphi$ from $B_{\min}$ to
    $(P_{\cl}^{+})_{\ell}$ are bijective.
\end{enumerate}
We use the notation $\lev(B)$ to indicate the level of the perfect 
crystal $B$.
\end{definition}

\subsection{The energy function and one-dimensional sums} \label{subsec:energy}

For an affine crystal $B$, we can define
an integer-valued function on $B$ called an energy function.
Let $B_1$ and $B_2$ be finite $U'_q(\mathfrak{g})$-crystals. Then
following~\cite[Section 4]{KKMMNN:1992a}
\begin{enumerate}
\item There is a unique isomorphism of $U'_q(\mathfrak{g})$-crystals
$R=R_{B_2,B_1}:B_2\otimes B_1\rightarrow B_1\otimes B_2$.
\item There is a function $H=H_{B_2,B_1}:B_2\otimes
B_1\rightarrow\Z$, unique up to global additive constant, such
that $H$ is constant on classical components and, for all $b_2\in B_2$
and $b_1\in B_1$, if $R(b_2\otimes b_1)=b_1'\otimes b_2'$, then
\begin{equation} \label{eq:local energy}
  H(\etil{0}(b_2\otimes b_1))=
  H(b_2\otimes b_1)+
  \begin{cases}
    -1 & \text{if $\vare_0(b_2)>\varphi_0(b_1)$ and
    $\vare_0(b_1')>\varphi_0(b_2')$} \\
    1 & \text{if $\vare_0(b_2)\le\varphi_0(b_1)$ and
    $\vare_0(b_1')\le \varphi_0(b_2')$} \\
    0 & \text{otherwise.}
  \end{cases}
\end{equation}
\end{enumerate}

The maps $R$ and $H$ are called the local isomorphism and local energy
function on $B_2\otimes B_1$.  The pair $(R,H)$ is called the
combinatorial $R$-matrix.

For a finite $U'_q(\mathfrak{g})$-crystal $B$, define $u(B)$ to be the
unique vector in $B$ such that $\wt(u(B)) = \la$, $u(B)$ is the only
vector in $B$ with weight $\la$, and for all $b\in B$, $\wt(b)$ is in
the convex hull of the Weyl group action on $\la$ \cite{SS:2004}.
Then with $B_{1}$ and $B_{2}$ as above,
\begin{equation*}
  R(u(B_2)\otimes u(B_1)) = u(B_1)\otimes u(B_2).
\end{equation*}
It is convenient to normalize the local energy function $H$ by
requiring that
\begin{equation*}
  H(u(B_2)\otimes u(B_1)) = 0.
\end{equation*}
With this convention it follows by definition that
\begin{equation*}
H_{B_1,B_2} \circ R_{B_2,B_1} = H_{B_2,B_1}
\end{equation*}
as $\mathbb{Z}$-valued functions on $B_2\otimes B_1$.

We wish to define an energy function $D_B:B\rightarrow \Z$ for tensor
products of perfect crystals of the form $B^{r,s}$~\cite[Section
3.3]{HKOTT:2001}.  Let $B=B^{r,s}$; we for the moment assume that it
is perfect.  Then there exists a unique element $b^\natural\in B$ such
that $\varphi(b^\natural)=\lev(B)\La_0$.  Define
$D_B:B\rightarrow\Z$ by
\begin{equation}\label{eq:D single}
  D_B(b) = H_{B,B}(b\otimes b^\natural) - H_{B,B}(u(B)\otimes
b^\natural).
\end{equation}

The intrinsic energy $D_B$ for the $L$-fold tensor product
$B=B_L\otimes\dotsm\otimes B_1$ where $B_j=B^{r_j,s_j}$ is given by
\begin{equation*}
D_B = \sum_{1\le i<j\le L} H_i R_{i+1}R_{i+2}\dotsm R_{j-1}
 + \sum_{j=1}^L D_{B_j} R_1R_2\dotsm R_{j-1},
\end{equation*}
where $H_i$ and $R_i$ are the local energy function and $R$-matrix on
the $i^{\mathrm{th}}$ and $i+1^{\mathrm{th}}$ tensor factor, respectively.

\begin{definition}
    \label{def:formulaX}
    Let $B$ be a tensor product of crystals as above and let $\la$ be
    a dominant weight.  The one-dimensional configuration sum,
    (referred to in \cite{HKOTT:2001, HKOTY:1999} as the $X$ formula)
    is
    $$
    X(B,\la;q) = \sum_{b\in\mathcal{P}(B,\la)}q^{D_{B}(b)},
    $$
    where $\mathcal{P}(B,\la)$ is simply the set of vertices in $B$ of
    weight $\la$.
\end{definition}

One of the motivations for explicitly describing the crystals
$B^{r,s}$ is to allow for a combinatorial proof of the equivalence
between the above sum and the fermionic formula arising from the Bethe
ansatz \cite{HKOTT:2001, HKOTY:1999}.

        \section{Preliminaries on Type $D_{n}$ crystals}
        \label{sec:KR:tabs}
        For the rest of Chapters \ref{sec:KR} and \ref{sec:B2s} we restrict
our attention to the finite Lie algebra of type $D_n$ and the affine
Kac-Moody algebra of type $D_n^{(1)}$.  Denote by $I=\{0,1,\ldots,n\}$
the index set of the Dynkin diagram for $D_n^{(1)}$ (see Figure
\ref{fig:DD}) and by $J=\{1,2,\ldots,n\}$ the Dynkin diagram
for type $D_n$.

\subsection{Dynkin data}
For type $D_n$, the simple roots are
\begin{equation}\label{eq:alpha}
\begin{split}
\alpha_i&=\epsilon_i-\epsilon_{i+1} \qquad \text{for $1\le i<n$}\\
\alpha_n&=\epsilon_{n-1}+\epsilon_n
\end{split}
\end{equation}
and the fundamental weights are
\begin{equation*}
\begin{aligned}
\om_i&=\epsilon_1+\cdots+\epsilon_i &\text{for $1\le i\le n-2$}\\
\om_{n-1}&=(\epsilon_1+\cdots+\epsilon_{n-1}-\epsilon_n)/2 &\\
\om_n&=(\epsilon_1+\cdots+\epsilon_{n-1}+\epsilon_n)/2&
\end{aligned}
\end{equation*}
where $\epsilon_i\in \mathbb{Z}^n$ is the $i$-th unit standard vector.
The central element for $D_n^{(1)}$ is
$$c=h_{0}+h_{1}+2h_{2}+\cdots+2h_{n-2}+h_{n-1}+h_{n}.$$

\subsection{Classical crystals}\label{sec:class crys}

Kashiwara and Nakashima~\cite{KN} described the crystal structure of
all classical highest weight crystals $B(\Lambda)$ of highest weight
$\Lambda$ explicitly.  Kirillov-Reshetikhin crystals corresponding to
the last two fundamental weights (the ``spin'' weights) are treated as
a special case in \cite{HKOTT:2001, HKOTY:1999}; thus in the sequel,
we only need to consider representations of $D_{n}$ that are quotients
of tensor powers of the vector representation, and omit spin
representations from our consideration here.  Recall that we may
interpret such a dominant weight $\La$ as a partition with $\lan
h_{i},\La\ran$ columns of height $i$.  The corresponding crystals can
be represented by tableaux on this partition using the partially
ordered alphabet
\begin{displaymath}
1<2<\cdots < n-1 < \stack{n}{\bar{n}} < \overline{n-1}<\cdots 
\bar{2}<\bar{1}
\end{displaymath}
where $n$ and $\bar{n}$ are incomparable, such that the following
conditions hold~\cite[page 202]{HK}:
\begin{criterion} \label{legal}
{}~
\begin{enumerate}
    \item \label{legal:1} If 
    $\begin{array}{|c|c|}
        \hline
        a & b  \\
        \hline
         
    \end{array}
    $ 
    is in the tableau, then $a\leq b$;
    
    \item \label{legal:2} If $\begin{array}{|c|}
        \hline
        a  \\
        \hline
        b  \\
        \hline
         
    \end{array}
    $ is in the tableau, then $b\nleq a$;
    
    \item \label{legal:2.5} If $a$ and $\bar{a}$ appear in a 
    column of length $N$, with $a$ in row $p$ and $\bar{a}$ in row 
    $q$, we have
    $$
    (q-p)+a > N;
    $$
    
    \item \label{legal:3} If the tableau contains a configuration of 
    the form
    $$
    \begin{array}{cc|c|c|}
        \cline{3-4}
        p\to&& a &  \vdots  \\
        \cline{3-4}
        q\to&\multicolumn{1}{c}{} & &b  \\
        \cline{4-4}
        &\multicolumn{1}{c}{} & &\vdots  \\
        \cline{4-4}
        r\to&\multicolumn{1}{c}{} & &\bar{b}  \\
        \cline{4-4}
        &\multicolumn{1}{c}{} & &\vdots  \\
        \cline{4-4}
        s\to&\multicolumn{1}{c}{} & &\bar{a}  \\
        \cline{4-4}
         
    \end{array}
    \qquad
    \textrm{ or }
    \qquad
    \begin{array}{c|c|c|c}
        \cline{2-2}
        p\to&a & \multicolumn{2}{c}{}  \\
        \cline{2-2}
        &\vdots &  \multicolumn{2}{c}{}   \\
        \cline{2-2}
        q\to&b &  \multicolumn{2}{c}{}   \\
        \cline{2-2}
        &\vdots &  \multicolumn{2}{c}{}   \\
        \cline{2-2}
        r\to&\bar{b} &  \multicolumn{2}{c}{}   \\
        \cline{2-3}
        s\to&\vdots & \bar{a} &   \\
        \cline{2-3}
         
    \end{array}
    $$
    with $1 \leq a < b < n$ then 
    $$
    (q-p) + (s-r) < b-a
    $$
    Note that we may have $p=q$ in the first case or $s=r$ in the 
    second case above;
    
    \item \label{legal:4} No configuration of the form
        $$
    \begin{array}{cc|c|c|}
        \cline{3-4}
        && a &  a  \\
        \cline{3-4}
        &\multicolumn{1}{c}{} & &\vdots  \\
        \cline{4-4}
        &\multicolumn{1}{c}{} & &\bar{a}  \\
        \cline{4-4}
         
    \end{array}
    \qquad
    \textrm{ or }
    \qquad
    \begin{array}{c|c|c|c}
        \cline{2-2}
        &a & \multicolumn{2}{c}{}  \\
        \cline{2-2}
        &\vdots &  \multicolumn{2}{c}{}   \\
        \cline{2-3}
        &\bar{a} & \bar{a} &   \\
        \cline{2-3}
         
    \end{array}
    $$
    with $1\leq a < n$ may appear in the tableau;
    
    \item \label{legal:5}  If the tableau contains a configuration of 
    the form
    $$
    \begin{array}{cc|c|c|}
        \cline{3-4}
        p\to&& a &  \vdots  \\
        \cline{3-4}
        &\multicolumn{1}{c}{} & &n  \\
        \cline{4-4}
        &\multicolumn{1}{c}{} & &\bar{n}  \\
        \cline{4-4}
        &\multicolumn{1}{c}{} & &\vdots  \\
        \cline{4-4}
        s\to&\multicolumn{1}{c}{} & &\bar{a}  \\
        \cline{4-4}
         
    \end{array}
    \qquad
    \textrm{ or }
    \begin{array}{c|c|c|c}
        \cline{2-2}
        p\to&a & \multicolumn{2}{c}{}  \\
        \cline{2-2}
        &\vdots &  \multicolumn{2}{c}{}   \\
        \cline{2-2}
        &n &  \multicolumn{2}{c}{}   \\
        \cline{2-2}
        &\bar{n} &  \multicolumn{2}{c}{}   \\
        \cline{2-3}
        s\to&\vdots & \bar{a} &   \\
        \cline{2-3}
         
    \end{array}
    \textrm{ or }
    \begin{array}{cc|c|c|}
        \cline{3-4}
        p\to&& a &  \vdots  \\
        \cline{3-4}
        &\multicolumn{1}{c}{} & &\bar{n}  \\
        \cline{4-4}
        &\multicolumn{1}{c}{} & &n  \\
        \cline{4-4}
        &\multicolumn{1}{c}{} & &\vdots  \\
        \cline{4-4}
        s\to&\multicolumn{1}{c}{} & &\bar{a}  \\
        \cline{4-4}
         
    \end{array}
    \qquad
    \textrm{ or }
    \begin{array}{c|c|c|c}
        \cline{2-2}
        p\to&a & \multicolumn{2}{c}{}  \\
        \cline{2-2}
        &\vdots &  \multicolumn{2}{c}{}   \\
        \cline{2-2}
        &\bar{n} &  \multicolumn{2}{c}{}   \\
        \cline{2-2}
        &n &  \multicolumn{2}{c}{}   \\
        \cline{2-3}
        s\to&\vdots & \bar{a} &   \\
        \cline{2-3}
         
    \end{array}    
    $$
    then 
    $$
    s-p-1 < n-a
    $$
    as with condition \ref{legal:3}, the number of blocks represented 
    by the $\vdots$ at the elbow of these configurations may be any 
    non-negative integer;    
    
    \item \label{legal:6} if the tableau contains an $n$ or 
    $\bar{n}$, then neither $n$ nor $\bar{n}$ may appear in the region
    strictly down and to the right of it;
    
    \item \label{legal:7} If the tableau has a configuration of the 
    form
    $$
    \begin{array}{cc|c|c|}
        \cline{3-4}
        p\to&& a &  \vdots  \\
        \cline{3-4}
        q\to& &\vdots &n  \\
        \cline{3-4}
        & & \vdots& \vdots  \\
        \cline{3-4}
        r\to& &\bar{n} &\vdots  \\
        \cline{3-4}
        s\to& & \vdots &\bar{a}  \\
        \cline{3-4}
         
    \end{array}
    \qquad
    \textrm{ or }
    \qquad
    \begin{array}{cc|c|c|}
        \cline{3-4}
        p\to&& a &  \vdots  \\
        \cline{3-4}
        q\to& &\vdots &\bar{n}  \\
        \cline{3-4}
        & & \vdots& \vdots  \\
        \cline{3-4}
        r\to& &n &\vdots  \\
        \cline{3-4}
        s\to& & \vdots &\bar{a}  \\
        \cline{3-4}
         
    \end{array}
    $$    
    with $r - q$ odd or  
    $$
    \begin{array}{cc|c|c|}
        \cline{3-4}
        p\to&& a &  \vdots  \\
        \cline{3-4}
        q\to& &\vdots &n  \\
        \cline{3-4}
        & & \vdots& \vdots  \\
        \cline{3-4}
        r\to& &n &\vdots  \\
        \cline{3-4}
        s\to& & \vdots &\bar{a}  \\
        \cline{3-4}
         
    \end{array}
    \qquad
    \textrm{ or }
    \qquad
    \begin{array}{cc|c|c|}
        \cline{3-4}
        p\to&& a &  \vdots  \\
        \cline{3-4}
        q\to& &\vdots &\bar{n}  \\
        \cline{3-4}
        & & \vdots& \vdots  \\
        \cline{3-4}
        r\to& &\bar{n} &\vdots  \\
        \cline{3-4}
        s\to& & \vdots &\bar{a}  \\
        \cline{3-4}
         
    \end{array}
    $$    
    with $r - q$ even, we have 
    $$
    s-p < n - a.
    $$
    
\end{enumerate}
\end{criterion}

A tableau that satisfies Criterion \ref{legal} is called a type
$D_{n}$ tableau.

The crystal $B(\om_1)$, (also called the vector representation) is
described pictorially by the crystal graph:

\vspace{2in}

\begin{center}\setlength{\unitlength}{1.8cm}
    \begin{picture}(0,0)
	\put(-3.2,0){\framebox{$1$}}
	\put(-2.9,.1){\vector(1,0){.4}}
	\put(-2.75,.16){\scriptsize{$1$}}
	\put(-2.45,0){\framebox{$2$}}
	\put(-2,.02){$\cdots$}
	\put(-1.6,.1){\vector(1,0){.4}}
	\put(-1.61,.16){\scriptsize{$n-2$}}
	\put(-1.15,0){\framebox{$n-1$}}
	\put(-.55,.3){\vector(1,1){.3}}
	\put(-.88,.56){\scriptsize{$n-1$}}
	\put(-.55,-.15){\vector(1,-1){.3}}
	\put(-.53,-.47){\scriptsize{$n$}}
	\put(-.19,.6){\framebox{$n$}}
	\put(.17,.6){\vector(1,-1){.26}}
	\put(.33,.56){\scriptsize{$n$}}
	\put(-.19,-.6){\framebox{$\overline{n}$}}
	\put(.17,-.44){\vector(1,1){.29}}
	\put(.33,-.47){\scriptsize{$n-1$}}
	\put(0.3,0){\framebox{$\overline{n-1}$}}
	\put(1.,.1){\vector(1,0){.4}}
	\put(1.01,.16){\scriptsize{$n-2$}}
	\put(1.62,.02){$\cdots$}
	\put(2.05,0){\framebox{$\overline{2}$}}
	\put(2.35,.1){\vector(1,0){.4}}
	\put(2.50,.16){\scriptsize{$1$}}
	\put(2.8,0){\framebox{$\overline{1}$}}	
    \end{picture}
\end{center}
\vspace{.5in} 

\begin{definition}
    \label{def:col_word}
    Let $T$ be a type $D_{n}$ tableau.  The column word of $T$, denoted by
    $w_{T}$, is the word on the alphabet
    $\{1,2,\ldots,n,\bar{n},\ldots,\bar{1}\}$ that results from
    reading the entries in the tableau beginning with the left-most
    column and proceeding to the right, reading each column from
    bottom to top.
\end{definition}

If we interpret the column word as a label for a basis vector of a
quotient of a tensor power of the vector representation, the action of
the Kashiwara operators $\ftil{i}$ and $\etil{i}$ on a tableau is
determined by the tensor product rule as defined in
section~\ref{sec:crystal bases}.

\begin{example}
    Let $n=4$.  Then the tableau
    $$
    T=\left.
    \begin{array}{|c|c|c|c|c|}
	\hline
	1 & 2 & 4 & \bar{3} & \bar{3}  \\
	\hline
	3 & \bar{4} & \bar{4} & \bar{2} & \bar{1}  \\
	\hline 
    \end{array}
    \right.
    $$
    has column word
    $w_{T}=31\bar{4}2\bar{4}4\bar{2}\bar{3}\bar{1}\bar{3}$.  The
    $2$-signature of $T$ is $+-+--$, derived from the subword
    $32\bar{2}\bar{3}\bar{3}$, and the reduced $2$-signature is a single
    $-$.  Therefore $e_{2}(T) = 0$ and
    $$
    \ftil{2}(T)=\left.
    \begin{array}{|c|c|c|c|c|}
	\hline
	1 & 2 & 4 & \bar{3} & \mathbf{\bar{2}}  \\
	\hline
	3 & \bar{4} & \bar{4} & \bar{2} & \bar{1}  \\
	\hline 
    \end{array}
    \right. ,
    $$
    since the rightmost $-$ in the reduced $2$-signature of $T$ comes from
    the northeastmost $\bar{3}$.  The $4$-signature of $T$ is $-++-++$,
    derived from the subword $3\bar{4}\bar{4}4\bar{3}\bar{3}$, and the
    reduced $4$-signature is $-+++$, from the subword
    $3\bar{4}\bar{3}\bar{3}$.  This tells us that
    $$
    \ftil{4}(T)=\left.
    \begin{array}{|c|c|c|c|c|}
	\hline
	1 & 2 & 4 & \bar{3} & \bar{3}  \\
	\hline
	\mathbf{\bar{4}} & \bar{4} & \bar{4} & \bar{2} & \bar{1}  \\
	\hline 
    \end{array}
    \right.
    \quad\textrm{ and }\quad
    \etil{4}(T)=\left.
    \begin{array}{|c|c|c|c|c|}
	\hline
	1 & 2 & 4 & \bar{3} & \bar{3}  \\
	\hline
	3 & \mathbf{3} & \bar{4} & \bar{2} & \bar{1}  \\
	\hline 
    \end{array}
    \right. .
    $$    
\end{example}

\subsection{Plactic monoid of type $D$}\label{sec:plac}
The plactic monoid for type $D$ is the free monoid generated by the
letters $\{1, \ldots, n, \bar{n},$ $\ldots, \bar{1}\}$, modulo certain
relations introduced by Lecouvey~\cite{L}.  These Lecouvey relations
generalize the Knuth relations governing the RSK correspondence.  Note
that we write the letters of our words in reverse of the order used
in~\cite{L}.  An admissible column is one whose column word $C=x_L
x_{L-1} \cdots x_1$ is such that $x_{i+1}\not\le x_i$ for
$i=1,\ldots,L-1$.  Note that the letters $n$ and $\bar{n}$ are the
only letters that may appear more than once in $C$.  Let $z\le n$ be a
letter in $C$.  Then $N(z)$ denotes the number of letters $x$ in $C$
such that $x\le z$ or $x\ge \bar{z}$.  A column $C$ is called
admissible if $L\le n$ and for any pair $(z,\bar{z})$ of letters in
$C$ with $z\le n$ we have $N(z)\le z$.  The Lecouvey $D$ equivalence
relations are given by:
\begin{enumerate}
    \item  If $x\neq \bar{z}$, then
    \begin{displaymath}
	xzy\equiv zxy \textrm{ for }x\leq y<z
	\textrm{ and }
	yzx\equiv yxz \textrm{ for }x<y\leq z	.
    \end{displaymath}
    \item  If $1<x<n$ and $x\leq y\leq\bar{x}$, then 
    \begin{displaymath}
	(x-1)(\overline{x-1})y\equiv\bar{x}xy\textrm{ and }
	y\bar{x}x\equiv y(x-1)(\overline{x-1}).
    \end{displaymath}
    \item  If $x\leq n-1$, then
    \begin{displaymath}
	\left\{
	\begin{array}{c}
	    n\bar{x}\bar{n}\equiv n\bar{n}\bar{x}\\
	    \bar{n}\bar{x}n\equiv \bar{n}n\bar{x}\\
	\end{array}
	\right.
	\textrm{ and }
	\left\{
	\begin{array}{c}
	    xn\bar{n}\equiv nx\bar{n}\\
	    x\bar{n}n\equiv \bar{n}xn\\
	\end{array}
	\right. .
    \end{displaymath}
    \item 
    \begin{displaymath}
	\left\{
	\begin{array}{c}
	    \bar{n}\bar{n}n\equiv\bar{n}(n-1)(\overline{n-1})\\
	    nn\bar{n}\equiv n(n-1)(\overline{n-1})\\
	\end{array}
	\right.
	\textrm{ and }
	\left\{
	\begin{array}{c}
	    (n-1)(\overline{n-1})\bar{n}\equiv n\bar{n}\bar{n}\\
	    (n-1)(\overline{n-1})n\equiv \bar{n}nn\\
	\end{array}
	\right. .
    \end{displaymath}
    \item 
    Consider $w$ a non-admissible column word such that every
    consecutive subword of $w$ with length less than the length of $w$
    is admissible.  Let $z$ be the lowest unbarred letter such that
    the pair $(z,\bar{z})$ occurs in $w$ and $N(z)>z$.  Then
    $w\equiv\tilde{w}$ is the column word obtained by erasing the pair
    $(z,\bar{z})$ in $w$ if $z<n$, by erasing a pair $(n,\bar{n})$ of
    consecutive letters otherwise.
\end{enumerate}
This monoid gives us a bumping algorithm similar to the Schensted
bumping algorithm.  It is noted in~\cite{L} that a general type $D$
sliding algorithm, if one exists, would be very complicated.  However,
for tableaux with no more than two rows, we can derive such an
algorithm, the details of which are given in Section
\ref{sec:B2s:tabs}.

\subsection{Dual crystals}\label{sec:dual crystal}
Let $\omega_0$ be the longest element in the Weyl group of $D_n$.
The action of $\omega_0$ on the weight lattice $P$ of $D_n$ is given by
\begin{equation*}
\begin{split}
\omega_0(\om_i) &= -\om_{\tau(i)}\\
\omega_0(\alpha_i) &= -\alpha_{\tau(i)}
\end{split}
\end{equation*}
where $\tau:J\to J$ is the identity if $n$ is even and interchanges
$n-1$ and $n$ and fixes all other Dynkin nodes if $n$ is odd. 

For any $D_{n}$-crystal $B$ there is a unique involution $\dual:B\to
B$, called the dual map, satisfying
\begin{equation*}
\begin{split}
\wt(b^\dual) &= \omega_0 \wt(b)\\
\etil{i}(b)^\dual &= \ftil{\tau(i)}(b^\dual)\\
\ftil{i}(b)^\dual &= \etil{\tau(i)}(b^\dual).
\end{split}
\end{equation*}
The involution $\dual$ sends the highest weight vector $u\in B(\La)$ to the 
lowest weight vector (the unique vector in $B(\La)$ of weight $\omega_0(\La)$).
We have
\begin{equation*}
(B_1\otimes B_2)^\dual \cong B_2\otimes B_1
\end{equation*}
with $(b_1\otimes b_2)^\dual \mapsto b_2^\dual \otimes b_1^\dual$.

Explicitly, on $B(\om_1)$ the involution $\dual$ is given by 
\begin{equation*}
 i \longleftrightarrow \overline{i}
\end{equation*}
except for $i=n$ with $n$ odd in which case $n \leftrightarrow n$ and
$\overline{n} \leftrightarrow \overline{n}$.  For $T\in B(\La)$ the
dual vertex $T^{\dual}$ is obtained by applying the $\dual$ map
defined for $B(\om_1)$ to each of the letters of
$w_{T}^{\textrm{rev}}$ (the word with the same letters as $T$ taken in
the reverse order), and then rectifying the resulting word by the
Lecouvey $D$ equivalence relations.

\begin{example}
If
\begin{equation*}
T=\begin{array}{|c|c|c|} \hline 1&1&2\\ \hline \overline{3} & \multicolumn{2}{c}{}\\
\cline{1-1}\end{array}
\in B(2\om_1+\om_2)
\end{equation*}
we have
\begin{equation*}
T^{*}=\begin{array}{|c|c|c|} \hline 3&\overline{1}&\overline{1}\\
\hline \overline{2} & \multicolumn{2}{c}{}\\ \cline{1-1}\end{array}.
\end{equation*}
\end{example}

\subsection{Properties of $B^{r,s}$}
As mentioned in the introduction, it was conjectured
in~\cite{HKOTT:2001,HKOTY:1999} that there are crystal bases $B^{r,s}$
associated with Kirillov--Reshetikhin modules $W^{r,s}$.  In addition
to the existence, Hatayama et al.~\cite{HKOTT:2001} conjectured
certain properties of $B^{r,s}$ which we state here.

\begin{conjecture}[\cite{HKOTT:2001}]\label{conj:Brs}
If the crystal $B^{r,s}$ of type $D_n^{(1)}$ exists, it has the following
properties:
\begin{enumerate}
\item As a classical crystal $B^{r,s}$ decomposes as
$\bks{r}\cong \bigoplus_{\La}B(\La)$, where the direct sum is taken 
over all partitions $\La$ whose complement in an $r\times s$ 
rectangle is tiled by vertical dominos.
\item $B^{r,s}$ is perfect of level $s$.
\item $B^{r,s}$ is equipped with an energy function
$D_{B^{r,s}}$ such that $D_{B^{r,s}}(b)=k-s$ if $b$ is in a $D_{n}$ 
component $B(\La)$ for which $\La$ differs from an $r\times s$ 
rectangle by $k$ vertical dominos.

\end{enumerate}
\end{conjecture}

\begin{remark}
    The 
    decomposition of these Kirillov-Reshetikhin modules as modules
    over the embedded classical algebra of type $D_{n}$ specified by
    Item 1 of Conjecture \ref{conj:Brs} was confirmed for the case of
    untwisted algebras in \cite{Chari} and for twisted algebras in
    \cite{NS:2005}; therefore if the crystals exist, this is known to
    be true.
\end{remark}

\begin{remark}
    In \cite{FSS}, it is proved that given the following assumptions, 
    the affine structure of $B^{r,s}$ is specified by the embedding of
    certain Demazure crystals in $B^{r,s}$.  We specialize their
    statement to the type $D$ case.
    \begin{enumerate}
	\item $B^{r,s}$ is a regular crystal (i.e., eliminating all
	edges of color $\{i_{1},\ldots,i_{m}\}\subset I$ results in a
	crystal over
	$\mathfrak{g}_{I\setminus\{i_{1},\ldots,i_{m}\}}$),
    
        \item There is a unique element $u\in B^{r,s}$ such that 
	$$
	\vare(u) = s\Lz\qquad\mathrm{and}\qquad\varphi(u) =
	s\La_{r(\mathrm{mod}2)},
	$$

        \item $B^{r,s}$ admits the automorphism corresponding to
	$\sigma$, the Dynkin diagram automorphism interchanging $0$
	and $1$, leaving all other nodes fixed.
    \end{enumerate}

\end{remark}

	\section[A Dynkin Diagram Automorphism and the Branching
	Component Graph]{A Dynkin Diagram Automorphism and the\\
	Branching Component Graph}
        \label{sec:KR:BC}
        \subsection{Preliminaries}

\begin{figure}
    \vspace{.4cm}
\begin{center}\setlength{\unitlength}{1cm}
    \begin{picture}(0,0)
    \put(-2,0){$\circ$}
    \put(-1.8,.1){$0$}
    \put(-1.9,0.04){\line(1,-1){.43}}
    \put(-2,-1){$\circ$}
    \put(-1.8,-1.1){$1$}
    \put(-1.87,-.87){\line(1,1){.42}}
    \put(-1.5,-.5){$\circ$}
    \put(-1.5,-.3){$2$}
    \put(-1.36,-.428){\line(1,0){.5}}
    \put(-.9,-.5){$\circ$}
    \put(-.9,-.3){$3$}
    \put(-.5,-.51){$\cdot$}
    \put(-.2,-.51){$\cdot$}
    \put(.1,-.51){$\cdot$}
    \put(.42,-.5){$\circ$}
    \put(.15,-.3){$n-3$}
    \put(.55,-.428){\line(1,0){.51}}
    \put(1.02,-.5){$\circ$}
    \put(1.2,-.5){$n-2$}
    \put(1.56,0.04){\line(-1,-1){.41}}
    \put(1.523,0){$\circ$}
    \put(1.723,.1){$n-1$}
    \put(1.523,-1){$\circ$}
    \put(1.723,-1.1){$n$}
    \put(1.56,-.9){\line(-1,1){.43}}
    \qbezier[200](-2,0)(-2.4,-.43)(-2,-.86)
    \put(-1.97,.02){\vector(1,1){0}}
    \put(-1.98,-.87){\vector(1,-1){0}}    
    \end{picture}
\end{center}
\vspace{1cm}
\caption{Dynkin diagram automorphism}\label{fig:DD}
\end{figure}
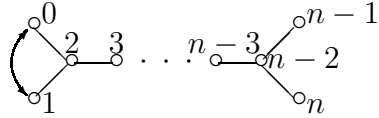

We know that there is a map $\sigma:\brs\to\brs$ such that
$\etil{0}=\sigma\etil{1}\sigma$ and $\ftil{0}=\sigma\ftil{1}\sigma$,
arising from the automorphism of $\uqsophat{2n}$ which interchanges
nodes $0$ and $1$ of its Dynkin diagram (Figure \ref{fig:DD}).  Let
$J\subset I$, and denote by $B_{J}$ the graph that results from
removing all $j$-colored edges from $\brs$ for $j\in J$.  Then as
directed graphs, $B_{\{0\}}$ (which has the structure of a $D_{n}$
crystal) is isomorphic to $B_{\{1\}}$.  This motivates the
introduction of the ``branching component graph'' of $\brs$ in
Definition \ref{def:BC} below, which records all the relevant
information about $B_{\{0,1\}}$.

It is easy to see that the connected components of $B_{\{0,1\}}$ will
be $D_{n-1}$-crystals whose highest weights correspond to partitions
as characterized by the propositions below.  The branching component
graph for $\brs$, which we denote $\BC(\brs)$, is determined by this
decomposition.  We will describe the decomposition into irreducible
$D_{n-1}$-crystals of a $D_{n}$-crystal indexed by a partition $\la$
whose complement in an $r\times s$ rectangle is tiled by vertical
dominos, calling the resulting graph $\BC(\la)$.  Given the result of
Chari specifying the decomposition of $B^{r,s}$ as a module over the
classical algebra, $\BC(\brs)$ is simply the union of these graphs.

\begin{definition}\label{def:BC}
    Let $\la$ be a partition labeling a $D_{n}$ highest weight.  The
    branching component graph $\BC(\la)$ has one vertex for each
    connected $D_{n-1}$-crystal that results from removing the $1$
    arrows from $B(\la)$.  We associate to each vertex the partition
    corresponding to the $D_{n-1}$ highest weight of that crystal.
    For a vertex $v\in\BC(\la)$, let $B(v)$ denote the set of crystal
    vertices in $B(\la)$ contained in the $D_{n-1}$-crystal indexed by
    $v$.  Then $\BC(\la)$ has an edge from $v$ to $w$ if there is a
    crystal vertex $b\in B(v)$ such that $\ftil{1}(b)\in B(w)$.
\end{definition}

This definition is essentially a refinement of the classically known
$D_{n}$ to $D_{n-1}$ branching rule, which is as follows: let $\La$ be
a $D_{n}$ dominant weight with no spin part (i.e., such that $\lan
h_{n},\La\ran = \lan h_{n-1},\La\ran = 0$) and let $V_{\La}$ denote
the irreducible $D_{n}$ module with highest weight $\La$.  Recall that
we may interpret $\La$ as a partition with $\lan h_{i},\La\ran$
columns of height $i$.  A horizontal strip is a subset $H$ of a
partition $\la$ such that distinct elements of $H$ are in distinct
columns of $\la$ and $\la\backslash H$ is a partition.  Then as a
$D_{n-1}$ module, $V_{\La}$ is the direct sum of $D_{n-1}$ modules
$V_{\mu}$ with $D_{n-1}$ highest weight $\mu$ for each $\mu$ that can
be obtained by successively removing two horizontal strips from $\La$.
As an illustration, see Example \ref{ex:branch}.  This rule can be
obtained by combining the $D_{n}$ to $B_{n-1}$ branching rule with the
$B_{n-1}$ to $D_{n-1}$ branching rule \cite{FH}.

This motivates the following definition, formulated by Schilling and
Shimozono \cite{SSh:2006}.

\begin{definition}
    \label{def:pmdiag}
    Consider a sequence of partitions $\la\subset\mu\subset\La$ such
    that $\mu/\la$ and $\La/\mu$ are horizontal strips.
    A $\pm$ diagram $D$ of shape $\La/\la$ is the partition $\La$ with
    a $+$'s in each box of $\mu/\la$ and $-$'s in each box of
    $\La/\mu$.  We call $\La=\mathrm{outer}(D)$ and
    $\la=\mathrm{inner}(D)$ the outer and inner shape of $D$,
    respectively.
\end{definition}

The set of $\pm$ diagrams with outer shape $\La$ is in natural
bijection with the $D_{n-1}$ highest weights of the $D_{n-1}$ to
$D_{n}$ branching components of the $D_{n}$ module $V(\La)$; i.e., as
$D_{n-1}$ modules we have the following isomorphism:
$$
V_{D_{n}}(\La) \cong \bigoplus_{
\stackrel{\pm\textrm{ diagrams }D}
{\mathrm{outer}(D)=\La}
}V_{D_{n-1}}(\mathrm{inner}(D)),
$$
where the subscript on $V$ indicates the algebra that defines the
representation.

\begin{example}
    \label{ex:branch}
    Consider the $D_{n}$ module with highest weight $\La_{2}+\La_{4}$, 
    which corresponds to the partition $(2, 2, 1, 1)$.  As a $B_{n-1}$
    module,  this is the multiplicity-free sum of the modules labeled
    by $(2, 2, 1, 1)$,  $(2, 2, 1, 0)$,  $(2, 1, 1, 1)$,  and $(2, 1, 1, 0)$.  As a
    $D_{n-1}$ module,  it is the sum of the modules associated with the
    following partitions with the given multiplicities.
    \begin{center}
    \begin{tabular}{c|c}
        Partition & Multiplicity \\
        \hline
        (2, 2, 1, 1) & 1  \\
        \hline
        (2, 2, 1) & 2  \\
        \hline
        (2, 1, 1, 1) & 2  \\
        \hline
        (2, 1, 1) & 4  \\
        \hline
        (2, 2) & 1  \\
        \hline
        (2, 1) & 2  \\
        \hline	
        (1, 1, 1, 1) & 1  \\
        \hline
        (1, 1, 1) & 2  \\
        \hline
        (1, 1) & 1  \\
    \end{tabular}
\end{center}    
\end{example}

This correspondence is made explicit by the following algorithm.

\begin{algorithm}
    \label{alg:pmbranch}
    
    Let $d$ be the difference between the number of $+$'s
    and the number of $-$'s in the $\pm$ diagram.

    \begin{enumerate}
    
        \item Place the maximal number of $1$'s and $\bar{1}$'s in the
	diagram such that:
	\begin{itemize}
	    \item the difference between the number of $1$'s and the
	    number of $\bar{1}$'s is $d$,
	
	    \item the $\bar{1}$'s are placed in boxes containing
	    $-$'s, starting from the right, and

	    \item the result is a legal tableau;
	\end{itemize}

	\item given that all letters $1, \bar{1}, 2, \bar{2}, \ldots,
	i-1, \overline{i-1}$ have been placed, fill the empty boxes of row
	$i-1$ with $i$'s;
    
	\item place the maximal number of $i$'s in row $i$ and
	$\bar{i}$'s in the unfilled boxes with $-$'s, respectively,
	such that
	\begin{itemize}
	    \item the resulting $i$-weight of the tableau 
	    agrees with the inner shape of the $\pm$ diagram,
	
	    \item the $\bar{i}$'s are placed in boxes containing
	    $-$'s, starting from the right, and

	    \item the result is a legal tableau.
	\end{itemize}    

    \end{enumerate}
    
\end{algorithm}

\begin{proposition}
    Algorithm \ref{alg:pmbranch} bijectively maps the set of $\pm$
    diagrams of shape $\La$ to the set of $D_{n-1}$ highest weight
    tableaux of shape $\La$.

\end{proposition}

\begin{proof}
    The output of the algorithm is easily verified to be a $D_{n-1}$
    highest weight tableaux.  The classical branching rule tells us
    that there are precisely as many such tableaux as $\pm$ diagrams.
    By the pigeonhole principle, this correspondence must be
    bijective.
\end{proof} 

The notion of ``$\alpha_{1}$-height'' of a $D_{n-1}$ irreducible
component $B(\la)$ as a subset of a $D_{n}$ crystal $B(\La)$ will play
a crucial role in our understanding of branching component graphs.
There are three equivalent ways to understand this statistic, each of
which is useful in certain contexts.

Let $T$ be an arbitrary element of $B(\la)\subset B(\La)$.  The
definition of the $\alpha_{1}$-height of $B(\la)$ is simply the value
of $\lan h_{1}, \mu\ran$ where $\mu$ is the $D_{n}$ weight of $T$.
This is equal to the difference between the number of columns of $\La$
and the number of $1$-arrows in a minimal path in the crystal graph
between the highest weight tableau and $T$.  Finally, we can also take
this to be the number of $1$'s minus the number of $\bar{1}$'s in $T$,
which is also the difference between $+$'s and $-$'s in the $\pm$
diagram corresponding to $B(\la)$.  In Algorithm \ref{alg:pmbranch},
this is the statistic $d$.

Because the edges of the branching component graph indicate the action
of $f_{1}$, which decreases the $\alpha_{1}$-height by $1$, this graph
is the Hasse diagram of a ranked poset.  To avoid confusion with
algebra rank, we use the term ``stratum'' to describe the poset ranks.
Uniformly labeling the strata by their $\alpha_{1}$-height, we have
that the stratum that is equidistant from the ``highest weight'' and
``lowest weight'' branching component vertices is designated ``stratum
$0$''; the strata in the direction of the highest weight vertex are
labeled by consecutive positive integers, those in the opposite
direction by negative integers.

\begin{proposition}
    \label{prop:BCsym}
    If $\la=(\la_{1}, \la_{2}, \ldots, \la_{\ell})$,  then $\BC(\la)$ has
    $2\la_{1}+1$ strata.  Furthermore,  the partitions labeling
    stratum $j$ also label stratum $-j$ for $0\leq j\leq \la_{1}$, 
    including multiplicities.
\end{proposition}

\begin{proof}
    The first of these follows from the observation that the
    $\alpha_{1}$-height ranges from $\la_{1}$ to $-\la_{1}$ in
    increments of $1$.  The second claim is a consequence of the
    symmetry of the weight lattice of any classical finite dimensional
    representation.
\end{proof}

It therefore suffices to describe strata $0$ through $\la_{1}$ for any
$\BC(\la)$, appealing to symmetry for the negative strata.

\begin{proposition}
    \label{prop:onebox}
    Let $v$ and $w$ be vertices of a branching component graph with an
    edge between them.  Then the partitions associated to $v$ and $w$ 
    are adjacent in Young's lattice; i.e., they differ by a single box.
\end{proposition}

\begin{proof}
    Observe that the Kashiwara operator $f_{1}$ only changes a single
    box of a tableau, either changing a $1$ to a $2$ or a $\bar{2}$ to
    a $\bar{1}$.  It follows that the edges of the branching component
    graph only connect partitions that differ by a single box.  In the
    case of changing a $1$ to a $2$, the partition associated with the
    new branching component vertex has one more box than the old one
    if the $-$ in the $2$-signature from the new $2$ is unbracketed;
    otherwise it has one less box.  Conversely, when changing a
    $\bar{2}$ to a $\bar{1}$, the partition associated with the new
    branching component vertex has one less box than the old one if
    the $+$ in the $2$ signature from the disappearing $\bar{2}$ is
    unbracketed; otherwise it has one more box.
\end{proof}

When characterizing branching component graphs, it therefore suffices
to show that all edges described in the statements of Proposition
\ref{prop:recs} and Conjecture \ref{prop:allparts} are in fact
present.  We do not need to show that pairs of vertices whose
partitions differ by more than one box are non-adjacent.

	\subsection{Rectangular partitions}

\begin{proposition}\label{prop:recs} 
    Let $\la=(s, s, \ldots, s)$ be a $k\times s$ rectangle.  Then
    $\BC(\la)$ is as follows: stratum $s$ of $\BC(\la)$ has a single
    vertex, which is labeled by a $(k-1)\times s$ rectangle.  The
    vertices of stratum $j-1$ for $j\geq 1$ are labeled by one of
    each partition contained in $\la$ that differs from some partition
    in stratum $j$ by a single box and whose first $k-2$ parts are
    $s$.  Two vertices in strata $j$ and $j-1$ labeled by partitions
    $\mu_{1}$ and $\mu_{2}$ have an edge between them if $\mu_{1}$ and
    $\mu_{2}$ are adjacent in Young's lattice.

\end{proposition}

Before proving this statement we examine the top two strata of
$\Sup(s\Lt)$ in detail.  The highest weight branching component vertex
(i.e., the branching component vertex containing the highest weight
vertex) is indexed by the one-part partition $(s)$.  To see that this
is true, simply observe that the highest weight tableau of $B(s\Lt)$
is $\underbrace{ \stack{1}{2}\cdots \stack{1}{2}}_{s}, $ and acting by
$\ftil{2}, \ldots, \ftil{n}$ in the most general possible way will
affect only the bottom row.  When we map these bottom row subtableaux
component\-wise by $a\mapsto a-1$ and $\bar{a}\mapsto \overline{a-1}$
to tableaux of shape $(s)$, and apply the same map to the colors of
the arrows, this is clearly isomorphic to the $D_{n-1}$-crystal with
highest weight $s\om_1$.

According to Algorithm \ref{alg:pmbranch} and the $D_{n}$ to $D_{n-1}$
branching rule, the two $D_{n-1}$ highest weight tableaux in stratum
$s-1$ are
$$
\begin{array}{ccccc}
    
    1 & 1 & \cdots & 1 & 2  \\
    
    2 & 2 & \cdots & 2 & 3  \\

\end{array}
\qquad
\textrm{ and }
\qquad
\begin{array}{ccccc}
    
    1 & 1 & \cdots & 1 & 2  \\
    
    2 & 2 & \cdots & 2 & \bar{2}  \\

\end{array}
.
$$
In the former case, the operators $f_{2},\ldots,f_{n}$ act freely on
the bottom row and the righmost box in the top row; in the latter
case, they act on all boxes in the bottom row except the rightmost,
since $\ftil{i}$, $\etil{i}$ for $i=2, \ldots, n$ do not act on the
last column.  It follows that the corresponding partitions are $(s,1)$
and $(s-1)$, with respective $\pm$ diagrams
$$
\begin{array}{|c|c|c|c|c|}
    \hline
     &  & \cdots &  &   \\
    \hline
    
     \phantom{+}& + & \cdots & + & +  \\
    \hline

\end{array}
\qquad
\textrm{ and }
\qquad
\begin{array}{|c|c|c|c|c|}
    \hline
    
     &  & \cdots &  & +  \\
    \hline
    
    + & + & \cdots & + & -  \\
    \hline

\end{array}
.
$$

Now, observe what can result from acting on a tableau
$T=\stack{a_1}{b_1}\cdots\stack{a_s}{b_s}$ in the highest weight
branching component vertex by $\ftil{1}$.  Since
$a_{1}=\cdots=a_{s}=1$, $a_{s}$ will turn into a $2$.  There are two
cases to consider: if $b_{s}=\bar{2}$, this results in a tableau with
a configuration $\stack{2}{\bar{2}}$ at the right end; otherwise, it
is a tableau with $a_{1}=\cdots=a_{s-1}=1$ where some element of
$D_{n-1}$ can act on the rightmost column.  It follows that there are
edges from the vertex in stratum $s$ to both of the vertices in
stratum $s-1$.

\begin{proof}[Proof of Proposition \ref{prop:recs}]

    The number of each partition and their strata follow from
    Algorithm \ref{alg:pmbranch} and the $D_{n}$ to $D_{n-1}$
    branching rule.  Explicitly, consider the content of stratum $j$,
    which is determined by the set of $\pm$ diagrams such that the
    difference between the number of $+$'s and the number of $-$'s is
    $j$.  Thanks to Proposition \ref{prop:BCsym}, we may assume $j\geq
    0$.  Because we currently only consider rectangular partitions,
    the possible configurations of $+$'s and $-$'s can be explicitly
    characterized.  First, we have $\ell$ columns with only a $+$ in
    the last row with $j\leq \ell \leq s$.  We furthermore have
    $\ell-j$ columns with only a $-$ in the last row, located
    immediately to the right of the columns with only a $+$.  Finally,
    we may have between $0$ and $s-2\ell+j$ columns with a $\pm$ pair
    in the last two rows, located immediately to the right of the
    columns with only a $-$ and extending to the right edge of the
    rectangle.  The choices of the number of $+$'s and $\pm$ columns
    are independent of each other up to the bounds described above.
    In other words, stratum $j$ has one vertex for each partition
    contained in a $k\times s$ rectangle such that the first $k-2$
    parts are $s$, the last part is no greater than $s-j$, and the
    difference between the size of the partition and $k\times s$ is
    congruent mod $2$ to $j$.  This is equivalent to the
    characterization of vertices prescribed by Proposition
    \ref{prop:recs}; i.e., those partitions contained in a $k\times s$
    rectangle such that the first $k-2$ parts are $s$ and which differ
    from a $(k-1)\times s$ rectangle by $s-j$ additions or removals of
    boxes in the last two rows.

    To show that the edges of this graph are as prescribed by
    Proposition \ref{prop:recs}, we explicitly construct tableaux with
    the following properties.  Each of these tableaux is in a
    $D_{n-1}$ crystal whose highest weight corresponds to a partition
    that can have a box added in at least one of the last two rows.
    We will then show that acting by the Kashiwara operator $f_{1}$
    produces a tableau in a $D_{n-1}$ crystal whose highest weight
    corresponds to the partition with one of the last two rows
    augmented thusly.

    The claims that $f_{1}$ may also have the effect of removing a box
    from one of the last two rows of the partition labeling its
    $D_{n-1}$ crystal's highest weight follows by the $*$-duality of
    the crystal.  To see this, let $v$ be a branching component vertex
    labeled by a partition $\mu$ such that $\mu_{k} \geq 1$ (and
    consequently $j$, the stratum of $v$, is at least $-s + 1$); we
    wish to show that there is an edge from $v$ to $w$, where $w$ is
    in stratum $j-1$ and its partition differs from $\mu$ only by
    diminishing the last part by one.  We know that there is a vertex
    $v'$ also labeled by $\mu$ in stratum $-j$ called the
    \label{comp_vert}
    complementary vertex of $v$; the vertices in $B(v')$ are precisely
    the image under the $*$ map of the vertices in $B(v)$.  We know
    from the above paragraph that $v'$ has an incoming edge from $w'$,
    the complementary vertex of $w$ .  This means that for some
    crystal vertex $b\in B(w')$, we know that $f_{1}b\in B(v')$.  It
    follows that $e_{1}(b^{*})\in B(v)$ and $b^{*}\in B(w)$, and we
    thus conclude that there is an edge from $v$ to $w$.  An identical
    argument shows that a branching component vertex $v$ in stratum
    $j\geq -s+1$ labeled by a partition $\mu$ such that $\mu_{k-1}
    \geq \mu_{k}+1$ has an edge to the vertex $w$ in stratum $j-1$
    whose partition differs from $\mu$ only by diminishing the
    penultimate part by one.
    
    We now construct tableaux such that the action of the Kashiwara
    operator $f_{1}$ has the described effect on the branching
    component vertices.  We begin by constructing an arbitrary
    $D_{n-1}$ highest weight tableau, as given by Algorithm
    \ref{alg:pmbranch}.  We have two cases to consider.  One is
    $$
    \begin{array}{cccccccccccc}
        
        1 & \cdots & 1 & 1 & \cdots & 1 & 2 & \cdots & 2 & 2 & \cdots & 
	2  \\
        
        2 & \cdots & 2 & 2 & \cdots & 2 & 3 & \cdots & 3 & 3 & \cdots & 
	3  \\
        
        \vdots & \cdots & \vdots & \vdots & \cdots & \vdots & \vdots & \cdots & \vdots & \vdots & \cdots & \vdots  \\
        
        k-1 & \cdots & k-1 & k-1 & \cdots & k-1 & k & \cdots & k & k & \cdots & 
	k  \\
        
        k & \cdots & k & k+1 & \cdots & k+1 & k+1 & \cdots & k+1 &
	\bar{1} & \cdots & \bar{1}  \\
        
	\multicolumn{3}{c}{\underbrace{\phantom{eeeeeeeeeeeeeeee}}_{s_{1}}} &
	\multicolumn{3}{c}{\underbrace{\phantom{eeeeeeeeeeeeeeee}}_{s_{2}}} &
	\multicolumn{3}{c}{\underbrace{\phantom{eeeeeeeeeeeeeeee}}_{s_{3}}} &
	\multicolumn{3}{c}{\underbrace{\phantom{eeeeeeeee}}_{s_{4}}}
         
    \end{array}
    $$
    with $D_{n-1}$ weight $s_{2}\varpi_{k-2} +
    (s_{1}+s_{4}-s_{2})\varpi_{k-1} + (s_{2}+s_{3})\varpi_{k}$.  This
    tableau corresponds by Algorithm \ref{alg:pmbranch} to the $\pm$
    diagram with $s_{1}$ columns with a $+$, $s_{4}-s_{2}$ columns
    with a $-$, and $s_{2}$ columns with a $\pm$ pair.  Note that this
    implies $s_{2}\leq s_{4}$.  The other case is
    $$
    \begin{array}{ccccccccccccc}
        
        1 & \cdots & 1 & 1 & \cdots & 1 & 1 & \cdots & 1 & 2 & 2 & \cdots & 
	2  \\
        
        2 & \cdots & 2 & 2 & \cdots & 2 & 2 & \cdots & 2 & 3 & 3 & \cdots & 
	3  \\
        
	\vdots & \cdots & \vdots & \vdots & \cdots & \vdots & \vdots &
	\vdots & \cdots & \vdots & \vdots & \cdots & \vdots \\
        
	k-1 & \cdots & k-1 & k-1 & \cdots & k-1 & k-1 & \cdots & k-1 &
	k & k & \cdots & k \\
        
        k & \cdots & k & k+1 & \cdots & k+1 & \bar{k} & \cdots & \bar{k} & 
	\bar{k} &
	\bar{1} & \cdots & \bar{1}  \\
        
	\multicolumn{3}{c}{\underbrace{\phantom{eeeeeeeeeeeeeeee}}_{s_{1}}} &
	\multicolumn{3}{c}{\underbrace{\phantom{eeeeeeeeeeeeeeee}}_{s_{2}}} &
	\multicolumn{3}{c}{\underbrace{\phantom{eeeeeeeeeeeeeeee}}_{s_{3}}} &
	\raisebox{-6pt}{$\delta$} &
	\multicolumn{3}{c}{\underbrace{\phantom{eeeeeeeee}}_{s_{4}}}
         
    \end{array}
    $$
    where $\delta$ may be $0$ or $1$.  This tableau has $D_{n-1}$
    weight $(s_{2}+2s_{3}+\delta)\varpi_{k-2} +
    (s_{1}+s_{4}-s_{2}-s_{3})\varpi_{k-1} + s_{2}\varpi_{k}$,
    corresponding by Algorithm \ref{alg:pmbranch} to the $\pm$
    diagram with $s_{1}$ columns with a $+$, $s_{4}-s_{2}-s_{3}$
    columns with a $-$, and $s_{2}$ columns with a $\pm$ pair.  In
    this case we must have $s_{2}+s_{3}\leq s_{4}$.  To make the
    second case exclusive of the first, we add the requirement that
    $\delta+s_{3}>0$.

    Consider the action on these tableaux of the sequence
    $f_{1}f_{2}\cdots f_{k-1}$.  In the first case, we assume
    $s_{2}\neq 0$, otherwise the $k-1^{\mathrm{st}}$ part of the
    associated partition has size $s$.  We therefore know that this
    sequence acts entirely on the $(s_{1}+s_{2})^{\mathrm{th}}$
    column, counting from the left.  It has the effect of increasing
    by 1 each entry in this column except the bottom entry, resulting
    in the
    $D_{n-1}$ highest weight tableau with $D_{n-1}$ weight
    $(s_{2}-1)\varpi_{k-2} + (s_{1}+s_{4}-s_{2}+1)\varpi_{k-1} +
    (s_{2}+s_{3})\varpi_{k}$.
    
    In the second case, we consider the subcases of $\delta = 0$ and
    $\delta = 1$.  If $\delta = 0$, we have $s_{3}\neq 0$, and the
    sequence $f_{1}f_{2}\cdots f_{k-1}$ acts entirely on the
    $(s_{1}+s_{2}+s_{3})^{\mathrm{th}}$ column, counting from the
    left, resulting in the 
    $D_{n-1}$ highest weight tableau with $D_{n-1}$ weight
    $(s_{2}+2s_{3}-1)\varpi_{k-2} +
    (s_{1}+s_{4}-s_{2}-s_{3}+1)\varpi_{k-1} + s_{2}\varpi_{k}$.  On
    the other hand,  if $\delta = 1$,  this sequence acts entirely on
    the letter in the bottom row of the ``$\delta$'' column,  giving
    us the 
    $D_{n-1}$ highest weight tableau with $D_{n-1}$ weight
    $(s_{2}+2s_{3})\varpi_{k-2} +
    (s_{1}+s_{4}-s_{2}-s_{3}+1)\varpi_{k-1} + s_{2}\varpi_{k}$.
    
    In any case,  we have constructed a tableau (namely, 
    $f_{2}f_{3}\cdots f_{k-1}$ applied to a $D_{n-1}$ highest weight
    tableau) such that acting by $f_{1}$ adds a box to the penultimate
    row of the partition associated to its branching component vertex,
    according to Algorithm \ref{alg:pmbranch}.

    Moving our attention to adding a box in the last row, consider the
    tableau that results from applying the sequence of operators
    $$f_{2}^{s_{1}+s_{4}-s_{2}}f_{3}^{s_{1}+s_{4}-s_{2}}
    \cdots f_{k}^{s_{1}+s_{4}-s_{2}}$$
    to
    our first case of highest weight tableau.  This produces
    
    \footnotesize{
    $$
    \begin{array}{ccccccccccccccc}
        
	1 & \cdots & 1 & 1 & \cdots & 1 & 2 & \cdots & 2 & 2 & \cdots
	& 2 & 3 & \cdots & 3 \\
        
	2 & \cdots & 2 & 3 & \cdots & 3 & 3 & \cdots & 3 & 3 & \cdots
	& 3 & 4 & \cdots & 4 \\
        
	\vdots & \cdots & \vdots & \vdots & \cdots & \vdots & \vdots &
	\cdots & \vdots & \vdots & \cdots & \vdots & \vdots & \cdots &
	\vdots \\
        
	k-1 & \cdots & k-1 & k & \cdots & k & k & \cdots & k & k &
	\cdots & k & k+1 & \cdots & k+1 \\
        
        k+1 & \cdots & k+1 & k+1 & \cdots & k+1 & k+1 & \cdots & k+1 &
	\bar{1} & \cdots & \bar{1} & \bar{1} & \cdots & \bar{1}  \\
        
	\multicolumn{3}{c}{\underbrace{\phantom{eeeeeeeeeeeeeeee}}_{s_{2}}} &
	\multicolumn{3}{c}{\underbrace{\phantom{eeeeeeeeeeeeeeee}}_{s_{1}}} &
	\multicolumn{3}{c}{\underbrace{\phantom{eeeeeeeeeeeeeeee}}_{s_{3}}} &
	\multicolumn{3}{c}{\underbrace{\phantom{eeeeeeeee}}_{s_{2}}} &
	\multicolumn{3}{c}{\underbrace{\phantom{eeeeeeeeeeeeeeee}}_{s_{4}-s_{2}}}
         
    \end{array}
    $$    
    }
    \normalsize{    Applying $f_{1}$ to this tableau gives us 
    }
    \footnotesize{
    $$
    \begin{array}{ccccccccccccccc}
        
	1 & \cdots & 1 & 1 & \cdots & 1 & 2 & \cdots & 2 & 2 & \cdots
	& 2 & 3 & \cdots & 3 \\
        
	2 & \cdots & 2 & 3 & \cdots & 3 & 3 & \cdots & 3 & 3 & \cdots
	& 3 & 4 & \cdots & 4 \\
        
	\vdots & \cdots & \vdots & \vdots & \cdots & \vdots & \vdots &
	\cdots & \vdots & \vdots & \cdots & \vdots & \vdots & \cdots &
	\vdots \\
        
	k-1 & \cdots & k-1 & k & \cdots & k & k & \cdots & k & k &
	\cdots & k & k+1 & \cdots & k+1 \\
        
        k+1 & \cdots & k+1 & k+1 & \cdots & k+1 & k+1 & \cdots & k+1 &
	\bar{1} & \cdots & \bar{1} & \bar{1} & \cdots & \bar{1}  \\
        
	\multicolumn{3}{c}{\underbrace{\phantom{eeeeeeeeeeeeeeee}}_{s_{2}}} &
	\multicolumn{3}{c}{\underbrace{\phantom{eeeeeeeeeeeeeeee}}_{s_{1}-1}} &
	\multicolumn{3}{c}{\underbrace{\phantom{eeeeeeeeeeeeeeee}}_{s_{3}+1}} &
	\multicolumn{3}{c}{\underbrace{\phantom{eeeeeeeee}}_{s_{2}}} &
	\multicolumn{3}{c}{\underbrace{\phantom{eeeeeeeeeeeeeeee}}_{s_{4}-s_{2}}}
         
    \end{array}
    $$    
    }
    
    \normalsize{ and on this tableau the sequence
    $$e_{k}^{s_{1}+s_{4}-s_{2}-1}e_{k-1}^{s_{1}+s_{4}-s_{2}-1}\cdots
    e_{2}^{s_{1}+s_{4}-s_{2}-1}$$ produces a $D_{n-1}$ highest weight
    tableau with $D_{n-1}$ weight $s_{2}\varpi_{k-2} +
    (s_{1}+s_{4}-s_{2}-1)\varpi_{k-1} + (s_{2}+s_{3}+1)\varpi_{k}$,}
    corresponding to adding a box to the last row of the associated
    partition.
    
    We now consider the sequence $$F =
    f_{2}^{s_{1}+s_{4}-s_{2}-s_{3}}
    f_{3}^{s_{1}+s_{4}-s_{2}-s_{3}} \cdots
    f_{k}^{s_{1}+s_{4}-s_{2}-s_{3}}$$ applied to our second case of
    highest weight tableau.  We consider separately the cases of
    whether $s_{1}$ is greater than or less than $2s_{3}+\delta$.
    
    If $s_{1} \leq 2s_{3}+\delta$,  the sequence $F$ applied to our
    highest weight tableau produces
    \scriptsize{
    $$
    \begin{array}{ccccccccccccccccc}
        
	1 & \cdots & 1 & 1 & \cdots & 1 & 1 & 1 & \cdots
	& 1 & 2 & 2 & \cdots & 2 & 3 & \cdots & 3 \\
        
	2 & \cdots & 2 & 2 & \cdots & 2 & 3 & 3 & \cdots
	& 3 & 3 & 3 & \cdots & 3 & 4 & \cdots & 4\\
        
	\vdots & \cdots & \vdots & \vdots & \cdots & \vdots & \vdots &
	\vdots & \cdots & \vdots & \vdots & \vdots &
	\cdots & \vdots & \vdots & \cdots &
	\vdots\\
        
	k-1 & \cdots & k-1 & k-1 & \cdots & k-1 & k & k &
	\cdots & k & k & k & \cdots & k & k+1 & \cdots & k+1\\
        
	k+1 & \cdots & k+1 & \bar{k} & \cdots & \bar{k} & \bar{k}  &
	\bar{2} & \cdots & \bar{2} & \bar{2} & \bar{1} & \cdots & \bar{1} &
	\bar{1} & \cdots & \bar{1} \\
        
	\multicolumn{3}{c}{\underbrace{\phantom{eeeeeeeeeeeeeeee}}_{s_{1}+s_{2}}} &
	\multicolumn{3}{c}{\underbrace{\phantom{eeeeeeeeeeeeeeee}}_{s_{3}-\Big\lfloor\frac{s_{1}-\delta}{2}\Big\rfloor}} &
	\raisebox{-6pt}{$(s_{1}-\delta)_{\mathrm{mod}2}$} &
	\multicolumn{3}{c}{\underbrace{\phantom{eeeeeeeee}}_{\Big\lfloor\frac{s_{1}-\delta}{2}\Big\rfloor}} &
	\raisebox{-6pt}{$\delta$} &
	\multicolumn{3}{c}{\underbrace{\phantom{eeeeeeeee}}_{s_{2}+s_{3}}} &
	\multicolumn{3}{c}{\underbrace{\phantom{eeeeeeeeeeeeeeee}}_{s_{4}-s_{2}-s_{3}}}
         
    \end{array}
    $$        

    }
    
    \normalsize{and if $s_{1} > 2s_{3}+\delta$,  we have}
    
    \footnotesize{
    $$
    \begin{array}{cccccccccccccccc}
        
	1 & \cdots & 1 & 1 & \cdots & 1 & 1 & \cdots
	& 1 & 2 & 2 & \cdots & 2 & 3 & \cdots & 3 \\
        
	2 & \cdots & 2 & 3 & \cdots & 3 & 3 & \cdots
	& 3 & 3 & 3 & \cdots & 3 & 4 & \cdots & 4\\
        
	\vdots & \cdots & \vdots & \vdots & \cdots & \vdots & \vdots &
	\cdots & \vdots & \vdots & \vdots & \cdots & \vdots & \vdots &
	\cdots & \vdots\\
        
	k-1 & \cdots & k-1 & k & \cdots & k & k &
	\cdots & k & k & k & \cdots & k & k+1 & \cdots & k+1\\
        
	k+1 & \cdots & k+1 & k+1 & \cdots & k+1 &
	\bar{2} & \cdots & \bar{2} & \bar{2} & \bar{1} & \cdots & \bar{1} &
	\bar{1} & \cdots & \bar{1} \\
        
	\multicolumn{3}{c}{\underbrace{\phantom{eeeeeeeeeeeeeeee}}_{s_{2}+2s_{3}+\delta}} &
	\multicolumn{3}{c}{\underbrace{\phantom{eeeeeeeeeeeeeeee}}_{s_{1}-(2s_{3}-\delta)}} &
	\multicolumn{3}{c}{\underbrace{\phantom{eeeeeeeee}}_{s_{3}}} &
	\raisebox{-6pt}{$\delta$} &
	\multicolumn{3}{c}{\underbrace{\phantom{eeeeeeeee}}_{s_{2}+s_{3}}} &
	\multicolumn{3}{c}{\underbrace{\phantom{eeeeeeeeeeeeeeee}}_{s_{4}-s_{2}-s_{3}}}
         
    \end{array}
    $$        

    }
    
    \normalsize{ } We treat the subcases of $\delta = 0$ and $\delta =
    1$ separately.  In both of the above tableaux, if $\delta = 1$,
    the Kashiwara operator $f_{1}$ acts on the $\bar{2}$ at the bottom
    of the $\delta$ column, and if $\delta = 0$, it acts on the $1$ at
    the top of the column immediately to the left of where the
    $\delta$ column would be.

    \normalsize{If $s_{1} \leq 2s_{3}+\delta$ and $\delta = 0$ we get}
    \scriptsize{
    $$
    \begin{array}{ccccccccccccccccc}
        
	1 & \cdots & 1 & 1 & \cdots & 1 & 1 & 1 & \cdots
	& 1 & 2 & 2 & \cdots & 2 & 3 & \cdots & 3 \\
        
	2 & \cdots & 2 & 2 & \cdots & 2 & 3 & 3 & \cdots
	& 3 & 3 & 3 & \cdots & 3 & 4 & \cdots & 4\\
        
	\vdots & \cdots & \vdots & \vdots & \cdots & \vdots & \vdots &
	\vdots & \cdots & \vdots & \vdots & \vdots &
	\cdots & \vdots & \vdots & \cdots &
	\vdots\\
        
	k-1 & \cdots & k-1 & k-1 & \cdots & k-1 & k & k &
	\cdots & k & k & k & \cdots & k & k+1 & \cdots & k+1\\
        
	k+1 & \cdots & k+1 & \bar{k} & \cdots & \bar{k} & \bar{k}  &
	\bar{2} & \cdots & \bar{2} & \bar{2} & \bar{1} & \cdots & \bar{1} &
	\bar{1} & \cdots & \bar{1} \\
        
	\multicolumn{3}{c}{\underbrace{\phantom{eeeeeeeeeeeeeeee}}_{s_{1}+s_{2}}} &
	\multicolumn{3}{c}{\underbrace{\phantom{eeeeeeeeeeeeeeee}}_{s_{3}-\Big\lfloor\frac{s_{1}-\delta}{2}\Big\rfloor}} &
	\raisebox{-6pt}{$(s_{1}-\delta)_{\mathrm{mod}2}$} &
	\multicolumn{3}{c}{\underbrace{\phantom{eeeeeeeee}}_{\Big\lfloor\frac{s_{1}-\delta}{2}\Big\rfloor-1}} &
	\raisebox{-6pt}{1} &
	\multicolumn{3}{c}{\underbrace{\phantom{eeeeeeeee}}_{s_{2}+s_{3}}} &
	\multicolumn{3}{c}{\underbrace{\phantom{eeeeeeeeeeeeeeee}}_{s_{4}-s_{2}-s_{3}}}
         
    \end{array}
    $$            
    }
    If $s_{1} \leq 2s_{3}+\delta$ and $\delta = 1$ we get
    \scriptsize{
    $$
    \begin{array}{cccccccccccccccc}
        
	1 & \cdots & 1 & 1 & \cdots & 1 & 1 & 1 & \cdots
	& 1 & 2 & \cdots & 2 & 3 & \cdots & 3 \\
        
	2 & \cdots & 2 & 2 & \cdots & 2 & 3 & 3 & \cdots
	& 3 & 3 & \cdots & 3 & 4 & \cdots & 4\\
        
	\vdots & \cdots & \vdots & \vdots & \cdots & \vdots & \vdots &
	\vdots & \cdots & \vdots & \vdots & \cdots &
	\vdots & \vdots & \cdots &
	\vdots\\
        
	k-1 & \cdots & k-1 & k-1 & \cdots & k-1 & k & k &
	\cdots & k & k & \cdots & k & k+1 & \cdots & k+1\\
        
	k+1 & \cdots & k+1 & \bar{k} & \cdots & \bar{k} & \bar{k}  &
	\bar{2} & \cdots & \bar{2} & \bar{1} & \cdots & \bar{1} &
	\bar{1} & \cdots & \bar{1} \\
        
	\multicolumn{3}{c}{\underbrace{\phantom{eeeeeeeeeeeeeeee}}_{s_{1}+s_{2}}} &
	\multicolumn{3}{c}{\underbrace{\phantom{eeeeeeeeeeeeeeee}}_{s_{3}-\Big\lfloor\frac{s_{1}-\delta}{2}\Big\rfloor}} &
	\raisebox{-6pt}{$(s_{1}-\delta)_{\mathrm{mod}2}$} &
	\multicolumn{3}{c}{\underbrace{\phantom{eeeeeeeee}}_{\Big\lfloor\frac{s_{1}-\delta}{2}\Big\rfloor}} &
	\multicolumn{3}{c}{\underbrace{\phantom{eeeeeeeee}}_{s_{2}+s_{3}+1}} &
	\multicolumn{3}{c}{\underbrace{\phantom{eeeeeeeeeeeeeeee}}_{s_{4}-s_{2}-s_{3}}}
         
    \end{array}
    $$      
    }

    \normalsize{On the other hand, if $s_{1} > 2s_{3}+\delta$ and
    $\delta = 0$, we have}
    
    \footnotesize{
    $$
    \begin{array}{cccccccccccccccc}
        
	1 & \cdots & 1 & 1 & \cdots & 1 & 1 & \cdots
	& 1 & 2 & 2 & \cdots & 2 & 3 & \cdots & 3 \\
        
	2 & \cdots & 2 & 3 & \cdots & 3 & 3 & \cdots
	& 3 & 3 & 3 & \cdots & 3 & 4 & \cdots & 4\\
        
	\vdots & \cdots & \vdots & \vdots & \cdots & \vdots & \vdots &
	\cdots & \vdots & \vdots & \vdots & \cdots & \vdots & \vdots &
	\cdots & \vdots\\
        
	k-1 & \cdots & k-1 & k & \cdots & k & k &
	\cdots & k & k & k & \cdots & k & k+1 & \cdots & k+1\\
        
	k+1 & \cdots & k+1 & k+1 & \cdots & k+1 &
	\bar{2} & \cdots & \bar{2} & \bar{2} & \bar{1} & \cdots & \bar{1} &
	\bar{1} & \cdots & \bar{1} \\
        
	\multicolumn{3}{c}{\underbrace{\phantom{eeeeeeeeeeeeeeee}}_{s_{2}+2s_{3}+\delta}} &
	\multicolumn{3}{c}{\underbrace{\phantom{eeeeeeeeeeeeeeee}}_{s_{1}-(2s_{3}-\delta)}} &
	\multicolumn{3}{c}{\underbrace{\phantom{eeeeeeeee}}_{s_{3}-1}} &
	\raisebox{-6pt}{1} &
	\multicolumn{3}{c}{\underbrace{\phantom{eeeeeeeee}}_{s_{2}+s_{3}}} &
	\multicolumn{3}{c}{\underbrace{\phantom{eeeeeeeeeeeeeeee}}_{s_{4}-s_{2}-s_{3}}}
         
    \end{array}
    $$        

    } \normalsize{and if $s_{1} > 2s_{3}+\delta$ and $\delta = 1$, we
    have}
    
    \footnotesize{
    $$
    \begin{array}{ccccccccccccccc}
        
	1 & \cdots & 1 & 1 & \cdots & 1 & 1 & \cdots
	& 1 & 2 & \cdots & 2 & 3 & \cdots & 3 \\
        
	2 & \cdots & 2 & 3 & \cdots & 3 & 3 & \cdots
	& 3 & 3 & \cdots & 3 & 4 & \cdots & 4\\
        
	\vdots & \cdots & \vdots & \vdots & \cdots & \vdots & \vdots &
	\cdots & \vdots & \vdots & \cdots & \vdots & \vdots &
	\cdots & \vdots\\
        
	k-1 & \cdots & k-1 & k & \cdots & k & k &
	\cdots & k & k & \cdots & k & k+1 & \cdots & k+1\\
        
	k+1 & \cdots & k+1 & k+1 & \cdots & k+1 &
	\bar{2} & \cdots & \bar{2} & \bar{1} & \cdots & \bar{1} &
	\bar{1} & \cdots & \bar{1} \\
        
	\multicolumn{3}{c}{\underbrace{\phantom{eeeeeeeeeeeeeeee}}_{s_{2}+2s_{3}+\delta}} &
	\multicolumn{3}{c}{\underbrace{\phantom{eeeeeeeeeeeeeeee}}_{s_{1}-(2s_{3}-\delta)}} &
	\multicolumn{3}{c}{\underbrace{\phantom{eeeeeeeee}}_{s_{3}}} &
	\multicolumn{3}{c}{\underbrace{\phantom{eeeeeeeee}}_{s_{2}+s_{3}+1}} &
	\multicolumn{3}{c}{\underbrace{\phantom{eeeeeeeeeeeeeeee}}_{s_{4}-s_{2}-s_{3}}}
         
    \end{array}
    $$        

    }

    \normalsize{   }     
    
    In all of these cases,  we find that applying the sequence $$E =
    e_{k}^{(s_{1}+s_{4}-s_{2}-s_{3}-1)}
    e_{3}^{(s_{1}+s_{4}-s_{2}-s_{3}-1)} \cdots
    e_{2}^{(s_{1}+s_{4}-s_{2}-s_{3}-1)}$$ has the effect of producing
    the tableau
    
    $$
    \begin{array}{ccccccccccccc}
        
        1 & \cdots & 1 & 1 & \cdots & 1 & 1 & \cdots & 1 & 2 & 2 & \cdots & 
	2  \\
        
        2 & \cdots & 2 & 2 & \cdots & 2 & 2 & \cdots & 2 & 3 & 3 & \cdots & 
	3  \\
        
	\vdots & \cdots & \vdots & \vdots & \cdots & \vdots & \vdots &
	\vdots & \cdots & \vdots & \vdots & \cdots & \vdots \\
        
	k-1 & \cdots & k-1 & k-1 & \cdots & k-1 & k-1 & \cdots & k-1 &
	k & k & \cdots & k \\
        
        k & \cdots & k & k+1 & \cdots & k+1 & \bar{k} & \cdots & \bar{k} & 
	\bar{k} &
	\bar{1} & \cdots & \bar{1}  \\
        
	\multicolumn{3}{c}{\underbrace{\phantom{eeeeeeeeeeeeeeee}}_{s_{1}-1}} &
	\multicolumn{3}{c}{\underbrace{\phantom{eeeeeeeeeeeeeeee}}_{s_{2}+1}} &
	\multicolumn{3}{c}{\underbrace{\phantom{eeeeeeeeeeeeeeee}}_{s_{3}-\tilde{\delta}}} &
	\raisebox{-6pt}{$1-\delta$} &
	\multicolumn{3}{c}{\underbrace{\phantom{eeeeeeeee}}_{s_{4}+\delta}}
         
    \end{array}
    $$
    This tableau is a $D_{n-1}$ highest weight tableau with $D_{n-1}$
    weight $(s_{2}+2s_{3}+\delta)\varpi_{k-2} +
    (s_{1}+s_{4}-s_{2}-s_{3}-1)\varpi_{k-1} + (s_{2}+1)\varpi_{k}$,
    which precisely corresponds to adding a box to the last row of the
    partition we started with.
\end{proof}

\begin{figure}
\begin{center}\setlength{\unitlength}{1in}
    \begin{picture}(0, .3)
	\put(0, 0){$\emptyset$}
	\put(.03, -.02){\vector(0, -1){.14}}
	\put(.09, -.18){\line(0, -1){.125}}
	\put(.09, -.18){\line(-1, 0){.125}}	
	\put(-.035, -.18){\line(0, -1){.125}}
	\put(-.035, -.305){\line(1, 0){.125}}	
	\put(-.035, -.32){\vector(-1, -2){.1}}
	\put(.09, -.32){\vector(1, -2){.1}}	
	\put(-.3, -.55){\line(1, 0){.25}}
	\put(-.3, -.55){\line(0, -1){.125}}
	\put(-.3, -.675){\line(1, 0){.25}}
	\put(-.175, -.55){\line(0, -1){.125}}	
	\put(-.05, -.55){\line(0, -1){.125}}		
	\put(.2, -.66){$\emptyset$}
	\put(-.175, -.69){\vector(0, -1){.14}}
	\put(.23, -.69){\vector(0, -1){.14}}
	\put(-.05, -.69){\vector(1, -1){.14}}
	\put(.29, -.85){\line(0, -1){.125}}
	\put(.29, -.85){\line(-1, 0){.125}}	
	\put(.165, -.85){\line(0, -1){.125}}
	\put(.165, -.975){\line(1, 0){.125}}
	\put(-.345, -.85){\line(1, 0){.375}}
	\put(-.345, -.85){\line(0, -1){.125}}
	\put(-.095, -.85){\line(0, -1){.125}}
	\put(.03, -.85){\line(0, -1){.125}}
	\put(-.345, -.975){\line(1, 0){.375}}
	\put(-.22, -.85){\line(0, -1){.125}}
	\put(-.175, -.99){\vector(0, -1){.14}}
	\put(.23, -.99){\vector(0, -1){.14}}
	\put(.09, -.99){\vector(-1, -1){.14}}
	\put(0, -1.3){$\vdots$}
	
    \end{picture}
\end{center}

\vspace{3cm}
\caption{Strata $3$ through $0$ of the branching component graph 
$\BC((3))$\label{fig:onerow}}
\end{figure}
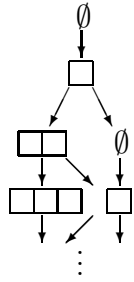

\begin{example}
    Figure \ref{fig:onerow} depicts the non-negative strata of $\BC((3))$.
\end{example}

\begin{figure}
\begin{center}\setlength{\unitlength}{1in}
    \begin{picture}(-2.2, 0)
	\put(-2.5, 0){\line(1, 0){.25}}
	\put(-2.5, 0){\line(0, -1){.125}}
	\put(-2.25, 0){\line(0, -1){.125}}
	\put(-2.5, -.125){\line(1, 0){.25}}
	\put(-2.375, 0){\line(0, -1){.125}}
	\put(-2.45, -.16){\vector(-1, -2){.1}}
	\put(-2.3, -.16){\vector(1, -2){.12}}
	\put(-2.75, -.38){\line(1, 0){.25}}
	\put(-2.75, -.38){\line(0, -1){.25}}
	\put(-2.625, -.38){\line(0, -1){.25}}
	\put(-2.75, -.505){\line(1, 0){.25}}
	\put(-2.5, -.38){\line(0, -1){.125}}
	\put(-2.75, -.63){\line(1, 0){.125}}
	\put(-2.85, -.68){\vector(-1, -2){.1}}
	\put(-2.65, -.68){\vector(0, -1){.2}}
	\put(-2.52, -.63){\vector(1, -1){.25}}
	\put(-2.25, -.44){\line(1, 0){.125}}
	\put(-2.25, -.44){\line(0, -1){.125}}
	\put(-2.125, -.44){\line(0, -1){.125}}
	\put(-2.25, -.565){\line(1, 0){.125}}
	\put(-2.09, -.68){\vector(1, -2){.1}}
	\put(-2.2, -.68){\vector(0, -1){.2}}
	\put(-2.35, -.63){\vector(-1, -1){.25}}
	\put(-3.1, -.9){\line(1, 0){.25}}
	\put(-3.1, -.9){\line(0, -1){.25}}
	\put(-2.975, -.9){\line(0, -1){.25}}
	\put(-3.1, -1.025){\line(1, 0){.25}}
	\put(-2.85, -.9){\line(0, -1){.25}}
	\put(-3.1, -1.15){\line(1, 0){.25}}
	\put(-2.95, -1.2){\vector(1, -2){.1}}
	\put(-2.75, -.96){\line(1, 0){.25}}
	\put(-2.75, -.96){\line(0, -1){.125}}
	\put(-2.5, -.96){\line(0, -1){.125}}
	\put(-2.75, -1.085){\line(1, 0){.25}}
	\put(-2.625, -.96){\line(0, -1){.125}}
	\put(-2.65, -1.2){\vector(0, -1){.2}}
	\put(-2.55, -1.2){\vector(1, -1){.25}}
	\put(-2.3, -.9){\line(1, 0){.125}}
	\put(-2.3, -.9){\line(0, -1){.25}}
	\put(-2.175, -.9){\line(0, -1){.25}}
	\put(-2.3, -1.025){\line(1, 0){.125}}
	\put(-2.3, -1.15){\line(1, 0){.125}}
	\put(-2.27, -1.2){\vector(-1, -1){.21}}
	\put(-2.2, -1.2){\vector(0, -1){.2}}
	\put(-2, -1.08){$\emptyset$}
	\put(-2, -1.2){\vector(-1, -2){.1}}
	\put(-2.75, -1.42){\line(1, 0){.25}}
	\put(-2.75, -1.42){\line(0, -1){.25}}
	\put(-2.625, -1.42){\line(0, -1){.25}}
	\put(-2.75, -1.545){\line(1, 0){.25}}
	\put(-2.5, -1.42){\line(0, -1){.125}}
	\put(-2.75, -1.67){\line(1, 0){.125}}
	\put(-2.55, -1.7){\vector(1, -2){.1}}
	\put(-2.25, -1.48){\line(1, 0){.125}}
	\put(-2.25, -1.48){\line(0, -1){.125}}
	\put(-2.125, -1.48){\line(0, -1){.125}}
	\put(-2.25, -1.608){\line(1, 0){.125}}
	\put(-2.2, -1.66){\vector(-1, -2){.12}}	
	\put(-2.5, -1.95){\line(1, 0){.25}}
	\put(-2.5, -1.95){\line(0, -1){.125}}
	\put(-2.25, -1.95){\line(0, -1){.125}}
	\put(-2.5, -2.075){\line(1, 0){.25}}
	\put(-2.375, -1.95){\line(0, -1){.125}}
    \end{picture}
\end{center}
\begin{center}\setlength{\unitlength}{1in}
    \begin{picture}(-5.1, .31)
	\put(-2.25, .06){\line(1, 0){.125}}
	\put(-2.25, .06){\line(0, -1){.125}}
	\put(-2.125, .06){\line(0, -1){.125}}
	\put(-2.25, -.065){\line(1, 0){.125}}
	\put(-2.125, -.12){\vector(1, -2){.12}}
	\put(-2.25, -.12){\vector(-1, -2){.1}}
	\put(-2.42, -.4){\line(1, 0){.125}}
	\put(-2.42, -.4){\line(0, -1){.25}}
	\put(-2.295, -.4){\line(0, -1){.25}}
	\put(-2.42, -.525){\line(1, 0){.125}}
	\put(-2.42, -.65){\line(1, 0){.125}}
	\put(-2.35, -.73){\vector(1, -2){.1}}
	\put(-2, -.58){$\emptyset$}
	\put(-2, -.68){\vector(-1, -2){.127}}	
	\put(-2.25, -.98){\line(1, 0){.125}}
	\put(-2.25, -.98){\line(0, -1){.125}}
	\put(-2.125, -.98){\line(0, -1){.125}}
	\put(-2.25, -1.105){\line(1, 0){.125}}
	\put(-0.9, -.58){$\emptyset$} 
   \end{picture}
\end{center}
\vspace{4cm}
\caption{Branching component graph 
$\BC(B^{2, 2})$\label{fig:sup}}
\end{figure}
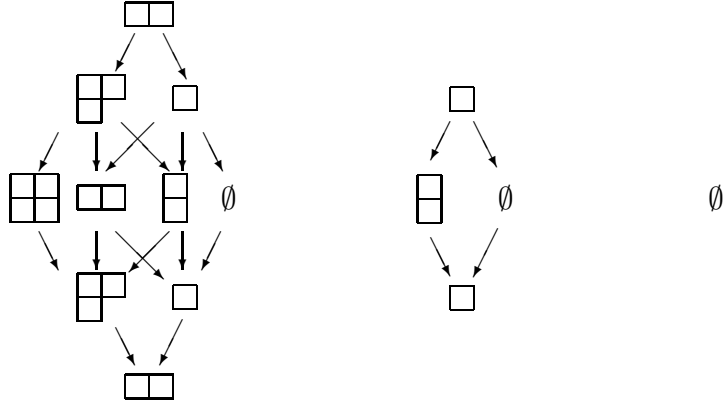

\begin{example}
    Figure \ref{fig:sup} depicts $\BC(B^{2, 2})$,  which is the
    union of $\BC(\emptyset)$,\\  $\BC((1, 1))$,  and $\BC((2, 2))$.
\end{example}

\subsection{Arbitrary partitions}

All partitions can be constructed by vertically juxtaposing some
number of rectangles all with different heights.  Just as we can build 
up these partitions,  we can also build up the branching component
graphs by combining the graphs for those constituent rectangles.

Note that while the structure of branching component graphs for
non-rectangular partitions is presented as a conjecture, the main
result of this chapter (Conjecture \ref{conj:brs}) does not depend on
its validity.  It is included here simply to present as much
information as possible about branching component graphs.  Also note
that the statements regarding the vertices and their associated
$D_{n-1}$ weights is proved; only the properties of edges are
conjectural.

\begin{conjecture}\label{prop:allparts}
    Let $\la^{(1)}, \la^{(2)}$ be partitions, let $\ell(\la^{(2)})$
    denote the number of parts of $\la^{(2)}$, and let $m_{1}^{(1)}$
    denote the multiplicity of $\la^{(1)}_{1}$ in $\la^{(1)}$.
    Suppose that $\ell(\la^{(2)})\leq m_{1}^{(1)}-2$.  Let
    $\la^{(1)}\la^{(2)}$ denote the horizontal concatenation of
    $\la^{(1)}$ and $\la^{(2)}$; i.e.
    $(\la^{(1)}\la^{(2)})_{i}=\la^{(1)}_{i}+\la^{(2)}_{i}$.  Then
    $$\BC(\la^{(1)}\la^{(2)})=\BC(\la^{(1)})\times\BC(\la^{(2)}), $$
    where $\times$ denotes the directed graph product of the two
    smaller graphs,  and the vertex ($\mu^{(1)}, \mu^{(2)}$) is labeled
    by $\mu^{(1)}\mu^{(2)}$.
\end{conjecture}

\begin{proof}
    The construction of $\pm$ diagrams treats each maximal rectangle in
    such a partition independently, and the difference between the
    number of $+$'s and $-$'s is additive, so the partitions, their
    multiplicities, and their strata are given precisely by the above 
    proposition.  
    
    Now, as in the proof of \ref{prop:recs}, we give sequences of
    operators that have the effect of adding a box to a row of a
    partition indexing a branching component vertex.  Suppose that
    $\mu = \mu^{(1)}\mu^{(2)}$ is a partition such that in
    $\BC(\la^{(1)})$ and $\BC(\la^{(2)})$, at least one of $\mu^{(1)}$
    or $\mu^{(2)}$ has an edge to the partition differing from it by
    the addition of a box in the $k^{\mathrm{th}}$ row.  Denote the
    partition differing from $\mu$ by the addition of a box in the
    $k^{\mathrm{th}}$ row by $\mu_{+k}$.  Our conjecture breaks into
    two cases.  If $k$ is not congruent modulo 2 to the number of
    parts of $\la$, we claim that the sequence
    $$
    f_{1}f_{2}\cdots f_{k}
    $$
    applied to the $D_{n-1}$ highest weight tableau of $B(\mu)$
    produces the $D_{n-1}$ highest weight tableau of $B(\mu_{+k})$.
    Otherwise, $k$ is congruent modulo 2 to the number of parts of
    $\la$.  Let $r$ be the number of columns of height $k$ in $\mu$.
    Then the sequence
    $$
    e_{k}^{r-1}e_{k-1}^{r-1}\cdots
    e_{2}^{r-1}f_{1}f_{2}^{r}f_{3}^{r}\cdots f_{k}^{r}
    $$
    applied to the $D_{n-1}$ highest weight tableau of $B(\mu)$
    produces the $D_{n-1}$ highest weight tableau of $B(\mu_{+k})$.
    Compare these statements to those of Proposition \ref{prop:recs}.
    
    The heuristic for these to be separate cases is that the first
    corresponds to removing the $+$ from a $\pm$ pair in a single
    column, whereas the second case corresponds to removing a $+$ from
    a column where it was alone.

\end{proof}

Proposition
\ref{prop:recs} and Conjecture \ref{prop:allparts} characterize all
branching component graphs in our current scope of interest, since for
a classical component $B(\mu)$ to be in $\brs$, the complement of
$\mu$ in an $r\times s$ rectangle must be tileable by vertical
dominos.

        \section{The Crystal Automorphism}
        \label{sec:KR:sigma}
        Recall that the automorphism $\sig:\brs\to\brs$, where $B^{r,s}$ is
taken as a $D_{n-1}$-crystal, corresponds to the $D_{n}^{(1)}$ Dynkin
diagram automorphism interchanging nodes $0$ and $1$.  As a result, we
know that $\sigma$ must be a crystal isomorphism when restricted to
any $D_{n-1}$-crystal corresponding to a vertex in $\BC(\brs)$.  We
may therefore describe $\sig$ by first specifying a shape-preserving
automorphism $\chsig$ of $\BC(\brs)$, then defining $\sig$ to be the
union of the corresponding crystal isomorphisms.

To explicitly determine $\sig(b)$ from $\chsig$ for any $b\in\brs$,
let $v$ be the branching component vertex such that $b\in B(v)$, and
let $T_{v}$ denote the $D_{n-1}$-highest weight vertex of $B(v)$.  We
know that for some finite sequence $i_{1},\ldots,i_{k}$ of integers in
$\{2,\ldots,n\}$, we have $\ftil{i_{1}}\cdots\ftil{i_{k}}T_{v}=b$.
Let $T_{\chsig(v)}$ be the highest weight vector of $B(\chsig(v))$.
For any such sequence of integers, we may define
$\sigma(b)=\ftil{i_{1}}\cdots\ftil{i_{k}}T_{\chsig(v)}$.

\begin{proposition}\label{prop:flip}
    Let $b\in B^{r,s}$, and suppose $\ftil{0}(b)\neq 0$.  Then the
    stratum of $\ftil{0}(b)$ must be one greater than the stratum of
    $b$,
    and therefore the
    map $\chsig$ must send a vertex from stratum $j$ to a vertex in
    stratum $-j$.
\end{proposition}

\begin{proof}
    Recall from the weight structure of type $D$ algebras that
    $\alpha_{0}=2\Lz-\Lt$ and
    $\mathrm{wt}(\ftil{0}(b))=\mathrm{wt}(b)-\alpha_{0}$.  Define
    $\mathrm{cw}(b)=\mathrm{wt}(b)-(\varphi_{0}(b)-\vare_{0}(b))\Lambda_0$.
    The above implies that
    $\mathrm{cw}(\ftil{0}(b))=\mathrm{cw}(b)+\Lt=\mathrm{cw}(b)+(e_1+e_2)$,
    where $e_{1}, e_{2}$ are the first and second elementary basis
    vectors in the weight space.  Similarly,
    $\mathrm{cw}(\ftil{i}(b))=\mathrm{cw}(b)-\alpha_i$ for
    $i=1,\ldots,n$, so that among $\ftil{1},\ldots,\ftil{n}$, only
    $\ftil{1}$ changes the $e_{1}$ component of a weight by $-1$.
    Since $\ftil{1}$ decreases the stratum by one and $\ftil{0}$
    changes the $e_{1}$ component of a weight by $+1$, it follows that
    $\ftil{0}$ increases the stratum by one.
\end{proof}

Recall that $\pm$ diagrams can be used to label vertices of
$\BC(B^{r,s})$.  We state a conjectural construction of $\chsig$ due
to Schilling and Shimozono \cite{SSh:2006} defined in terms of $\pm$
diagrams that satisfies the condition of Proposition \ref{prop:flip}.

\begin{conjecture}
    \label{conj:pmsigma}
    Let $D$ be a $\pm$ diagram of shape $\La/\la$ such that the
    complement of $\La$ in an $r\times s$ rectangle can be tiled by
    vertical dominos.  Let $m_{i}^{\circ}$, $m_{i}^{+}$, $m_{i}^{-}$
    and $m_{i}^{\pm}$ denote the number of columns of height $i$ in
    $D$ with no symbol, a $+$, a $-$, and a $\pm$ pair, respectively.
    If $r$ is even, set $m_{0}^{\circ} = s - \La_{1}$.  Note that if
    $i\equiv r+1\; (\mathrm{mod}\; 2)$ then $m_{i}^{*} = 0$ for
    $*\in\{\circ,+,-,\pm\}$.
    Then $\chsig(D)$ is the $\pm$ diagram with $m_{i+2}^{\pm}$ empty
    columns of height $i$, $m_{i}^{-}$ columns of height $i$ with a
    $+$, $m_{i}^{+}$ columns of height $i$ with a $-$, and
    $m_{i}^{\circ}$ columns of height $i+2$ with a $\pm$ pair for all
    $0\leq i\leq r$.

\end{conjecture}

This conjectural $\chsig$ also respects a symmetry that appears in the
poset of $D_{n}$ components of $B^{r,s}$.  To state this precisely, we
need to know which $D_{n}$ components of $B^{r,s}$ have $\mu$ as a
$D_{n-1}$ highest weight.  In other words, we want to know for which
outer shapes $\La$ there is a $\pm$ diagram with inner shape $\mu$.
The smallest such partition, which we denote by $\muz$, is the outer
shape of the $\pm$ diagram with inner shape $\mu$ none of whose
columns have a $\pm$ pair.  Complementarily, the largest such
partition, which we denote by $\muo$, is the outer shape of the $\pm$
diagram of which $\mu$ is the inner shape and such that every column
of height less than $r$ has at least one $+$ or one $-$ in it.  It is
also the result of adding a vertical domino with a $\pm$ pair to those
columns of $\muz$ that have fewer than $r$ boxes and do not have a $+$
or $-$ in them.

\begin{remark}
    In the above, $\muz$ and $\muo$ depend on the parity of $r$.
\end{remark}

\begin{example}\label{ex:run1}
    Consider those classical components $B(\la)$ of $B^{4,4}$ for
    which $\mu=(2,1)$ is associated with a vertex of $\BC(\la)$.  The
    smallest such $\la$ is $\muz=(2,2)$; the largest is
    $\muo=(4,4,1,1)$.  The intermediate partitions are $(2,2,1,1)$,
    $(3,3)$, $(3,3,1,1)$, and $(4,4)$.
\end{example}

Given a partition $\la$ whose complement is tiled by vertical dominos,
let $\bar{\la}$ be the partition whose columns have length equal to
the floor of half the length of the columns of $\la$.  In words,
$\bar{\la}$ is the result of including only every other part of $\la$,
omitting the first part if $\la$ has an odd number of parts.  It is
easy to see that $\overline{\muz}$ and $\overline{\muo}$ differ by a
horizontal strip, and that the partitions containing $\mu$ as
described above form the poset interval between them in Young's
lattice.  Consider the multiset of rows with boxes in
$\overline{\muo}/\overline{\muz}$.  Taking the complement in this
multiset provides a canonical symmetry of this interval.  Explicitly,
define the map $t$ from this interval to itself by saying that both
$\la/\overline{\muz}$ and $\overline{\muo}/t(\la)$ must have $m_{i}$
boxes in row $i$.

Conjecture \ref{conj:pmsigma} implies that $\chsig$ respects this
symmetry.  Observe that the outer shape of $\chsig(D)$ differs from
outer($D$) by replacing $m_{i}^{\pm}$ columns with $\pm$ pairs by
$c_{i-2}-m_{i}^{\pm}$ many such columns when $i\equiv r
(\mathrm{mod}2)$.  This is precisely the same involution of the
interval as $t$ in the above paragraph.

\begin{figure}
\begin{center}\setlength{\unitlength}{1in}
    \begin{picture}(0,.3)
	\put(-.25,0){\line(1,0){.5}}
	\put(-.25,-.125){\line(1,0){.5}}	
	\put(-.25,0){\line(0,-1){.25}}
	\put(-.125,0){\line(0,-1){.25}}
	\put(0,0){\line(0,-1){.125}}
	\put(.125,0){\line(0,-1){.125}}
	\put(.25,0){\line(0,-1){.125}}	
	\put(-.25,-.25){\line(1,0){.125}}
	\put(-.13,-.3){\line(-1,-1){.2}}
	\put(.13,-.3){\line(1,-1){.2}}	
	\put(-.55,-.55){\line(1,0){.375}}
	\put(-.55,-.675){\line(1,0){.375}}	
	\put(-.55,-.55){\line(0,-1){.25}}
	\put(-.425,-.55){\line(0,-1){.25}}
	\put(-.3,-.55){\line(0,-1){.125}}
	\put(-.175,-.55){\line(0,-1){.125}}
	\put(-.55,-.8){\line(1,0){.125}}
	\put(.065,-.61){\line(1,0){.5}}
	\put(.065,-.735){\line(1,0){.5}}	
	\put(.065,-.61){\line(0,-1){.125}}
	\put(.19,-.61){\line(0,-1){.125}}
	\put(.315,-.61){\line(0,-1){.125}}
	\put(.44,-.61){\line(0,-1){.125}}
	\put(.565,-.61){\line(0,-1){.125}}
	\put(-.3,-.8){\line(0,-1){.2}}
	\put(.3,-.8){\line(0,-1){.2}}	
	\put(-.2,-.8){\line(2,-1){.41}}		
	\put(-.425,-1.1){\line(1,0){.25}}
	\put(-.42,-1.225){\line(1,0){.25}}	
	\put(-.425,-1.1){\line(0,-1){.25}}
	\put(-.3,-1.1){\line(0,-1){.25}}
	\put(-.175,-1.1){\line(0,-1){.125}}
	\put(-.425,-1.35){\line(1,0){.125}}
	\put(.125,-1.16){\line(1,0){.375}}
	\put(.125,-1.285){\line(1,0){.375}}	
	\put(.125,-1.16){\line(0,-1){.125}}
	\put(.25,-1.16){\line(0,-1){.125}}
	\put(.375,-1.16){\line(0,-1){.125}}
	\put(.5,-1.16){\line(0,-1){.125}}
	\put(-.3,-1.4){\line(1,-1){.2}}
	\put(.3,-1.4){\line(-1,-1){.2}}	
	\put(-.125,-1.7){\line(1,0){.25}}
	\put(-.125,-1.825){\line(1,0){.25}}	
	\put(-.125,-1.7){\line(0,-1){.125}}
	\put(0,-1.7){\line(0,-1){.125}}
	\put(.125,-1.7){\line(0,-1){.125}}
    \end{picture}
\end{center}

\vspace{4.5cm}
\caption{Interval in Young's lattice between $(2)$ and $(4,1)$\label{fig:run}}
\end{figure}
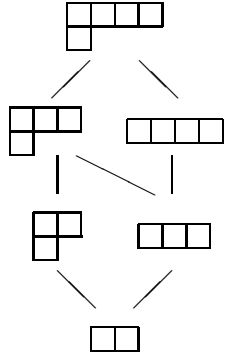

\begin{example}
    Consider $\mu=(2,1)$ in $B^{4,4}$ as in Example
    \ref{ex:run1} above.  The interval between $\overline{\muz}$ and
    $\overline{\muo}$ is shown in figure \ref{fig:run}.
\end{example}

In the case of $r=2$ these intervals are simply linear orderings, and
Conjecture \ref{conj:pmsigma} has been proved in this special case
(Theorem \ref{thm:b2s}; see Section \ref{sec:B2s:r2}).  This symmetry
appears to be the only generally consistent procedure for determining
what the classical component of $\chsig(v)$ should be, and it is the
most natural generalization of the behavior of $\chsig$ in the case of
$B^{2,s}$.

We have now defined $\chsig$ on $\BC(B^{r,s})$, which in turn defines 
$\sigma$ on $B^{r,s}$.

\begin{conjecture}[Restatement of Conjecture \ref{conj:brs}]
    Define a $D_{n}^{(1)}$-crystal $\tbrs$ by adding $0$ arrows to the
    classical crystals according to the rule
    $$\ftil{0}=\sig\ftil{1}\sig.$$ Then if $\brs$ exists,
    $\tbrs\simeq\brs$.
\end{conjecture}

    %
    %

    \newchap{The Family $B^{2,s}$ of Type $D$ Kirillov-Reshetikhin 
    Modules}
    \label{sec:B2s}
    
	\section{Specialization of previous work to $r=2$}
        \label{sec:B2s:tabs}

In this chapter we restrict our attention to $B^{2,s}$, proving
Theorems \ref{thm:b2s} and \ref{thm:b2stabs}.  We begin by
specializing the statements and constructions of Chapter \ref{sec:KR}
to the case of $r=2$.

\subsection{Classical tableaux}

First, we note that given the assumption that our tableaux are height 
$2$ rectangles, the following simpler Criterion replaces 
Criterion \ref{legal}.

\begin{criterion} \label{legal2}
{}~
\begin{enumerate}
    \item \label{legal2:1} If $ab$ is in the filling, then $a\leq b$;
    \item \label{legal2:2} If $\stack{a}{b}$ is in the filling, then $b\nleq a$;
    \item \label{legal2:3} No configuration of the form $\stack{a}{\phantom{\bar{a}}}
      \stack{a}{\bar{a}}$ or $\stack{a}{\bar{a}}\stack{\phantom{a}}{\bar{a}}$ appears;
    \item \label{legal2:4} No configuration of the form 
      $\stack{n-1}{n}\ldots\stack{n}{\overline{n-1}}$ or 
      $\stack{n-1}{\bar{n}}\ldots\stack{\bar{n}}{\overline{n-1}}$ appears;
    \item \label{legal2:5} No configuration of the form $\stack{1}{\bar{1}}$ appears.
\end{enumerate}
\end{criterion}
Note that for $k\geq 2$, condition $5$ follows from conditions $1$ and
$3$.  Also, observe that Criterion \ref{legal2} is unchanged by
replacing condition 4 with the following:
\begin{enumerate}
    \item[(4a)] {\it No configuration of the form }
    $\stack{n-1}{n}\stack{n}{\overline{n-1}}$ {\it or}
    $\stack{n-1}{\bar{n}}\stack{\bar{n}}{\overline{n-1}}$ {\it
    appears.}
\end{enumerate}
    
To see this equivalence, observe that by conditions $1$ and $2$ the
only columns that can appear between $\stack{n-1}{n}$ and
$\stack{n}{\overline{n-1}}$ are $\stack{n-1}{n}$,
$\stack{n-1}{\overline{n-1}}$, and $\stack{n}{\overline{n-1}}$, and if
present, they must appear in that order from left to right.  If a
column of the form $\stack{n-1}{\overline{n-1}}$ appears, we have a
configuration of the form $\stack{n-1}{\phantom{\overline{n-1}}}
\stack{n-1}{\overline{n-1}}$, which is forbidden by condition $3$.  On
the other hand, if no column of the form $\stack{n-1}{\overline{n-1}}$
appears, the columns $\stack{n-1}{n}$ and $\stack{n}{\overline{n-1}}$
are adjacent, which is disallowed by condition 4a.

\subsection{Sliding algorithm}\label{sec:slide}

In section \ref{sec:KR:tabs}, we saw the Lecouvey relations of the
Type $D$ plactic monoid.  Here we have the derived sliding algorithm
for the case of height $2$ rectangles.
\begin{enumerate}

    \item  If $x\neq \bar{z}$, then
    \begin{eqnarray*}
	\left.
	\begin{array}{c|c|c|}
	    \cline{3-3}
	     \multicolumn{1}{c}{} & & y \\
	    \cline{2-3}
	    & x & z  \\
	    \cline{2-3}
	\end{array}
	\equiv
	\begin{array}{c|c|c|}
	    \cline{2-3}
	      & x& y \\
	    \cline{2-3}
	    \multicolumn{1}{c}{}&  & z  \\
	    \cline{3-3}
	\end{array}	
	\equiv
	\begin{array}{|c|c|c}
	    \cline{1-2}
	     x &  y \\
	    \cline{1-2}
	    z & \multicolumn{2}{c}{} \\
	    \cline{1-1}
	\end{array}
	\right. \:
	\textrm{ for }x\leq y<z,
	\\
	\textrm{ and }
	\left.
	\begin{array}{c|c|c|}
	    \cline{3-3}
	     \multicolumn{1}{c}{} & & x \\
	    \cline{2-3}
	    & y & z  \\
	    \cline{2-3}
	\end{array}
	\equiv
	\begin{array}{|c|c|}
	    \cline{1-1}
	       x& \multicolumn{1}{c}{} \\
	    \cline{1-2}
	     y & z  \\
	    \cline{1-2}
	\end{array}	
	\equiv
	\begin{array}{|c|c|c}
	    \cline{1-2}
	     x &  z \\
	    \cline{1-2}
	    y & \multicolumn{2}{c}{} \\
	    \cline{1-1}
	\end{array}
	\right. \:	
	\textrm{ for }x<y\leq z	.
    \end{eqnarray*}
    
    \item  If $1<x<n$ and $x\leq y\leq\bar{x}$, then 
    \begin{eqnarray*}
	\left.
	\begin{array}{c|c|c|}
	    \cline{3-3}
	     \multicolumn{1}{c}{} & & y \\
	    \cline{2-3}
	    & x-1 & \overline{x-1}\rule{0pt}{2.5ex}  \\
	    \cline{2-3}
	\end{array}
	\equiv
	\begin{array}{c|c|c|}
	    \cline{2-3}
	      & x-1& y \\
	    \cline{2-3}
	    \multicolumn{1}{c}{}&  & \overline{x-1}\rule{0pt}{2.5ex}  \\
	    \cline{3-3}
	\end{array}	
	\equiv	
	\begin{array}{|c|c|c}
	    \cline{1-2}
	     x &  y \\
	    \cline{1-2}
	    \bar{x} & \multicolumn{2}{c}{} \\
	    \cline{1-1}
	\end{array}
	\right. \:
	\\
	\textrm{ and }
	\left.
	\begin{array}{c|c|c|}
	    \cline{3-3}
	     \multicolumn{1}{c}{} & & x \\
	    \cline{2-3}
	    & y & \overline{x}  \\
	    \cline{2-3}
	\end{array}
	\equiv
	\begin{array}{|c|c|}
	    \cline{1-1}
	       x-1& \multicolumn{1}{c}{} \\
	    \cline{1-2}
	     y & \overline{x-1}\rule{0pt}{2.5ex}  \\
	    \cline{1-2}
	\end{array}	
	\equiv	
	\begin{array}{|c|c|c}
	    \cline{1-2}
	     x-1 &  \overline{x-1}\rule{0pt}{2.5ex} \\
	    \cline{1-2}
	    y & \multicolumn{2}{c}{} \\
	    \cline{1-1}
	\end{array}
	\right. \qquad\quad	
	.
    \end{eqnarray*}
    
    \item  If $x\leq n-1$, then
    \begin{eqnarray*}
	&&\phantom{\textrm{ and }}
	\left\{
	\begin{array}{c}
	    \left.
	    \begin{array}{c|c|c|}
		\cline{3-3}
		 \multicolumn{1}{c}{} & & \bar{n} \\
		\cline{2-3}
		& n & \overline{x}  \\
		\cline{2-3}
	    \end{array}
	    \equiv
	    \begin{array}{|c|c|}
		\cline{1-1}
		   \bar{n} & \multicolumn{1}{c}{} \\
		\cline{1-2}
		 n & \bar{x}  \\
		\cline{1-2}
	    \end{array}	
	    \equiv	    
	    \begin{array}{|c|c|c}
		\cline{1-2}
		 \bar{n} &  \bar{x} \\
		\cline{1-2}
		n & \multicolumn{2}{c}{} \\
		\cline{1-1}
	    \end{array}
	    \right. \:
	    \\ 
	    \\
	    \left.
	    \begin{array}{c|c|c|}
		\cline{3-3}
		 \multicolumn{1}{c}{} & & n \\
		\cline{2-3}
		& \bar{n} & \overline{x}  \\
		\cline{2-3}
	    \end{array}
	    \equiv
	    \begin{array}{|c|c|}
		\cline{1-1}
		   n& \multicolumn{1}{c}{} \\
		\cline{1-2}
		 \bar{n} & \bar{x}  \\
		\cline{1-2}
	    \end{array}	
	    \equiv	    
	    \begin{array}{|c|c|c}
		\cline{1-2}
		 n &  \bar{x} \\
		\cline{1-2}
		\bar{n} & \multicolumn{2}{c}{} \\
		\cline{1-1}
	    \end{array}
	    \right. \:
	    \\
	\end{array}
	\right.
	\\&&
	\textrm{ and }
	\left\{
	\begin{array}{c}
	    \left.
	    \begin{array}{c|c|c|}
		\cline{3-3}
		 \multicolumn{1}{c}{} & & \bar{n} \\
		\cline{2-3}
		& x & n  \\
		\cline{2-3}
	    \end{array}
	    \equiv
	    \begin{array}{c|c|c|}
		\cline{2-3}
		  & x& \bar{n} \\
		\cline{2-3}
		\multicolumn{1}{c}{}&  & n  \\
		\cline{3-3}
	    \end{array}	
	    \equiv	    
	    \begin{array}{|c|c|c}
		\cline{1-2}
		 x &  \bar{n} \\
		\cline{1-2}
		n & \multicolumn{2}{c}{} \\
		\cline{1-1}
	    \end{array}
	    \right. \:
	    \\ 
	    \\
	    \left.
	    \begin{array}{c|c|c|}
		\cline{3-3}
		 \multicolumn{1}{c}{} & & n \\
		\cline{2-3}
		& x & \overline{n}  \\
		\cline{2-3}
	    \end{array}
	    \equiv
	    \begin{array}{c|c|c|}
		\cline{2-3}
		  & x& n \\
		\cline{2-3}
		\multicolumn{1}{c}{}&  & \bar{n}  \\
		\cline{3-3}
	    \end{array}	
	    \equiv	    
	    \begin{array}{|c|c|c}
		\cline{1-2}
		 x &  n \\
		\cline{1-2}
		\bar{n} & \multicolumn{2}{c}{} \\
		\cline{1-1}
	    \end{array}
	    \right. \:
	    \\
	\end{array}
	\right. .
    \end{eqnarray*}
    
    \item 
    \begin{eqnarray*}
	&&\phantom{\textrm{ and }}
	\left\{
	\begin{array}{c}
	    \left.
	    \begin{array}{c|c|c|}
		\cline{3-3}
		 \multicolumn{1}{c}{} & & n \\
		\cline{2-3}
		& \bar{n} & \bar{n}  \\
		\cline{2-3}
	    \end{array}
	    \equiv
	    \begin{array}{|c|c|}
		\cline{1-1}
		  n-1 & \multicolumn{1}{c}{} \\
		\cline{1-2}
		 \bar{n} & \overline{n-1}\rule{0pt}{2.5ex}  \\
		\cline{1-2}
	    \end{array}	
	    \equiv	    
	    \begin{array}{|c|c|c}
		\cline{1-2}
		 n-1 &  \overline{n-1}\rule{0pt}{2.5ex} \\
		\cline{1-2}
		\bar{n} & \multicolumn{2}{c}{} \\
		\cline{1-1}
	    \end{array}
	    \right. \:
	    \\ 
	    \\
	    \left.
	    \begin{array}{c|c|c|}
		\cline{3-3}
		 \multicolumn{1}{c}{} & & \bar{n} \\
		\cline{2-3}
		& n & n  \\
		\cline{2-3}
	    \end{array}
	    \equiv
	    \begin{array}{|c|c|}
		\cline{1-1}
		   n-1& \multicolumn{1}{c}{} \\
		\cline{1-2}
		 n & \overline{n-1}\rule{0pt}{2.5ex}  \\
		\cline{1-2}
	    \end{array}	
	    \equiv	    
	    \begin{array}{|c|c|c}
		\cline{1-2}
		 n-1 &  \overline{n-1}\rule{0pt}{2.5ex} \\
		\cline{1-2}
		n & \multicolumn{2}{c}{} \\
		\cline{1-1}
	    \end{array}
	    \right. \:
	    \\	    
	\end{array}
	\right.
	\\
	&&
	\textrm{ and }
	\left\{
	\begin{array}{c}
	    \left.
	    \begin{array}{c|c|c|}
		\cline{3-3}
		 \multicolumn{1}{c}{} & & \bar{n} \\
		\cline{2-3}
		& n-1 & \overline{n-1} \rule{0pt}{2.5ex} \\
		\cline{2-3}
	    \end{array}
	    \equiv
	    \begin{array}{c|c|c|}
		\cline{2-3}
		  & n-1& \bar{n} \\
		\cline{2-3}
		\multicolumn{1}{c}{}&  & \overline{n-1}\rule{0pt}{2.5ex}  \\
		\cline{3-3}
	    \end{array}	
	    \equiv	    
	    \begin{array}{|c|c|c}
		\cline{1-2}
		 \bar{n} &  \bar{n} \\
		\cline{1-2}
		n & \multicolumn{2}{c}{} \\
		\cline{1-1}
	    \end{array}
	    \right. \:
	    \\ 
	    \\
	    \left.
	    \begin{array}{c|c|c|}
		\cline{3-3}
		 \multicolumn{1}{c}{} & & n \\
		\cline{2-3}
		& n-1 & \overline{n-1}\rule{0pt}{2.5ex}  \\
		\cline{2-3}
	    \end{array}
	    \equiv
	    \begin{array}{c|c|c|}
		\cline{2-3}
		  & n-1& n \\
		\cline{2-3}
		\multicolumn{1}{c}{}&  & \overline{n-1}\rule{0pt}{2.5ex}  \\
		\cline{3-3}
	    \end{array}	
	    \equiv	    
	    \begin{array}{|c|c|c}
		\cline{1-2}
		 n &  n \\
		\cline{1-2}
		\bar{n} & \multicolumn{2}{c}{} \\
		\cline{1-1}
	    \end{array}
	    \right. \:
	    \\		        
	\end{array}
	\right. .
    \end{eqnarray*}
\end{enumerate}
If a word is composed entirely of barred letters or entirely of
unbarred letters, only relation ($1$) (the Knuth relation) applies,
and the type $A$ \textit{jeu de taquin} may be used.

\subsection{Technical definitions for operations on tableaux}

\begin{definition} \label{def:null}
A null-configuration of size $k$ is
\begin{equation*}
    \framebox{$\emptyset_{k}$} = 
    \begin{cases}
	\underbrace{
	\stack{1}{\bar{2}}
	\cdots
	\stack{1}{\bar{2}}
	}_{\frac{k}{2}}
	\;
	\underbrace{
	\stack{2}{\bar{1}}
	\cdots
	\stack{2}{\bar{1}}
	}_{\frac{k}{2}} & \textrm{ if $k$ is even }\\
	\underbrace{
	\stack{1}{\bar{2}}
	\cdots
	\stack{1}{\bar{2}}
	\,
	\stack{2}{\bar{2}}
	}_{\frac{k+1}{2}}
	\underbrace{
	\stack{2}{\bar{1}}
	\cdots
	\stack{2}{\bar{1}}
	}_{\frac{k-1}{2}} & \textrm{ if $k$ is odd. }
    \end{cases}
\end{equation*}

\end{definition}

Null-configurations are named thus because $\etil{i}$ and $\ftil{i}$
for $i=2, \ldots, n$ send $b$ to $0$, where $b$ is the $2\times s$
tableau which is a null-configuration of size $s$.  Therefore, $b$ is
the basis vector for the trivial representation of $D_{n-1}$ in the
$D_{n}$ crystal $B(s\Lt)$.  Put another way, inserting a
null-configuration into a tableau $T$ has no effect on $\vare_{i}(T)$
or $\varphi_{i}(T)$ for $i=2, \ldots, n$.

\begin{definition}
    \label{def:redtab}
    Let $T\in B(s\Lt)$.  The reduced form of $T$, denoted by $\Tred$,
    is the skew tableau that results from removing all $1$'s,
    $\bar{1}$'s, and any null-configuration from $T$.  The
    rectification of $\Tred$ by the sliding algorithm in section
    \ref{sec:slide} is called the completely reduced form of $T$.
\end{definition}

Since the type $D$ plactic operations cannot increase the length of a 
column of a skew tableau,  we know that the completely reduced form of 
$T$ must have no more than two rows.

\subsection{Properties of $B^{r,s}$}

The properties of $B^{r,s}$ in Conjecture \ref{conj:Brs} can be
described more simply in the case of $r=2$.  We may restate the
classical direct sum decomposition (item 1) as
\begin{equation}
    \label{eq:decomp}
    B^{2,s} \simeq \bigoplus_{i=0}^{s}i\Lt
\end{equation}
and restate the description of the energy function as
$D_{B^{2,s}}(b)=k$ if $b$ is in a $D_{n}$ component with highest
weight $k\Lt$.

\subsection{Specialization of $\hsig$.}\label{sec:hsig}
In the special case of $r=2$ the definition of $\hsig$ on the set of
branching component vertices has a very easy combinatorial
description, since there is only one pair of rows in which $+$'s and
$-$'s can appear in its $\pm$ diagrams.

Suppose $v\in\BC(k\Lt)$ is labeled by the partition $(\la_{1},
\la_{2})$.  The corresponding $\pm$ diagram must have $k-\la_{1}$
columns with a $\pm$ pair, so $\chsig(v)\in\BC(s-(k-\la_{1}))$.  By
Proposition \ref{prop:recs}, we know that the strata of $\BC(k\Lt)$
are multiplicity free for all $k\geq 0$, so Proposition
\ref{prop:flip} specifies $\chsig(v)$ by telling us that it is in the
opposite stratum of $v$.

\begin{example}
    The action of $\hsig$ on $\Sup(\tbt)$ is given in Figure
    \ref{fig:sigma}.

\end{example}

\begin{figure}
\begin{center}\setlength{\unitlength}{1in}
    \begin{picture}(-2.2,0)
	\put(-2.5,0){\line(1,0){.25}}
	\put(-2.5,0){\line(0,-1){.125}}
	\put(-2.25,0){\line(0,-1){.125}}
	\put(-2.5,-.125){\line(1,0){.25}}
	\put(-2.375,0){\line(0,-1){.125}}
%

	\put(-2.75,-.38){\line(1,0){.25}}
	\put(-2.75,-.38){\line(0,-1){.25}}
	\put(-2.625,-.38){\line(0,-1){.25}}
	\put(-2.75,-.505){\line(1,0){.25}}
	\put(-2.5,-.38){\line(0,-1){.125}}
	\put(-2.75,-.63){\line(1,0){.125}}
%

	\put(-2.25,-.44){\line(1,0){.125}}
	\put(-2.25,-.44){\line(0,-1){.125}}
	\put(-2.125,-.44){\line(0,-1){.125}}
	\put(-2.25,-.565){\line(1,0){.125}}
%

	\put(-3.1,-.9){\line(1,0){.25}}
	\put(-3.1,-.9){\line(0,-1){.25}}
	\put(-2.975,-.9){\line(0,-1){.25}}
	\put(-3.1,-1.025){\line(1,0){.25}}
	\put(-2.85,-.9){\line(0,-1){.25}}
	\put(-3.1,-1.15){\line(1,0){.25}}
%

	\put(-2.75,-.96){\line(1,0){.25}}
	\put(-2.75,-.96){\line(0,-1){.125}}
	\put(-2.5,-.96){\line(0,-1){.125}}
	\put(-2.75,-1.085){\line(1,0){.25}}
	\put(-2.625,-.96){\line(0,-1){.125}}
%

	\put(-2.3,-.9){\line(1,0){.125}}
	\put(-2.3,-.9){\line(0,-1){.25}}
	\put(-2.175,-.9){\line(0,-1){.25}}
	\put(-2.3,-1.025){\line(1,0){.125}}
	\put(-2.3,-1.15){\line(1,0){.125}}
%

	\put(-2,-1.08){$\emptyset$}

	\put(-2.75,-1.42){\line(1,0){.25}}
	\put(-2.75,-1.42){\line(0,-1){.25}}
	\put(-2.625,-1.42){\line(0,-1){.25}}
	\put(-2.75,-1.545){\line(1,0){.25}}
	\put(-2.5,-1.42){\line(0,-1){.125}}
	\put(-2.75,-1.67){\line(1,0){.125}}
%

	\put(-2.25,-1.48){\line(1,0){.125}}
	\put(-2.25,-1.48){\line(0,-1){.125}}
	\put(-2.125,-1.48){\line(0,-1){.125}}
	\put(-2.25,-1.608){\line(1,0){.125}}
%
	
	\put(-2.5,-1.95){\line(1,0){.25}}
	\put(-2.5,-1.95){\line(0,-1){.125}}
	\put(-2.25,-1.95){\line(0,-1){.125}}
	\put(-2.5,-2.075){\line(1,0){.25}}
	\put(-2.375,-1.95){\line(0,-1){.125}}
    \end{picture}
\end{center}
\begin{center}\setlength{\unitlength}{1in}
    \begin{picture}(-5.1,.3)
	\put(-2.25,.06){\line(1,0){.125}}
	\put(-2.25,.06){\line(0,-1){.125}}
	\put(-2.125,.06){\line(0,-1){.125}}
	\put(-2.25,-.065){\line(1,0){.125}}
%

	\put(-2.42,-.4){\line(1,0){.125}}
	\put(-2.42,-.4){\line(0,-1){.25}}
	\put(-2.295,-.4){\line(0,-1){.25}}
	\put(-2.42,-.525){\line(1,0){.125}}
	\put(-2.42,-.65){\line(1,0){.125}}
%

	\put(-2,-.58){$\emptyset$}
%

	\put(-2.25,-.98){\line(1,0){.125}}
	\put(-2.25,-.98){\line(0,-1){.125}}
	\put(-2.125,-.98){\line(0,-1){.125}}
	\put(-2.25,-1.105){\line(1,0){.125}}
	\put(-0.9,-.58){$\emptyset$} 
%
%
	\qbezier[2000](-0.95,-.45)(-2.125,-0.2)(-3.36,-.45)
	\put(-0.94,-.45){\vector(4,-1){0}}
	\put(-3.37,-.45){\vector(-4,-1){0}}
	\put(-2.3,-1){\vector(-4,3){1.2}}
	\put(-2.3,-1){\vector(4,-3){0}}
	\put(-2.3,-.1){\vector(-4,-3){1.2}}
	\put(-2.3,-.1){\vector(4,3){0}}
	\qbezier[200](-2.46,-.65)(-3,-0.85)(-3.58,-.65)
	\put(-2.44,-.64){\vector(4,1){0}}
	\put(-3.6,-.64){\vector(-4,1){0}}
	\put(-3.85,.3){\vector(0,-1){1.72}}
	\put(-3.85,.32){\vector(0,1){0}}
	\qbezier[150](-4.15,-.18)(-4.324,-.53)(-4.15,-.88)
	\put(-4.14,-.16){\vector(1,2){0}}
	\put(-4.14,-.9){\vector(1,-2){0}}

   \end{picture}
\end{center}
\vspace{4cm} 
\caption{Definition of $\hsig$ on $\Sup(\tbt)$\label{fig:sigma}}
\end{figure}

	\section[Realization of $B^{2,s}$ as a Set of Young-like
	Tableaux]{Realization of $B^{2,s}$ as a Set of\\ Young-like
	Tableaux}
        \label{sec:B2s:tabs2}
        \subsection{Bijection between affine and classical
tableaux}\label{sec:bij}

We now define a $D_n$-crystal with vertices labeled by the set
$\mathcal{T}(s)$ of tableaux of shape $(s,s)$ satisfying conditions
\ref{legal2:1}, \ref{legal2:2}, and \ref{legal2:4} of Criterion
\ref{legal2}, but not condition \ref{legal2:3}.  We will construct a
bijection between $\mathcal{T}(s)$ and the vertices of
$\bigoplus_{i=0}^{s}B(i\Lt)$, so that $\mathcal{T}(s)$ may be viewed
as a $D_n$-crystal with the classical decomposition \eqref{eq:decomp}.
In section \ref{sec:B2s:tabs2} we will consider $\ftil{0}$ and
$\etil{0}$ on $\mathcal{T}(s)$ as given by $\chsig$ to give it the
structure of a perfect $D_n^{(1)}$-crystal, which we call $\tbs$; we
will then prove that $\tbs$ is the only crystal satisfying Conjecture
\ref{conj:Brs} for the case of $r=2$.

\begin{definition}
    \label{def:affine_tab}
    The set $\mathcal{T}(s)$ of affine tableaux of width $s$ is the
    set of tableaux satisfying conditions \ref{legal2:1},
    \ref{legal2:2}, and \ref{legal2:4} of Criterion \ref{legal2}.
\end{definition}

\begin{proposition}\label{ablock}
Let $T\in\mathcal{T}(s)\setminus B(s\Lt)$ with $T\neq
\stack{1}{\bar{1}}\cdots \stack{1}{\bar{1}}$, and define
$\bar{\bar{i}}=i$ for $1\leq i \leq n$.  Then there is a unique
$a\in\{1,\ldots,n,\bar{n}\}$ and $m\in\mathbb{Z}_{>0}$ such that $T$
contains one of the following configurations (called an
$a$-configuration) in consecutive columns:
\begin{eqnarray*}
\stack{a}{b_{1}}
\underbrace{\stack{a}{\bar{a}} \cdots \stack{a}{\bar{a}}}_{m}
\stack{c_{1}}{d_{1}} ,&&
\textrm{where $b_{1}\neq \bar{a}$, and  
($c_{1}\neq a$ or $d_{1}\neq\bar{a}$)};\\
\stack{b_{2}}{c_{2}}
\underbrace{\stack{a}{\bar{a}} \cdots \stack{a}{\bar{a}}}_{m}
\stack{d_{2}}{\bar{a}} ,&&
\textrm{where $d_{2}\neq a$, and  
($b_{2}\neq a$
or $c_{2}\neq\bar{a}$)};\\
\stack{b_{3}}{c_{3}}
\underbrace{\stack{a}{\bar{a}} \cdots \stack{a}{\bar{a}}}_{m+1}
\stack{d_{3}}{e_{3}},&&
\textrm{where $b_{3}\neq a$ and $e_{3}\neq\bar{a}$}.
\end{eqnarray*}
\end{proposition}

\begin{proof}
If $s=1$, the set $\mathcal{T}(s)\setminus B(s\Lt)$ contains only
$\stack{1}{\bar{1}}$, so the statement of the proposition is empty.
Assume that $s\ge 2$.  The existence of an $a$-configuration for some
$a\in\{1,\ldots,n,\bar{n}\}$ follows from the fact that $T$ violates
condition~\ref{legal2:3} of Criterion~\ref{legal2}.  The conditions on
$b_{i}, c_{i}, d_{i}$ for $i=1,2,3$ and $e_{3}$ mean that $m$ is
chosen to maximize the size of the $a$-configuration.
Condition~\ref{legal2:1} of Criterion~\ref{legal2} and the conditions
on the parameters $b_i,c_i,d_i,e_3$ imply that there can be no other
$a$-configurations in $T$.
\end{proof}

The map $D_{2,s}:\mathcal{T}(s)\to \bigoplus_{k=0}^{s}B(k\Lt)$, called
the height-two drop map, is defined as follows for $T\in
\mathcal{T}(s)$.  If $T=\stack{1}{\bar{1}}\cdots \stack{1}{\bar{1}}$,
then $D_{2,s}(T)=\es\in B(0)$.  If $T\in B(s\Lt)$, $D_{2,s}(T)=T$.
Otherwise $T$ contains a unique $a$-configuration by
Proposition~\ref{ablock}, and $D_{2,s}(T)$ is obtained from $T$ by
removing $\underbrace{\stack{a}{\bar{a}} \cdots
\stack{a}{\bar{a}}}_{m}$ from it.

\begin{theorem}
    Let $T\in \mathcal{T}(s)$.  Then $D_{2,s}(T)$ satisfies Criterion
    \ref{legal2}, and is therefore a tableau in
    $\bigoplus_{k=0}^{s}B(k\Lt)$.
\end{theorem}

\begin{proof}
    Condition $1$ is satisfied since the relation $\leq$ on our
    alphabet is transitive.  Conditions $2$ and $5$ are automatically
    satisfied, since the columns that remain are not changed.
    Condition $3$ is satisfied since by Proposition \ref{ablock},
    there can be no more than one $a$-configuration in $T$.  Condition
    $4$ is satisfied since $D_{2,s}$ does not remove any columns of
    the form $\stack{n-1}{n}$, $\stack{n-1}{\bar{n}}$,
    $\stack{n}{\overline{n-1}}$, or $\stack{\bar{n}}{\overline{n-1}}$.
\end{proof}

Proposition \ref{filling} shows that $D_{2,s}$ is a bijection by
constructing its inverse.

\begin{example}
We have
\begin{displaymath}
T=    
\left.
\begin{array}{|c|c|c|c|}
    \hline
    1 & 2 & 3& 3 \\
    \hline
    \overline{4} & \overline{2} & \overline{2} & \overline{1} \\
    \hline
\end{array}
\right. ,\quad
D_{2,4}(T)=
\left.
\begin{array}{|c|c|c|}
    \hline
    1 & 3& 3 \\
    \hline
    \overline{4} & \overline{2} & \overline{1} \\
    \hline
\end{array}
\right.  .
\end{displaymath}
\end{example}

The inverse of $D_{2,s}$ is the height-two fill map
$F_{2,s}:\bigoplus_{k=0}^{s}B(k\Lt) \to\mathcal{T}(s)$.  Let
$t=\stack{a_{1}}{b_{1}}\cdots\stack{a_{k}}{b_{k}} \in B(k\Lt)$.  If
$k=s$, $F_{2,s}(t)=t$.  If $k<s$, then $F_{2,s}(t)$ is obtained by
finding a subtableau $\stack{a_{i}}{b_{i}}\stack{a_{i+1}}{b_{i+1}}$ in
$t$ such that
\begin{criterion}\label{eq:crit}
    $$b_{i}\leq\bar{a}_{i}\leq b_{i+1}\qquad\textrm{ or }\qquad
    a_{i}\leq
    \bar{b}_{i+1}\leq a_{i+1}.$$
\end{criterion}
(Recall that $\bar{\bar{i}}=i$ for $i\in\{1,\ldots,n\}$.)  Note that
the first pair of inequalities imply that $a_{i}$ is unbarred, and the
second pair of inequalities imply that $b_{i+1}$ is barred.  We may
therefore insert between columns $i$ and $i+1$ of $t$ either the
configuration $\underbrace{\stack{a_{i}}{\bar{a}_{i}} \cdots
\stack{a_{i}}{\bar{a}_{i}}}_{s-k}$ or
$\underbrace{\stack{\bar{b}_{i+1}}{b_{i+1}} \cdots
\stack{\bar{b}_{i+1}}{b_{i+1}}}_{s-k}$, depending on which part of
Criterion \ref{eq:crit} is satisfied.  We say that $i$ is the filling
location of $t$.  If no such subtableau exists, then $F_{2,s}$ will
either prepend $\underbrace{\stack{\bar{b}_{1}}{b_{1}} \cdots
\stack{\bar{b}_{1}}{b_{1}}}_{s-k}$ to $t$ or append
$\underbrace{\stack{a_{k}}{\bar{a}_{k}} \cdots
\stack{a_{k}}{\bar{a}_{k}}}_{s-k}$ to the end of $t$.  In these cases
the filling locations are $k$ and $0$, respectively.

\begin{proposition} \label{filling}
The map $F_{2,s}$ is well-defined on $\bigoplus_{i=0}^{s}B(i\Lt)$.
\end{proposition}

The proof of this proposition follows from the next three lemmas.

\begin{lemma} 
Suppose that $t\in\bigoplus_{k=0}^{s-1}B(k\Lt)$ has no subtableaux
$\stack{a_{i}}{b_{i}}\stack{a_{i+1}}{b_{i+1}}$ satisfying Criterion
\ref{eq:crit}.  Then either appending
$\underbrace{\stack{a_{k}}{\bar{a}_{k}} \cdots
\stack{a_{k}}{\bar{a}_{k}}}_{s-k}$ or prepending
$\underbrace{\stack{\bar{b}_{1}}{b_{1}} \cdots
\stack{\bar{b}_{1}}{b_{1}}}_{k}$ to $t$ will produce a tableau in
$\mathcal{T}(s)\setminus B(s\Lt)$.
\end{lemma}

\begin{proof} Suppose
$t=\stack{a_{1}}{b_{1}}\cdots\stack{a_{k}}{b_{k}}\in B(k\Lt)$ is as
above for $k<s$.  We will show that if prepending
$\stack{\bar{b}_{1}}{b_{1}} \cdots \stack{\bar{b}_{1}}{b_{1}}$ to $t$
does not produce a tableau in $\mathcal{T}(s)\setminus B(s\Lt)$, then
appending $\stack{a_{k}}{\bar{a}_{k}} \cdots
\stack{a_{k}}{\bar{a}_{k}}$ to $t$ will produce a tableau in
$\mathcal{T}(s)\setminus B(s\Lt)$.  There are two reasons we might not
be able to prepend $\stack{\bar{b}_{1}}{b_{1}} \cdots
\stack{\bar{b}_{1}}{b_{1}}$; $b_{1}$ may be unbarred, or we may have
$a_{1}<\bar{b}_{1}$.

First, suppose $b_{1}$ is unbarred.  If $b_{k}$ is also unbarred, then
$b_{k}$ is certainly less than $\bar{a}_{k}$, so we may append 
$\stack{a_{k}}{\bar{a}_{k}} \cdots \stack{a_{k}}{\bar{a}_{k}}$ to
$t$.  
Hence, suppose that $b_{k}$ is barred. We will show that $a_{k}$ is
unbarred
and $\bar{a}_{k}>b_{k}$.

We know that $t$ has a subtableau of the form $\stack{a_{i}}{b_{i}}
\stack{a_{i+1}}{b_{i+1}}$ such that $b_{i}$ is unbarred and $b_{i+1}$
is barred.  It follows that $a_{i}$ is unbarred, and therefore
$\bar{a}_{i}>b_{i}$.  Since our hypothesis states that
$\stack{a_{i}}{b_{i}} \stack{a_{i+1}}{b_{i+1}}$ must not satisfy
Criterion~\ref{eq:crit}, this means that $\bar{a}_{i}>b_{i+1}$, which
is equivalent to $\bar{b}_{i+1}>a_{i}$.  Once again observing that
$\stack{a_{i}}{b_{i}}\stack{a_{i+1}}{b_{i+1}}$ does not satisfy
Criterion~\ref{eq:crit}, this implies that $\bar{b}_{i+1}>a_{i+1}$;
i.e., $a_{i+1}$ is unbarred, and $\bar{a}_{i+1}>b_{i+1}$.

We proceed with an inductive argument on $i<j<k$.  Suppose that 
$\stack{a_j}{b_j}\stack{a_{j+1}}{b_{j+1}}$ is a subtableau of $t$
such that $b_j$
and $b_{j+1}$ are barred, $a_j$ is unbarred, and $\bar{a}_j>b_j$.  By
reasoning identical 
to the above, we conclude that
\begin{eqnarray}
     & \bar{a}_j>b_{j+1}
\Rightarrow
\bar{b}_{j+1}>a_j
\Rightarrow
\bar{b}_{j+1}>a_{j+1}
\Rightarrow
\bar{a}_{j+1}>b_{j+1}, &  \label{eq:imp}
\end{eqnarray}
which once again means that $a_{j+1}$ is unbarred.

This inductively shows that $a_{k}$ is unbarred and
$\bar{a}_{k}>b_{k}$, so we may append $\stack{a_{k}}{\bar{a}_{k}}
\cdots \stack{a_{k}}{\bar{a}_{k}}$ to $t$ to get a tableau in
$\mathcal{T}(s)\setminus B(s\Lt)$.  By a symmetrical argument, we
conclude that if $a_{k}$ is barred, then we may prepend
$\stack{\bar{b}_{1}}{b_{1}} \cdots \stack{\bar{b}_{1}}{b_{1}}$ to $t$.

Now, suppose that $b_{1}$ is barred and $\bar{b}_{1}>a_{1}$.  This
means that $a_{1}$ is unbarred and $\bar{a}_{1}>b_{1}$, so the
induction carried out in equation \ref{eq:imp} applies.  It follows
that 
$a_{k}$ is unbarred and $\bar{a}_{k}>b_{k}$, so once again we may
append 
$\stack{a_{k}}{\bar{a}_{k}} \cdots \stack{a_{k}}{\bar{a}_{k}}$ to
$t$.  
Also, by a symmetrical argument, when $a_{k}$ is unbarred and
$b_{k}>\bar{a}_{k}$, 
we may prepend $\stack{\bar{b}_{1}}{b_{1}} \cdots
\stack{\bar{b}_{1}}{b_{1}}$ to $t$. 
Thus, when no subtableau of $t$ satisfy Criterion~\ref{eq:crit},
either appending $\stack{a_{k}}{\bar{a}_{k}} \cdots
\stack{a_{k}}{\bar{a}_{k}}$ or 
prepending $\stack{\bar{b}_{1}}{b_{1}} \cdots
\stack{\bar{b}_{1}}{b_{1}}$ to $t$ will 
produce a tableau in $\mathcal{T}(s)\setminus B(s\Lt)$.
\end{proof}

\begin{lemma} 
Any tableau
$t=\stack{a_{1}}{b_{1}}\cdots\stack{a_{k}}{b_{k}}\in\bigoplus_{k=0}^{s-1}B(k\Lt)$
has no more than two filling locations.  If it has two, they are
consecutive integers, and $F_{2,s}(t)$ doesn't depend on this choice.
\end{lemma}

\begin{proof} Let $0\le i_*\le k$ be minimal such that $i_*$ is a
filling location of $t$.  First assume that $0<i_*<k$.  This implies
the existence of a subtableau
$\stack{a_{i_*}}{b_{i_*}}\stack{a_{i_*+1}}{b_{i_*+1}}$ which satisfies
Criterion~\ref{eq:crit}.

Suppose that the first condition $b_{i_{*}}\leq\bar{a}_{i_{*}}\leq
b_{i_{*}+1}$ of Criterion~\ref{eq:crit} is satisfied, and consider
whether $i_{*}+1$ can be a filling location.  If
$b_{i_{*}+1}\leq\bar{a}_{i_{*}+1}\leq b_{i_{*}+2}$, we have
$$b_{i_{*}+1}\leq\bar{a}_{i_{*}+1}\leq\bar{a}_{i_{*}}\leq
b_{i_{*}+1},$$
which implies that $\bar{a}_{i_{*}}=\bar{a}_{i_{*}+1}=b_{i_{*}+1}$, so
that $t$ violates part \ref{legal2:3} of Criterion~\ref{legal2}.
Similarly, if $a_{i_{*}+1}\leq\bar{b}_{i_{*}+2}\leq a_{i_{*}+2}$, then
we have $$\bar{a}_{i_{*}+1}\leq\bar{a}_{i_{*}}\leq b_{i_{*}+1}\leq
b_{i_{*}+2}\leq\bar{a}_{i_{*}+1},$$
which also implies that
$\bar{a}_{i_{*}}=\bar{a}_{i_{*}+1}=b_{i_{*}+1}$, once again violating
part~\ref{legal2:3} of Criterion \ref{legal2}.  We conclude that if
$i_{*}$ is a filling location for which Criterion \ref{eq:crit} is
satisfied by $b_{i_{*}}\leq\bar{a}_{i_{*}}\leq b_{i_{*}+1}$, then
$i_{*}+1$ is not a filling location.  Furthermore, this argument shows
that $a_{i_{*}+1}>a_{i_{*}}$ or $b_{i_{*}+1}>\bar{a}_{i_{*}}$.  By the
partial ordering on our alphabet, it follows that $t$ has no other
filling locations.

Now, suppose for the filling location $i_{*}$, Criterion \ref{eq:crit}
is satisfied by $a_{i_{*}}\leq\bar{b}_{i_{*}+1}\leq a_{i_{*}+1}$.  The
condition $a_{i_{*}+1}\leq\bar{b}_{i_{*}+2}\leq a_{i_{*}+2}$ for
$i_{*}+1$ to be a filling location implies that
$$\bar{b}_{i_{*}+2}\leq\bar{b}_{i_{*}+1}\leq
a_{i_{*}+1}\leq\bar{b}_{i_{*}+2},$$
which as above leads to a violation of part~\ref{legal2:3} of Criterion
\ref{legal2}.  However, $i_{*}+1$ may be a filling location if
Criterion \ref{eq:crit} is satisfied by
$b_{i_{*}+1}\leq\bar{a}_{i_{*}+1}\leq b_{i_{*}+2}$.  Note that this
inequality implies that $a_{i_{*}+1}\leq\bar{b}_{i_{*}+1}$, which
tells us that $a_{i_{*}+1}=\bar{b}_{i_{*}+1}$.  Thus, choosing to
insert
$\stack{\bar{b}_{i_{*}+1}}{b_{i_{*}+1}}\cdots\stack{\bar{b}_{i_{*}+1}}
{b_{i_{*}+1}}$ between columns $i_{*}$ and $i_{*}+1$ or to insert
$\stack{a_{i_{*}+1}}{\bar{a}_{i_{*}+1}}\cdots\stack{a_{i_{*}+1}}{\bar{a}_{i_{*}+1}}$
between columns $i_{*}+1$ and $i_{*}+2$ does not change $F_{2,s}(t)$.
Since $i_{*}+1$ is a filling location with Criterion~\ref{eq:crit}
satisfied by $b_{i_{*}}\leq\bar{a}_{i_{*}}\leq b_{i_{*}+1}$, the
preceding paragraph implies that there are no other filling locations
in $t$.

Finally, suppose that $i_*=0$ is a filling location for $t$; i.e.,
$b_{1}$ is barred, $a_{1}$ is unbarred, and $\bar{b}_{1}\leq a_{1}$.
If $1$ is a filling location, Criterion \ref{eq:crit} is satisfied by
$b_{1}\leq\bar{a}_{1}\leq b_{2}$; otherwise, part~\ref{legal2:3} of
Criterion~\ref{legal2} is violated.  Put together, this means that
$\bar{a}_{1}=b_{1}$, so prepending $\stack{\bar{b}_{1}}{b_{1}}\cdots
\stack{\bar{b}_{1}}{b_{1}}$ to $t$ and inserting
$\stack{a_{1}}{\bar{a}_{1}} \cdots\stack{a_{1}}{\bar{a}_{1}}$ between
columns $1$ and $2$ results in the same tableau.  As in the above
cases, part~\ref{legal2:3} of Criterion \ref{legal2} and the partial
order on the alphabet prohibit any other filling locations.
\end{proof}

\begin{example}
Let $s=4$.  Then
\begin{displaymath}
t=
\left.
\begin{array}{|c|c|c|}
    \hline
    1 & 2& 3 \\
    \hline
    \overline{4} & \overline{2} & \overline{1} \\
    \hline
\end{array}
\right.  ,\quad
F_{2,4}(t)=
\left.
\begin{array}{|c|c|c|c|}
    \hline
    1 & 2 & 2& 3 \\
    \hline
    \overline{4} & \overline{2} & \overline{2} & \overline{1} \\
    \hline
\end{array}
\right. .
\end{displaymath}
While we could choose either column two or column three as the
filling 
location, either choice results in the same tableau.
\end{example}

\begin{lemma} 
If a filling location of
$t=\stack{a_{1}}{b_{1}}\cdots\stack{a_{k}}{b_{k}}\in
\bigoplus_{i=0}^{s-1}B(i\Lt)$ satisfies Criterion \ref{eq:crit} with
both inequalities, then $F_{2,s}(t)$ is independent of this choice.
\end{lemma}

\begin{proof} Suppose that $i_{*}\neq 0,k$ is a filling location for
$t$ where both parts of Criterion \ref{eq:crit} are satisfied.  This
means that the subtableau
$\stack{a_{i_{*}}}{b_{i_{*}}}\stack{a_{i_{*}+1}}{b_{i_{*}+1}}$
satisfies both $\bar{a}_{i_{*}}\leq b_{i_{*}+1}$ and $a_{i_{*}}\leq
\bar{b}_{i_{*}+1}$.  The latter of these implies that
$b_{i_{*}+1}\leq\bar{a}_{i_{*}}$, so we have
$\bar{a}_{i_{*}}=b_{i_{*}+1}$ and $\bar{b}_{i_{*}+1}=a_{i_{*}}$.
Thus, filling with either $\stack{a_{i_{*}}}{\bar{a}_{i_{*}}} \cdots
\stack{a_{i_{*}}}{\bar{a}_{i_{*}} }$ or
$\stack{\bar{b}_{i_{*}+1}}{b_{i_{*}+1}} \cdots
\stack{\bar{b}_{i_{*}+1}}{b_{i_{*}+1}}$ between columns $i_{*}$ and
$i_{*}+1$ results in the same tableau $F_{2,s}(t)$.
\end{proof}

\begin{example}
To illustrate, for
\begin{displaymath}
t=
\left.
\begin{array}{|c|c|c|}
    \hline
    2 & 3& 3 \\
    \hline
    \overline{4} & \overline{2} & \overline{1} \\
    \hline
\end{array}
\right.  \quad \text{we have} \quad
F_{2,s}(t)=
\left.
\begin{array}{|c|c|c|c|}
    \hline
    2 & 2 & 3& 3 \\
    \hline
    \overline{4} & \overline{2} & \overline{2} & \overline{1} \\
    \hline
\end{array}
\right. .
\end{displaymath}
\end{example}

By identifying $\mathcal{T}(s)$ with $\bigoplus_{i=0}^{s}B(i\Lt)$ via
the maps $D_{2,s}$ and $F_{2,s}$, we have defined a $U_q(D_n)$-crystal
with the decomposition \eqref{eq:decomp}, with vertices labeled by the
$2\times s$ tableaux of $\mathcal{T}(s)$.  The action of the Kashiwara
operators $\etil{i}$, $\ftil{i}$ for $i\in\{1,\ldots,n\}$ on this
crystal is defined in terms of the above bijection, explicitly
\begin{equation} \label{eq:drop fill}
\begin{split}
    \etil{i}(T) &= F_{2,s}(\etil{i}(D_{2,s}(T)))\\
    \ftil{i}(T) &= F_{2,s}(\ftil{i}(D_{2,s}(T))), 
\end{split}
\end{equation}
for $T\in \mathcal{T}(s)$, where the $\etil{i}$ and $\ftil{i}$ on the
right are the standard Kashiwara operators on $U_q(D_n)$-crystals
\cite{KN}.  In section~\ref{sec:B2s:tabs2} we will discuss the action
of $\etil{0}$ and $\ftil{0}$ on $\mathcal{T}(s)$, which makes
$\mathcal{T}(s)$ into an affine crystal called $\tbs$.

\begin{remark}
    Using the filling and dropping map we obtain a natural inclusion
    of $\mathcal{T}(s')$ into $\mathcal{T}(s)$ for $s'<s$.
    
\end{remark}

\begin{definition} \label{def:upsilon}
For $s'<s$, the map 
$\Upsilon_{s'}^{s}:\mathcal{T}(s')\hookrightarrow\mathcal{T}(s)$ is 
defined by $\Upsilon_{s'}^{s}=F_{2,s}\circ D_{2,s'}$.
\end{definition}

\subsection{Combinatorial construction of $\sigma$.}

Conjecture \ref{conj:pmsigma} defines $\chsig$ in terms of $\pm$
diagrams, which gives a definition of $\sigma$ via $D_{n-1}$ highest
weight tableaux.  In the case of $B^{2,s}$, we can also give a direct
combinatorial description of $\sigma(T)$ for any crystal vertex $T$.
As an auxilliary construction which will be useful in its own right
later on, we combinatorially describe
$\iota_{i_{1}}^{i_{2}}:B(i_{1}\Lt)\hookrightarrow B(i_{2}\Lt)$, the unique
crystal embedding that agrees with
$\iotil_{i_{1}}^{i_{2}}:\Sup(i_{1}\Lt)\hookrightarrow\Sup(i_{2}\Lt)$.

\begin{remark}
It will often be useful to identify $B(k\Lt)$ with 
its image in $\tbs$.  We will use the notation $T\in 
B(k\Lt)\subset\tbs$ to indicate this identification.
\end{remark}

Let $i_{1},i_{2}\in\{0,\ldots,s\}$, and $i_{1}<i_{2}$, so
$\iota_{i_{1}}^{i_{2}}$ denotes the embedding of $B(i_{1}\Lt)$ in
$B(i_{2}\Lt)$.  Let $T\in B(i_{1}\Lt)\subset \tbs$.  This embedding
can be combinatorially understood through the following observations:
\begin{remark} \label{rem:comb iota}
    \begin{tabular}{c}
    \end{tabular}
    \begin{itemize}
	\item $\varphi_{k}(T)=\varphi_{k}(\iota_{i_{1}}^{i_{2}}(T))$
	and $\vare_{k}(T)=\vare_{k}(\iota_{i_{1}}^{i_{2}}(T))$ for
	$k=2,\ldots,n$;
	
	\item $D_{2,s}(\iota_{i_{1}}^{i_{2}}(T))$ has $(i_{2}-i_{1})$
	more columns than $D_{2,s}(T)$ (section \ref{sec:bij});
	
	\item $v(T)$ and $v(\iota_{i_{1}}^{i_{2}}(T))$ are in the
	same stratum, where $v(T)$ is the branching component vertex
	containing $T$.

    \end{itemize}
\end{remark}
In other words, we know that $T$ has a maximal $a$-configuration of
size $s-i_{1}$ (section \ref{sec:bij}), and has completely reduced
form $\Tred$ (Definition \ref{def:redtab}).  Furthermore, let $r_{1T}$
be the number of $1$'s in the first row of $D_{2,s}(T)$ to the left of
a null-configuration, let $r_{2T}$ be the size of the
null-configuration therein, (possibly $r_{2T} = 0$), and let $r_{3T}$
be the number of $\bar{1}$'s to the right of the null-configuration.
We then have $t_{1T} = r_{1T}+r_{2T}$ and $t_{2T} = r_{2T}+r_{3T}$,
and the stratum of $v(T)$ and $v(\iota_{i_{1}}^{i_{2}}(T))$ is
$t_{1T}-t_{2T}$.

\begin{remark}
    Observe that $t_{1T}$ and $t_{2T}$ are the number of $+$'s and
    $-$'s, respectively, in the $\pm$ diagram associated with the
    $D_{n-1}$ highest weight tableaux of $v(T)$, but these numbers can
    be extracted directly from any tableau, not just one that is
    $D_{n-1}$ highest weight.

\end{remark}

We wish to construct a tableau $S$ with an $a$-configuration of size
$s-i_{2}$ such that $S^{\#}=\Tred$ and $t_{1S}-t_{2S}=t_{1T}-t_{2T}$.
Based on properties of the height $2$ type $D$ sliding algorithm of
section \ref{sec:slide}, these conditions can only be satisfied when
$t_{jS}=t_{jT}+(i_{2}-i_{1})$ for $j=1,2$.

We can calculate $\iota_{i_{1}}^{i_{2}}(T)$ by the following
algorithm:
\begin{algorithm}\label{alg:iota}\mbox{}
\begin{enumerate}
    \item Remove the $a$-configuration of size $s-i_{1}$ from $T$ and
    slide it to get a $2\times i_{1}$ tableau.

    \item Remove the $1$'s, $\bar{1}$'s and the null-configuration 
    from the result to get a skew tableau of shape 
    $(i_{1},i_{1}-t_{2T})/(t_{1T})$.
    
    \item Using the type $D$ sliding algorithm, produce a skew tableau
    of shape
   
    $$((i_{2}),(i_{2})-(t_{2T}+(i_{2}-i_{1}))/(t_{1T}+(i_{2}-i_{1})))$$.

    \item Fill this tableau with $1$'s, $\bar{1}$'s, and a 
    null-configuration so that the result is a $2\times i_{2}$ tableau.
    
    \item Use the height $2$ fill map $F_{2,s}$ (section
    \ref{sec:bij}) to insert $s-i_{2}$ columns into the tableau.
\end{enumerate}
\end{algorithm}
This produces the unique tableau satisfying the three properties of
Remark~\ref{rem:comb iota}.

\begin{example}
Let 
$$
T=
\begin{array}{|c|c|c|c|c|c|c|}
    \hline
    1 & 1 & 2 & 2 & 2 & \bar{3} & \bar{2}  \\
    \hline
    2 & 2 & 3 & \bar{2} & \bar{2} & \bar{2} & \bar{1}  \\
    \hline
     
\end{array}
\in B(5\Lt)\subset \tilde{B}^{2,7}.
$$
Running through the steps of our algorithm (using relation (2) of
section \ref{sec:plac} for step (3)) gives us
\begin{enumerate}
    \item 
    $
    \phantom{\iota_{5}^{6}(T)=}\:        \phantom{\iota(T)=}   
    \begin{array}{|c|c|c|c|c|}
	\hline
	1 & 1 & 2 & \bar{3} & \bar{2}  \\
	\hline
	2 & 2 & 3 & \bar{2} & \bar{1}  \\
	\hline
     
    \end{array}
    $
    \vspace{.1in}

    \item 
    $
    \phantom{\iota_{5}^{6}(T)=}\:    \phantom{\iota(T)=}
    \begin{array}{|c|c|c|c|c|c}
	\cline{3-5}
	 \multicolumn{1}{c}{} & & 2 & \bar{3} & \bar{2}  \\
	\hline
	2 & 2 & 3 & \bar{2} &  \multicolumn{1}{c}{} \\
	\cline{1-4}
     
    \end{array}
    $
    \vspace{.1in}

    \item 
    $
    \phantom{\iota_{5}^{6}(T)=}\:	\phantom{\iota(T)=}
    \begin{array}{|c|c|c|c|c|c|}
	\cline{4-6}
	 \multicolumn{2}{c}{}  &  & 3 & \bar{3} & \bar{2}  \\
	\hline
	2 & 2 & 3 & \bar{3} &  \multicolumn{2}{c}{}  \\
	\cline{1-4}
     
    \end{array}
    $    
    \vspace{.1in}
    
    \item 
    $
    \phantom{\iota_{5}^{6}(T)=}\:	\phantom{\iota(T)=}
    \begin{array}{|c|c|c|c|c|c|}
	\hline
	1 & 1 & 1 & 3 & \bar{3} & \bar{2}  \\
	\hline
	2 & 2 & 3 & \bar{3} & \bar{1} & \bar{1} \\
	\hline
     
    \end{array}
    $    
    \vspace{.1in}
    
    \item 
    $
    \phantom{\iota(T)=}
    \iota_{5}^{6}(T)=
    \begin{array}{|c|c|c|c|c|c|c|}
	\hline
	1 & 1 & 1 & 3 & 3 & \bar{3} & \bar{2}  \\
	\hline
	2 & 2 & 3 & \bar{3} & \bar{3} & \bar{1} & \bar{1} \\
	\hline
     
    \end{array}
    \in B(6\Lt)\subset \tilde{B}^{2,7}.
    $ 
\end{enumerate}
\end{example}

We can also define a map $\iota_{i_{1}}^{i_{2}}:B(i_{1}\Lt)\to
B(i_{2}\Lt)\cup\{0\}$ for $i_{2}<i_{1}$ by
$$
    \iota_{i_{1}}^{i_{2}}(T)=\begin{cases}
    (\iota_{i_{2}}^{i_{1}})^{-1}(T) & \textrm{ if
    }T\in\iota_{i_{2}}^{i_{1}}(B(j\Lt)), \\
    0 & \textrm{ otherwise }.
\end{cases}
$$
Reversing Algorithm \ref{alg:iota} makes this map explicit.  Lastly,
we define $\iota_{i}^{i}$ to be the identity map on $B(i\Lt)$, so
$\iota_{i_{1}}^{i_{2}}$ is defined for all
$i_{1},i_{2}\in\{0,\ldots,s\}$.

\begin{definition}
    \label{def:bcdual}
    Let $T\in B(k\Lt)\subset\tbs$ be a tableau whose branching
    component vertex is in stratum $j$.  Then $T^{\bcdual}$ is the
    tableau in $B(k\Lt)$ with the same completely reduced form as $T$
    whose branching component vertex is in stratum $-j$.
\end{definition}

Alternatively, we can interpret this involution as follows.  Let $T\in
B(v)$, where $v$ is a branching component vertex with $D_{n-1}$
highest weight tableau $u$, and whose complementary vertex is $v'$.
(Recall from section \ref{sec:KR:BC} that the complementary vertex of
$v$ is defined to be the vertex in the opposite stratum from $v$ that
is associaated with the same partition as $v$.)  Let
$i_{1},\ldots,i_{m}$ be a sequence such that
$T=\ftil{i_{1}}\cdots\ftil{i_{m}}u$.  Then
$T^{\bcdual}=\ftil{i_{1}}\cdots\ftil{i_{m}}u'$, where $u'$ is the
$D_{n-1}$-highest weight tableau of $B(v')$.  Alternatively, this map
is the composition of $\dual$ with the ``local $\dual$'' map, which
applies only to the tableaux in $B(v)$ viewed as a $D_{n-1}$-crystal.

We now define $\sigma(T)$ in terms of the above combinatorial
operations.

\begin{definition}
    \label{def:combsig}
    Suppose $T\in B(k\Lt)\subset\tbs$, and $\ell$ be minimal such that
    $\iota_{k}^{s}(T)\in \iota_{\ell}^{s}(B(\ell\Lt))$ (i.e., $\ell$
    is the first part of the partition corresponding to the highest
    weight of the $D_{n-1}$-crystal containing $T$).  Then
    \begin{equation}\label{eq:sigma}
	\sigma(T)=\iota_{k}^{s+\ell-k}(T^{\bcdual}).
    \end{equation}

\end{definition}

\begin{remark}
    Since $\iota_{i_{1}}^{i_{2}}$ commutes with $\bcdual$, we also
    have $$\sigma(T)=(\iota_{k}^{s+\ell-k}(T))^{\bcdual}.$$

\end{remark}

\subsection{Properties of $\ftil{0}$ and $\etil{0}$.}

This explicit combinatorial construction immediately gives us useful
information about this crystal, such as the following lemma.

\begin{lemma} \label{lem:uk} For $k=0,1,\ldots,s$, let $u_{k}$ denote
the $D_{n}$ highest weight vector of the classical component
$B(k\Lt)\subset\tbs$.  Then
$$
\ftil{0}(u_{k})=\begin{cases}
    u_{k+1} & \textrm{if }k<s, \\
    0 & \textrm{if }k=s.
\end{cases}
$$
\end{lemma}

\begin{proof} Observe that
    $$
    u_{k}= \underbrace{
    \stack{1}{2}
    \cdots
    \stack{1}{2}}_{k}
    \underbrace{ 
    \stack{1}{\bar{1}}
    \cdots
    \stack{1}{\bar{1}}}_{s-k};
    $$
    We wish to calculate $\ftil{0}(u_{k})=\sigma\ftil{1}\sigma(u_{k})$.
    
    Note that $u_{k}\notin\iota_{k-1}^{k}(B((k-1)\Lt))$, so $\ell=k$ 
    in the combinatorial definition of $\sigma$ above.  It follows 
    that $\sigma(u_{k})=\iota_{k}^{s}(u_{k}^{\bcdual})$, which is
    \begin{equation*}
	\iota_{k}^{s}(u_{k}^{\bcdual}) 
	=
	\iota_{k}^{s}\Bigg( 
	\underbrace{
	\stack{1}{\bar{1}}
	\cdots
	\stack{1}{\bar{1}}}_{s-k}
	\underbrace{ 
	\stack{2}{\bar{1}}
	\cdots
	\stack{2}{\bar{1}}}_{k}\Bigg)
	=
	\framebox{$\emptyset_{s-k}$}
	\underbrace{ 
	\stack{2}{\bar{1}}
	\cdots
	\stack{2}{\bar{1}}}_{k},
    \end{equation*}
    where $\framebox{$\emptyset_{i}$}$ denotes a null-configuration of
    size $i$ (see Definition~\ref{def:null}).  
    If $k=s$, $\ftil{1}$ kills this tableau, as claimed in
    the second case of the Lemma.  Otherwise, acting by $\ftil{1}$
    will decrease the size of the null-configuration by $1$ and add
    another $\stack{2}{\bar{1}}$ to the columns on the right.  It
    follows that $\iota_{s}^{k}$ kills this tableau, but
    $\iota_{s}^{k+1}$ does not, so now $\ell=k+1$ in the combinatorial
    definition of $\sigma$.  Thus,
    \begin{equation*}
	\sigma\ftil{1}\sigma(u_{k})=\iota_{s}^{k+1}\Bigg(
	\underbrace{
	\stack{1}{2}
	\cdots
	\stack{1}{2}	}_{k+1}
	\framebox{$\emptyset_{s-k-1}$}
	\Bigg)
	=
	\underbrace{
	\stack{1}{2}
	\cdots
	\stack{1}{2}}_{k+1}
	\underbrace{ 
	\stack{1}{\bar{1}}
	\cdots
	\stack{1}{\bar{1}}}_{s-k-1}
	=u_{k+1}.
    \end{equation*}
\end{proof}

\begin{corollary} \label{cor:uk} Let $u_{k}$ be as above for $k>0$.  Then
    $$
    \etil{0}(u_{k})=	u_{k-1}     .
    $$
\end{corollary}

A similar combinatorial analysis can be carried out on lowest weight
tableaux to show that $\ftil{0}(u_{k}^{\dual})=u_{k-1}^{\dual}$ and
$\etil{0}(u_{k}^{\dual})=u_{k+1}^{\dual}$ for appropriate values of $k$.
Since $u_{0}=u_{0}^{\dual}$, this gives us the following Corollary:

\begin{corollary} \label{cor:phiuk} For $D_{n}$ highest weight vectors
$u_{k}$ and $D_{n}$ lowest weight vectors $u_{k}^{\dual}$, we have
$$
\varphi_{0}(u_{k})=\vare_{0}(u_{k}^{\dual})=s-k
\qquad \text{and} \qquad
\varphi_{0}(u_{k}^{\dual})=\vare_{0}(u_{k})=s+k
$$
\end{corollary}	
	    
        \section{Proof of Theorem \ref{thm:b2s}}
        \label{sec:B2s:r2}
        \subsection{Overview} To show that $\tbs$ is perfect, it must be shown that
all criteria of Definition~\ref{def:perfect} are satisfied with $\ell=s$.  
We have taken part 3 of Definition~\ref{def:perfect} as part of our hypothesis 
for Theorem \ref{thm:b2s}, so we do not attempt to prove this here.

Part 2 of Definition~\ref{def:perfect} is satisfied by simply noting
that $\la=\om_2=s\Lambda_2-2s\Lz$ is a weight in $P_{\cl}$ such that
$B_{\la}=\{u_{s}\}$ contains only one crystal vertex and all other
vertices in $\tbs$ have ``lower'' weights, in the ordering designated
by part 2 of Definition~\ref{def:perfect}.

In section \ref{sec:per con}, we show that $\tbs\otimes \tbs$ is
connected, proving that part 1 of Definition~\ref{def:perfect} is
satisfied.  Parts 4 and 5 of Definition~\ref{def:perfect} will be
dealt with simultaneously in sections~\ref{sec:sur} and~\ref{sec:inj}
by examining the levels of tableaux combinatorially.  We will see that
the level of a generic tableau is at least $s$ and the tableaux of
level $s$ are in bijection with the level $s$ weights.  In
section~\ref{sec:unique} we show that $\tbs$ is the unique affine
crystal satsifying the properties of Conjecture~\ref{conj:Brs} thereby
proving Theorem~\ref{thm:b2s}.

\subsection{Connectedness of $\tbs$} \label{sec:per con}

\begin{lemma}[Part 1 of Definition~\ref{def:perfect}]
The crystal $\tbs\otimes \tbs$ is connected as a $D_{n}^{(1)}$-crystal.
\end{lemma}

\begin{proof} (This proof is very similar to that of~\cite[Proposition
5.1]{Koga:1999}.)  For $k=0,1,\ldots,s$, let $u_{k}$ denote the
highest weight vector of the classical component $B(k\Lt)\subset\tbs$,
as in Lemma \ref{lem:uk}.  We will show that an arbitrary vertex
$b\otimes b'\in\tbs\otimes \tbs$ is connected to $u_{0}\otimes u_{0}$.

We know that for some $j\in\{0,\ldots,s\}$, we have $b'\in B(j\Lt)$. 
Then for some pair of sequences $i_{1}, i_{2},\ldots,i_{p}$ (with
entries in $\{1,\ldots,n\}$) and $m_{1}, m_{2},\ldots,m_{p}$ (with
entries in $\mathbb{Z}_{>0}$) and some $b^{1}\in\tbs$, we have
$\etil{i_{1}}^{m_{1}}\etil{i_{2}}^{m_{2}}\cdots\etil{i_{p}}^{m_{p}}(b\otimes 
b')=b^{1}\otimes u_{j}$.

By Corollary \ref{cor:phiuk}, $\varphi_{0}(u_{j})=s-j$, so if
$\vare_{0}(b^{1})\leq s-j$, Lemma \ref{lem:uk} tells us that
$\etil{0}^{j}(b^{1}\otimes u_{j})=b^{1}\otimes u_{0}$.  If
$\vare_{0}(b^{1})=r>s-j$, then $\etil{0}^{r-s+j}(b^{1}\otimes
u_{j})=b^{2}\otimes u_{j}$, where $\vare_{0}(b^{2})=r-(r-s+j)=s-j$, so
$\etil{0}^{j}(b^{2}\otimes u_{j})=b^{2}\otimes u_{0}$.  In either
case, our arbitrary $b\otimes b'$ is connected to an element of the
form $b''\otimes u_{0}$.

Let $j'$ be such that $b''\in B(j'\Lt)$.  Since $u_{0}$ is the unique
element of $B(0)$, the crystal for the trivial representation of
$U_q(D_n)$, we know that $B(j'\Lt)\otimes B(0)\simeq B(j'\Lt)$. 
Therefore, $b''\otimes u_{0}$ is connected to $u_{j'}\otimes u_{0}$. 
Finally, we note that $\varphi_{0}(u_{0})=s<s+j'=\vare_{0}(u_{j'})$
for $j'\neq 0$, so $\etil{0}^{j'}(u_{j'}\otimes
u_{0})=\etil{0}^{j'}(u_{j'})\otimes u_{0}=u_{0}\otimes u_{0}$,
completing the proof.
\end{proof}

\subsection{Preliminary observations}\label{per15}
We first make a few observations. 

\begin{proposition} \label{prop:eps1 eps0}
    Let $T\in B(k\Lt)\subset\tbs$, and set $T^{(m)}=\iota_{k}^{m}(T)$
    for $m\in\{\ell,\ell+1,\ldots,s\}$, where $\ell$ is minimal such
    that $\iota_{k}^{\ell}(T)\neq 0$.  If $\ell\neq s$, we have for
    $\ell\leq m<s$
    \begin{eqnarray*}
	&&\vare_{1}(T^{(m+1)})=\vare_{1}(T^{(m)})+1 
	\textrm{ and } 
	\vare_{0}(T^{(m+1)})=\vare_{0}(T^{(m)})-1, \\
	&&\varphi_{1}(T^{(m+1)})=\varphi_{1}(T^{(m)})+1 
	\textrm{ and } 
	\varphi_{0}(T^{(m+1)})=\varphi_{0}(T^{(m)})-1.
    \end{eqnarray*}
\end{proposition}

\begin{proof}
    Let $\ell\leq m\leq s-1$, so $\iota_{m}^{m+1}$ is defined.  We
    first consider the difference between the reduced $1$-signatures
    of $D_{2,s}(T^{(m)})$ and $D_{2,s}(T^{(m+1)})=
    D_{2,s}(\iota_{m}^{m+1}(T^{(m)}))$, since the action of $\etil{1}$
    on these tableaux is defined by the action of the classical
    $\etil{1}$ on their image under $D_{2,s}$.  Let $-^{M}+^{P}$ be
    the reduced $1$-signature of $D_{2,s}(T^{(m)})$.
    Let $r_{1}$ denote the number of $1$'s in $D_{2,s}(T^{(m)})$,
    $r_{3}$ the number of $\bar{1}$'s, $r_{2}$ the size of the
    null-configuration, and $t_{1}=r_{1}+r_{2}$, $t_{2}=r_{2}+r_{3}$.
    Then there is a contribution $-^{r_{2}}+^{r_{2}}$ to the
    1-signature from the null-configuration, and the remaining $-$'s
    and $+$'s come from $1$'s with a letter greater than $2$ below
    them and $\bar{1}$'s with a letter less than $\bar{2}$ above them,
    respectively.

    We now have two cases.
    If $t_{1}+t_{2}\geq s$, $\iota_{m}^{m+1}$ simply increases the
    size of the null-configuration in $D_{2,s}(T^{(m)})$ by $1$.  It
    follows that the reduced $1$-signature of $D_{2,s}(T^{(m+1)})$ is
    $-^{M+1}+^{P+1}$, as we wished to show.  On the other hand, if
    $t_{1}+t_{2}<s$, after step $(2)$ of Algorithm \ref{alg:iota} for
    $\iota_{m}^{m+1}$ we have a tableau of shape
    $(m,m-r_{3})/(r_{1})$.  In step (3), we slide this into shape
    $(m+1,m+1-(r_{3}+1))/(r_{1}+1)$.  We claim that the rightmost
    ``uncovered'' letter in the second row of this tableau is greater
    than $2$ and the leftmost ``unsupported'' letter in the first row
    is less than $\bar{2}$.  As observed in the preceding paragraph,
    this implies that after refilling the empty spaces as in step
    $(4)$ of Algorithm \ref{alg:iota} the reduced $1$-signature of our
    tableau is $-^{M+1}+^{P+1}$ in this case as well.
    
    Let us first consider the leftmost ``unsupported'' letter.  After 
    step $(2)$, our tableau is of the form
    $$
    \begin{array}{|c|c|c|c|c|c|c|c|c|}
        \cline{4-9}
         \multicolumn{2}{c}{} & & a_{r_{1}+1} & \cdots & a_{m-r_{3}} & 
         a_{m-r_{3}+1} & 
         \cdots & a_{s}  \\
        \cline{1-9}
        b_{1} & \cdots & b_{r_{1}} & b_{r_{1}+1} & \cdots & b_{m-r_{3}} &  
        \multicolumn{3}{c}{}  \\
        \cline{1-6}
    \end{array}
    $$
    and its column word is unchanged by the slide
    $$
    \begin{array}{|c|c|c|c|c|c|c|c|c|c|}
	\cline{4-5}	\cline{7-10}
	 \multicolumn{2}{c}{} & & a_{r_{1}+1} & \cdots & & a_{m-r_{3}} & 
	 a_{m-r_{3}+1} & 
	 \cdots & a_{s}  \\
	\cline{1-10}
	b_{1} & \cdots & b_{r_{1}} & b_{r_{1}+1} & \cdots & b_{m-r_{3}} &  
	\multicolumn{4}{c}{}  \\
	\cline{1-6}
    \end{array},
    $$    
    so we have $a_{m-r_{3}}<b_{m-r_{3}}\leq\bar{2}$.
    
    The second row of this tableau has $m-r_{3}$ boxes just as it did
    before sliding, so the boxes in the bottom row will never be
    moved.  It follows that this sliding procedure only changes
    L-shaped subtableaux into $\Le$-shapes \big(i.e., \quad $
    \begin{array}{|c|c|c}
        \cline{1-1}
         &  \multicolumn{2}{c}{}  \\
        \cline{1-2}
         &  &   \\
        \cline{1-2}
    \end{array}
    $
    into
    $
    \begin{array}{c|c|c|}
	\cline{3-3}
	   \multicolumn{1}{c}{} & & \\
	\cline{2-3}
	 &  &   \\
	\cline{2-3}
    \end{array}
    $
    \; \big) and never involves any $\Gamma$- or $\Lee$-shapes.
    According to the Lecouvey $D$-equivalence relations from
    section~\ref{sec:plac}, such moves can only be made when the
    letters in the bottom row are strictly greater than $2$.  
    Specifically, in relations $(3)$ and $(4)$, the letter which is 
    ``uncovered'' is either $n$ or $\bar{n}$, while in relations $(1)$ 
    and $(2)$ only the second case of each relation applies.  This 
    proves our claim, and thus the first half of the proposition.

Since $\etil{0}=\sigma\circ \etil{1}\circ\sigma$, we can derive the statements
about $\vare_0$ and $\varphi_0$ from the corresponding statements about
$\vare_1$ and $\varphi_1$. More precisely, $\vare_0(T)=\vare_1(\sigma(T))$
and $\varphi_0(T)=\varphi_1(\sigma(T))$ and by \eqref{eq:sigma} we have
\begin{equation*}
\sigma(T^{(m)})=\bigl(\iota_m^{s+\ell-m}\circ\iota_k^m(T)\bigr)^{\bcdual}
=\iota_k^{s+\ell-m}(T^{\bcdual}).
\end{equation*}
Hence
\begin{multline*}
\vare_0(T^{(m+1)})=\vare_1(\sigma(T^{(m+1)}))
=\vare_1(\iota_k^{s+\ell-m-1}(T^{\bcdual}))\\
=\vare_1(\iota_k^{s+\ell-m}(T^{\bcdual}))-1
=\vare_1(\sigma(T^{(m)}))-1=\vare_0(T^{(m)})-1.
\end{multline*}
A similar computation can be carried out for $\varphi_0$.
\end{proof}

\begin{corollary} \label{cor:lz plus lo}
    Given the above hypotheses, we have
    \begin{equation*}
	 \lan h_{0}+h_{1},\vare(T^{(s)})\ran= \lan
	 h_{0}+h_{1},\vare(T^{(s-1)})\ran= \cdots 
	 =\lan h_{0}+h_{1},\vare(T^{(\ell)})\ran\neq 0.
    \end{equation*}
\end{corollary}

The following observation is an immediate consequence of Remark 
\ref{rem:comb iota}:

\begin{corollary} \label{cor:hi}
    For $i=2,\ldots,n$,
    \begin{equation*}
	\lan h_{i},\vare(T^{(s)})\ran= \lan h_{i},\vare(T^{(s-1)})\ran= 
	\cdots  = \lan h_{i},\vare(T^{(\ell)})\ran.
    \end{equation*}
\end{corollary}

\begin{lemma}\label{lem:ups}
    The map
    $\Upsilon_{s-1}^{s}:\mathcal{T}(s-1)\hookrightarrow\mathcal{T}(s)$
    sending each summand $B(k\Lt)\subset\tilde{B}^{2,s-1}$ to
    $B(k\Lt)\subset\tbs$ for $k=0,\ldots,s-1$ increases the level of
    all tableaux by exactly $1$.
\end{lemma}

\begin{proof} 
    Let $T\in\tilde{B}^{2,s-1}$.
    We have $\varphi_{i}(\Upsilon_{s-1}^{s}(T))=\varphi_{i}(T)$ for
    $i=1,\ldots,n$.  To calculate the change in $\varphi_{0}(T)$, we
    must consider the difference between
    $\varphi_{1}(\sigma_{s-1}(T))$ and
    $\varphi_{1}(\sigma_{s}(\Upsilon_{s-1}^{s}(T)))$.  By our
    descriptions of maps on crystals, we have
    $\sigma_{s}(\Upsilon_{s-1}^{s}(T))=\Upsilon_{s-1}^{s}(\iota_{j}^{j+1}(\sigma_{s-1}(T)))$,
    where $j$ is determined by $\sigma_{s-1}(T)\in
    B(j\Lt)\subset\tilde{B}^{2,s-1}$.  By Proposition \ref{prop:eps1
    eps0},
    $\varphi_{1}(\iota_{j}^{j+1}(\sigma_{s-1}(T)))=\varphi_{1}(\sigma_{s-1}(T))+1$.
\end{proof}

\subsection{Surjectivity} \label{sec:sur}

Given a weight $\la\in(P_{\cl}^{+})_{s}$, we construct a tableau
$T_{\la}\in\tbs$ such that $\vare(T_{\la})=\varphi(T_{\la})=\la$.
This amounts to constructing $T_{\la}$ so that its reduced
$i$-signature is $-^{\vare_{i}(T_{\la})}+^{\vare_{i}(T_{\la})}$.  Note
that such a tableau is invariant under the $\dual$-involution, so its
symmetry allows us to define it beginning with the middle columns,
proceeding outwards.

For $i=0,\ldots,n$, let $k_{i}=\lan h_{i},\la\ran$.  We first
construct a tableau $T_{\la'}$ corresponding to the weight
$\la'=\sum_{i=2}^{n}k_{i}\Lambda_{i}$.  We begin with the middle
$k_{n-1}+k_{n}$ columns of $T_{\la'}$.  If $k_{n-1}+k_{n}$ is even and
$k_n\geq k_{n-1}$, these columns of $T_{\la'}$ are
\begin{displaymath}
\underbrace{
\stack{n-2}{n-1}
\cdots
\stack{n-2}{n-1}
}_{k_{n-1}}
\underbrace{
\stack{n-1}{n}
\cdots
\stack{n-1}{n}
}_{(k_{n}-k_{n-1})/2}
\underbrace{
\stack{\overline{n}}{\overline{n-1}}
\cdots
\stack{\overline{n}}{\overline{n-1}}
}_{(k_{n}-k_{n-1})/2}
\underbrace{
\stack{\overline{n-1}}{\overline{n-2}}
\cdots
\stack{\overline{n-1}}{\overline{n-2}}
}_{k_{n-1}}
\end{displaymath}
If $k_{n-1}+k_{n}$ is odd and $k_{n}\geq k_{n-1}$, we have
\begin{displaymath}
\underbrace{
\stack{n-2}{n-1}
\cdots
\stack{n-2}{n-1}
}_{k_{n-1}}
\underbrace{
\stack{n-1}{n}
\cdots
\stack{n-1}{n}
}_{(k_{n}-k_{n-1}-1)/2}
\stack{\bar{n}}{n}
\underbrace{
\stack{\overline{n}}{\overline{n-1}}
\cdots
\stack{\overline{n}}{\overline{n-1}}
}_{(k_{n}-k_{n-1}-1)/2}
\underbrace{
\stack{\overline{n-1}}{\overline{n-2}}
\cdots
\stack{\overline{n-1}}{\overline{n-2}}
}_{k_{n-1}}
\end{displaymath}
In either case, if $k_{n}<k_{n-1}$, interchange $n$ with $\bar{n}$ and
$k_n$ with $k_{n-1}$ in the above configurations.

Next we put a configuration of the form
$$
\underbrace{
\stack{1}{2}
\cdots
\stack{1}{2}
}_{k_{2}}
\underbrace{
\stack{2}{3}
\cdots
\stack{2}{3}
}_{k_{3}}
\qquad
\cdots
\qquad
\underbrace{
\stack{n-3}{n-2}
\cdots
\stack{n-3}{n-2}
}_{k_{n-2}}
$$
on the left, and a configuration of the form
$$
\underbrace{
\stack{\overline{n-2}}{\overline{n-3}}
\cdots
\stack{\overline{n-2}}{\overline{n-3}}
}_{k_{n-2}}
\underbrace{
\stack{\overline{n-3}}{\overline{n-4}}
\cdots
\stack{\overline{n-3}}{\overline{n-4}}
}_{k_{n-3}}
\qquad
\cdots
\qquad
\underbrace{
\stack{\bar{2}}{\overline{1}}
\cdots
\stack{\bar{2}}{\overline{1}}
}_{k_{2}}
$$
on the right, completing $T_{\la'}$.  Denote the set of tableaux
constructed by the procedure up to this point by $\T(s')$, where $s'$ 
denotes the level of $\la'$.

Observe that the reduced $1$-signature of $T_{\la'}$ is empty, so
$\lan h_{1},\varphi(T_{\la'})\ran=0$.  Furthermore, since $\la'$ has
the same number of $1$'s as $\bar{1}$'s, it is fixed by $\sigma$, so
$\lan h_{0},\varphi(T_{\la'})\ran=0$ as well.  Thus Proposition
\ref{prop:eps1 eps0} implies that $T_{\la'}\in \tbs_{\min}\cap
B(s'\Lt)\setminus\iota_{s'-1}^{s'}(B((s'-1)\Lt))$ as a subset of
$\tilde{B}^{2,s'}$.  Recall the embedding
$\Upsilon_{s'}^{s}:\tilde{B}^{2,s'} \hookrightarrow\tbs$ from
Definition~\ref{def:upsilon}.  Since $s'$ is the minimal $m$ for which
$\iota_{s'}^{m}(T_{\la'})\neq 0$, Lemma \ref{lem:ups} and its proof
tell us that $\vare_{0}(F_{2,s}(T_{\la'}))=s-s'$, where the fill map
$F_{2,s}$ inserts an $a$-configuration to increase the width of
$T_{\la'}$ to $s$.  By the same proposition the desired tableau is
$T_{\lambda}=\iota_{s'}^{s'+k_1}\circ F_{2,s}(T_{\lambda'})$.  We
denote by $\Mmin{s}$ the set of tableaux constructed by this
procedure.

\subsection{Injectivity}\label{sec:inj}

In this subsection we show that the tableaux in $\Mmin{s}$ are all the
minimal tableaux in $\tbs$.

We first introduce some useful notation.  Observe that any tableau $T$
can be written as $T=T_{1}T_{2}T_{3}T_{4}T_{5}$, where the block
$T_{i}$ has width $w_{i}$, and all letters in $T_{1}$ (resp.  $T_{5}$)
are unbarred or $\bar{n}$ in the second row (resp.  barred or $n$ in
the first row), all columns in $T_{2}$ (resp.  $T_{4}$) are of the
form $\stack{a}{\bar{b}}$ where $a<b\leq n-1$ (resp.  $b<a\leq n-1$),
and all columns in $T_{3}$ are of the form $\stack{a}{\bar{a}}$ for
some $a$.  Note that for a tableau in $B(s\Lt)\subset\tbs$ we have
$0\leq w_{3}\leq 1$.  Also note that $T_{2}$ and $T_{4}$ do not contain 
any $n$'s or $\bar{n}$'s.

\begin{theorem} \label{thm:injective}
    We have $\Mmin{s}=\tbs_{\min}$.
\end{theorem}

\begin{proof}
    Our proof is by induction on $s$.  For the base case, we have
    explicitly verified by exhaustion that the statement of the
    theorem is true for $s=0,1,2$.
    
    By our induction hypothesis,
    $\Mmin{s-1}=\tilde{B}^{2,s-1}_{\min}$.  By Lemma \ref{lem:ups},
    $\Upsilon_{s-1}^{s}$ increases the level of a tableau by $1$, and
    therefore
    $\Upsilon_{s-1}^{s}(\tilde{B}^{2,s-1}_{\min})=(\tbs\setminus
    B(s\Lt))_{\min}$.  By Corollaries \ref{cor:lz plus lo} and
    \ref{cor:hi}, $\iota_{s-1}^{s}$ does not change the level of a
    tableau, so it suffices to show that
    \begin{equation}     \label{eq:minimal}
    \T(s)=\tbs_{\min}\cap(B(s\Lt)\setminus\iota_{s-1}^{s}(B((s-1)\Lt))).
    \end{equation}
    By Lemmas \ref{lem:plus minus} and \ref{lem:mixed} below, if $T$ is a
    minimal tableau not in the image of $\iota_{s-1}^{s}$, it has
    $w_{2}=w_{4}=0$, and if $w_{3}=1$, then $T_{3}=\stack{n}{\bar{n}}$
    or $T_{3}=\stack{\bar{n}}{n}$.  By Lemma \ref{lem:unmixed},
    equation \eqref{eq:minimal} follows.
\end{proof}

The following lemmas are used in the proof of
Theorem~\ref{thm:injective}.  The proofs of these lemmas rely on
induction on $s$, as did the proof of Theorem~\ref{thm:injective}.
The base cases $s=0,1,2$ have been checked explicitly.

\begin{lemma}\label{lem:lb}
For all $t\in \tbs$ we have $\lan c,\varepsilon(t)\ran \ge s$ and 
$\lan c,\varphi(t)\ran \ge s$.
\end{lemma}

\begin{proof}
    By the induction hypothesis,
    $\tilde{B}^{2,s-1}_{\min}=\Mmin{s-1}$.

Observe that by Corollaries \ref{cor:lz plus lo} and \ref{cor:hi}
$\iota_{i}^{j}$ is level preserving, and by Lemma \ref{lem:ups}, the
map $\Upsilon_{s-1}^{s}$ increases the level of a tableau by one.  Our
induction hypothesis therefore allows us to assume that $t\in
B(s\Lt)\setminus \iota_{s-1}^{s}(B((s-1)\Lt))$.  Combinatorially, we
may characterize such tableaux as being those which are legal in the
classical sense and for which removing all $1$'s, $\bar{1}$'s, and
null configurations produces a tableau which is Lecouvey
$D$-equivalent to a tableau whose first row has width $s$.  This
characterization follows from the combinatorial description of
$\iota_{s-1}^{s}$ in Algorithm~\ref{alg:iota}.

We may further restrict our attention by the observation that if $T$
is minimal, then so is $T^{\dual}$.  We may therefore assume that the
stratum of the branching component vertex of $T$ is non-negative.
In particular, this means that $T$ has no more $\bar{1}$'s than
$1$'s.

Our approach is to consider the tableau $T'$ that
results from removing the leftmost column from $T$.  We will show that
if $T'$ is minimal, the level of $T$ exceeds the level of $T'$ by at
least $2$, and if $T'$ is not minimal, the level of $T$ is at least as
great as the level of $T'$.

First consider the case when $T'$ is minimal.  Since $T$ is assumed to
be such that removing all $1$'s and $\bar{1}$'s produces a tableau
which is Lecouvey $D$-equivalent to a tableau whose first row has
width $s$, it is the case that removing all $1$'s and $\bar{1}$'s from
$T'$ produces a tableau which is Lecouvey $D$-equivalent to a tableau
whose first row has width $s-1$.  The minimal tableaux of
$\tilde{B}^{2,s-1}$ with this property are precisely $\T(s-1)$.  By
properties of $\T(s-1)$, we know that $\varphi_{0}(T')=0$, so
$\varphi_{0}(T)\geq\varphi_{0}(T')$.  Since our base case is $s=2$, we
know that the first column of $T$ is $\stack{a}{b}$, where $a$ and $b$ 
are both unbarred.  Observe that 
$\varphi(T)=\varphi(T')+2\Lambda_{b}+\textrm{non-negative weight}$ if
$b\neq n-1$ or  
$\varphi(T)=\varphi(T')+\Lambda_{n-1}+\Lambda_{n}+\textrm{non-negative 
weight}$ if $b=n-1$.  Hence, the level is increased by at least $2$.

Now suppose $T'$ is not minimal.  The level of the $i$-signatures
(that is to say, the level of the sum of the weights $\varphi_{i}(T')$
which depend on $i$-signatures) cannot have a net decrease for
$i=1,\ldots,n$, but there is now a possibility that
$\varphi_{0}(T)<\varphi_{0}(T')$.  We will show that when
$\varphi_{0}(T)<\varphi_{0}(T')$, the level of the $i$-signatures goes
up by at least $\varphi_{0}(T')-\varphi_{0}(T)$.

First, suppose $T$ has no $1$'s.  Then by one of our hypotheses, it
also has no $\bar{1}$'s, and is therefore fixed by $\sigma$: it
follows that $\varphi_{0}(T)=\varphi_{1}(T)$, so we may assume the 
upper-left entry of $T$ to be $1$.

We know that $\varphi_{0}(T)$ is equal to the number of $-$'s in the 
reduced $1$-signature of $\sigma(T)$.  Consider the following 
tableaux:
\begin{eqnarray*}
    T' & = & \underbrace{\stack{1}{b_{1}}\cdots\stack{1}{
    b_{m_{1}}}}_{m_{1}}
    \stack{a_{m_{1}+1} }{ b_{m_{1}+1}}\cdots\stack{a_{s-m_{2}} }{ b_{s-m_{2}}}
    \underbrace{\stack{a_{s-m_{2}+1} }{ \bar{1}}\cdots\stack{a_{s} }{ 
    \bar{1}}}_{m_{2}}\\
    T & = &  \stack{1 }{b}\underbrace{\stack{1}{b_{1}}\cdots\stack{1 }{ 
    b_{m_{1}}}}_{m_{1}}
    \stack{a_{m_{1}+1} }{ b_{m_{1}+1}}\cdots\stack{a_{s-m_{2}} }{ b_{s-m_{2}}}
    \underbrace{\stack{a_{s-m_{2}+1} }{\bar{1}}\cdots\stack{a_{s}}{ 
    \bar{1}}}_{m_{2}}\\
    \sigma(T') & = & \underbrace{\stack{1 }{ 
    b_{1}}\cdots\stack{1}{b_{m_{2}-1}}\stack{1 }{ 
    b'_{m_{2}}}}_{m_{2}}
    \stack{a'_{m_{2}+1} }{ b'_{m_{2}+1}}\cdots\stack{a'_{s-m_{1}} }{ 
    b'_{s-m_{1}}}
    \underbrace{\stack{a'_{s-m_{1}+1} }{ \bar{1}}\cdots\stack{a'_{s} }{ 
    \bar{1}}}_{m_{1}} \\
    \sigma(T) & = & \underbrace{\stack{1 }{ b}\stack{1 }{ b_{1}}\cdots\stack{1 }{ 
    b_{m_{2}-1}}}_{m_{2}}
    \stack{a''_{m_{2}} }{ b''_{m_{2}}}\cdots\stack{a''_{s-m_{1}-1} }{ 
    b''_{s-m_{1}-1}}
    \underbrace{\stack{a''_{s-m_{1}} }{ \bar{1}}\cdots\stack{a''_{s}}{ 
    \bar{1}}}_{m_{1}+1} .
\end{eqnarray*}
Note that our assumption that $m_{2}\leq m_{1}+1$ ensures that the 
absence of primes on $b,b_{1},\ldots,b_{m_{2}-1}$ is accurate.  

Let us consider all possible ways for the number of $-$'s in the
$1$-signature to be smaller for $\sigma(T)$ than for $\sigma(T')$.
The number of $1$'s is the same, so the only way this contribution
could be decreased is by having more $2$'s in the first $m_{2}$
letters of the bottom row of $\sigma(T)$.  This can only come about by
having $b=2$, and only one $-$ may be removed in this way.

The other possibility is for the number of $-$'s contributed by
$\bar{2}$'s to be decreased.  The only Lecouvey relation which removes
a $\bar{2}$ assumes the presence of a column $\stack{2}{\bar{2}}$,
which we disallow, as it is necessarily part of a null-configuration.
To decrease this contribution therefore requires an additional $+$ in
the $1$-signature of $\sigma(T)$ compared to that of $\sigma(T')$,
which will bracket one of the $-$'s from a $\bar{2}$.  The additional
$+$ may come from one of the additional $\bar{1}$'s, or from a $2$
that is ``pushed out'' from under the $1$'s at the beginning in the
case $b=2$.  Note that this second possibility is mutually exclusive
with having more $2$'s bracketing $1$'s at the beginning.

In any case, we see that $\varphi_{0}(T)-\varphi_{0}(T')\leq2$, and
that when this value is $2$, the first column of $T$ is $\stack{1}{ 2}$.
This column adds no $+$'s to the $i$-signatures, but does provide a
new $-$ in the $2$-signature.  Since $\Lambda_2$ is a level $2$ weight, the
level stays the same in this case.

If $\varphi_{0}(T)-\varphi_{0}(T')=1$ and the first column of $T$ is 
$\stack{1}{ 2}$, we in fact have a net increase in level.  If $b\neq 2$, 
the $i$-signature levels go up by at least 1, so still the total level 
cannot decrease.

\end{proof}

\begin{lemma} \label{lem:plus minus}
If $T\in (\tbs\cap B(s\Lt))_{\min}$, we have $w_{1}+w_{2}=w_{4}+w_{5}=\lfloor
s/2 \rfloor$ with $w_1,w_2,w_4$ and $w_5$ as defined in the beginning of this
subsection.
\end{lemma}

\begin{proof}
    We first establish that $\lan c,\varphi(T)\ran\geq
2w_{1}+2w_{2}+w_{3}$, and thus by $\dual$-duality, $\lan
c,\vare(T)\ran\geq w_{3}+2w_{4}+2w_{5}$ as well.  Recall that $0\leq 
w_{3}\leq 1$.

First, observe that every letter in the bottom row of $T_{1}$
contributes:
\begin{itemize}
    \item a $-$ to the
    reduced $a$-signature if $2\leq a\leq n-2$ is in the bottom row;

    \item a $-$ to both the $(n-1)$-signature and the $n$-signature if
    $n-1$ is in the bottom row;

    \item a $-$ to the $n$-signature (resp.  $(n-1)$-signature) if $n$
    (resp.  $\bar{n}$) is the bottom row.
\end{itemize}
Suppose $T_{1}$ has a column of the form $\stack{a}{b}$ with $b\neq
a+1$, or $b=\bar{n}$ and $a\neq n-1$.  For the $-$ in the
$a$-signature of $T$ contributed by this $a$ to be bracketed, we must
have a column of the form $\stack{a'}{a+1}$ to the left of this column
in $T_{1}$, with $a'<a$.  Applying this observation recursively, we
see that to bracket as many $-$'s as possible we must eventually have
a column of the form $\stack{1}{c}$ for some $c\neq 2$.  Note that in
the case of columns of the form $\stack{n-1}{n}$ (resp.
$\stack{n-1}{\bar{n}}$) the $-$ in the $n$-signature (resp.
$(n-1)$-signature) from $n-1$ that is not bracketed by the $+$ from
the $n$ (resp.  $\bar{n}$) below it cannot be bracketed, since $n$ and
$\bar{n}$ may not appear in the same row.

Now, consider a column $\stack{a}{\bar{b}}$ in $T_{2}$, so we have
$a<b$, and thus also $\bar{a}>\bar{b}$.  Recall that $T_{2}$ has no
$n$'s or $\bar{n}$'s, so $b\leq n-1$.  This column contributes $-$'s
to the $a$-signature and the $(b-1)$-signature of $T$.  In this case,
these $-$'s may be bracketed.  Due to the conditions that the rows and
columns of $T$ are increasing, the $-$ from the $a$ can only be
bracketed by an $a+1$ in the bottom row of $T_{1}$ and the $-$ from
the $\bar{b}$ can only be bracketed by a $b$ in the bottom row of
$T_{1}$.  Furthermore, the letter above these must be strictly less
than $a$ and $b-1$, respectively.  By the reasoning in the previous
paragraph, we see that to bracket every $-$ engendered by the column
$\stack{a}{\bar{b}}$ we must have two columns of the form
$\stack{1}{a'}$, with each $a'\neq 2$.

If $w_{3}=1$, $T$ has a column of the form $\stack{a}{\bar{a}}$.  We
have two cases; $2\leq a\leq n-1$, and $a=n$ (resp.  $a=\bar{n}$).  In
the first case, we have a $-$ in the $(a-1)$-signature from the
$\bar{a}$ in this column.  Because of the prohibition against
configurations of the form $\stack{a}{}\stack{a}{\bar{a}}$, this $-$
can only be bracketed by a $+$ from an $a$ in the bottom row of
$T_{1}$.  Therefore, this column engenders another column of the form
$\stack{1}{a'}$.  In the case of $a=n$ (resp.  $a=\bar{n}$), we have a
$-$ in the $(n-1)$-signature (resp.  $n$-signature) which cannot be
bracketed.

To bracket the maximal number of $-$'s (i.e., to minimize $\lan
c,\varphi(T)\ran$) we see that unless $T_{3}=\stack{n}{\bar{n}}$ or
$T_{3}=\stack{\bar{n}}{n}$, we must have
\begin{equation}    \label{eq:t1t2t3}
    T_{1}T_{2}T_{3}=
    \underbrace{
    \stack{1}{b_{1}}
    \cdots
    \stack{1}{b_{2w_{2}+w_{3}}}
    }_{2w_{2}+w_{3}}
    \underbrace{
    \stack{a_{2w_{2}+w_{3}+1}}{b_{2w_{2}+w_{3}+1}}
    \cdots
    \stack{a_{w_{1}}}{b_{w_{1}}}
    }_{w_{1}-(2w_{2}+w_{3})}
    \underbrace{
    \stack{a_{w_{1}+1}}{\bar{b}_{w_{1}+1}}
    \cdots
    \stack{a_{w_{1}+w_{2}+w_{3}}}{\bar{b}_{w_{1}+w_{2}+w_{3}}}
    }_{w_{2}+w_{3}},
\end{equation}
where each column in the first block contributes $3$ to $\lan
c,\varphi(T)\ran$, each column in the second block contributes $2$ to
$\lan c,\varphi(T)\ran$, and the third block contributes nothing.  In
the case $T_{3}=\stack{n}{\bar{n}}$ or $T_{3}=\stack{\bar{n}}{n}$, we
have $w_{2}=0$, so we simply have $T_{1}T_{2}T_{3}=T_{1}T_{3}$, where
each column in $T_{1}$ increases $\lan c,\varphi(T)\ran$ by at least
$2$ and $T_{3}$ increases $\lan c,\varphi(T)\ran$ by $1$.  We
therefore have in the first case $\lan c,\varphi(T)\ran\geq
3(2w_{2}+w_{3})+2(w_{1}-2w_{2}-w_{3})=2w_{1}+2w_{2}+w_{3}$, and in the
second case $\lan c,\varphi(T)\ran\geq
2w_{1}+w_{3}=2w_{1}+2w_{2}+w_{3}$ , as we wished to show.

Since by Lemma \ref{lem:lb} elements in $\tbs$ have level at least
$s$, it follows that when $T\in(\tbs\cap B(s\Lt))_{\min}$, we have
$w_{1}+w_{2}\leq \lfloor s/2\rfloor$, and by $\dual$-duality that
$w_{4}+w_{5}\leq \lfloor s/2\rfloor$.  Furthermore, since
$s=w_{1}+w_{2}+w_{3}+w_{4}+w_{5}$, it follows that
$w_{1}+w_{2}=w_{4}+w_{5}=\lfloor\frac{s}{2}\rfloor$ and $w_{3}=0$ if
$s$ is even and $w_{3}=1$ if $s$ is odd.
\end{proof}

\begin{lemma} \label{lem:mixed}
    Suppose $T\in (\tbs\cap B(s\Lt))_{\min}$ has both an unbarred
    letter and a barred letter in a single column $\stack{a}{
    \bar{b}}$ other than $\stack{n}{\bar{n}}$, $\stack{\bar{n}}{ n}$,
    $\stack{n-1}{\bar{n}}$, or $\stack{n}{ \overline{n-1}}$.  Then
    $T\in\iota_{s-1}^{s}(B((s-1)\Lt)$.
\end{lemma}

\begin{proof}
    By using the reverse of Algorithm \ref{alg:iota}, it
suffices to show the following:
    \begin{enumerate}
	\item $T$ has a $1$;
    	\item $T$ has a $\bar{1}$;
    	\item after removing all $1$'s and $\bar{1}$'s, applying the
	Lecouvey $D$ relations will reduce the width of $T$.
    \end{enumerate}
    
The proof of Lemma \ref{lem:plus minus} shows that if $w_{2}\neq 0$,
or $w_{3}=1$ and $T_{3}\neq\stack{n}{\bar{n}}, \stack{\bar{n}}{n}$,
then $T$ has a $1$.  By $\dual$-duality, if $w_{4}\neq 0$, or the same
condition is placed on $w_{3}$ and $T_{3}$, then $T$ has a $\bar{1}$.
We will show that if $w_{2}+w_{3}\neq 0$, then $w_{3}+w_{4}\neq 0$,
which will prove statements (1) and (2) above.

If $w_{3}=1$ this statement is trivial, so we assume $w_{3}=0$.  We
show that the assumptions $w_{2}\neq0$ and $w_{4}=0$ lead to a
contradiction.  From the proof of Lemma \ref{lem:plus minus}, we know
that for $T$ to be minimal, every $-$ from $T_{5}$ must be bracketed.
Because of the increasing conditions on the rows and columns of $T$,
the $-$'s from the bottom row of $T_{5}$ cannot be bracketed by $+$'s
from $T_{5}$, so there must be at least $w_{5}$ $+$'s from
$T_{1}T_{2}$.  Inspection of \eqref{eq:t1t2t3} shows us that the first
block contributes no $+$'s, the second block contributes
$w_{1}-2w_{2}$ many $+$'s, and the third block contributes $2w_{2}$
many $+$'s.  We thus have $w_{1}\geq w_{5}$; but Lemma \ref{lem:plus
minus} tells us that $w_{1}+w_{2}=w_{4}+w_{5}$, contradicting our
assumption that $w_{2}\neq 0$ and $w_{4}=0$.

For the proof of statement (3), we must show that every configuration
$\stack{}{ a}\stack{c}{ b}$ in $T$ avoids the following patterns
(recall the Lecouvey $D$ sliding algorithm from section
\ref{sec:plac}): $\stack{}{x}\stack{n}{\bar{n}} $ and
$\stack{}{x}\stack{\bar{n}}{n} $ with $x\leq n-1$;
$\stack{}{n-1}\stack{\bar{n}}{\overline{n-1}} $ ;
$\stack{}{n-1}\stack{n}{\overline{n-1}} $ ; and $c\geq a$, unless
$c=a=\bar{b}$.  If $T$ has any of these patterns, the top row will not
slide over.

First, simply observe that the first four specified configurations
exclude the possibility of having a column of the form
$\stack{a}{\bar{b}}$ other than $\stack{n}{\bar{n}}$ or
$\stack{\bar{n}}{n}$.  It therefore suffices to show that the presence
of a column $\stack{d}{\bar{e}}$, $2\leq d,e\leq n-1$ implies that $T$
avoids $c\geq a$, unless $c=a=\bar{b}$.  We break our analysis of this
criterion into several special cases:

{\bf Case 1:} $a$ and $b$ are barred, $c$ is unbarred: trivial.

{\bf Case 2:} $a$ is unbarred, $b$ and $c$ are barred: This excludes 
the possibility of having $\stack{d}{\bar{e}}\in T$.

{\bf Case 3:} $a$ and $c$ are unbarred, $b$ is barred: We know the $-$
in the $c$-signature from $c$ must be bracketed; if it is by $b$, we
have $b=\bar{c}$.  As we saw in the proof of Lemma \ref{lem:plus
minus}, we must have the $-$ in the $(c-1)$-signature from $\bar{c}$
bracketed by a $c$ in the bottom row.  This forces $a\geq c$.  If the
$-$ in the $c$-signature from $c$ is bracketed by a $c+1$, it also
must be in the bottom row, forcing $a> c$.

{\bf Case 4:} $a$, $b$, $c$ all unbarred: Suppose $c\geq a$.  Let $d$
be the leftmost unbarred letter weakly to the right of $c$ which does
not have its $-$ bracketed by the letter immediately below it.  (Such
a letter exists, since we assume the occurence of
$\stack{d}{\bar{e}}\in T$, except when $d=e$; this case will be
treated below.)  This letter $d$ must be bracketed by a $d+1$ in the
bottom row, and it must be weakly to the left of $a$.  But we have
$d+1>d\geq c\geq a\geq d+1$; contradiction.

If instead we have a $\stack{d}{ \bar{d}}$ column, we must have the
$-$ from the $\bar{d}$ bracketed by a $d$ in the bottom row weakly to
the left of $a$.  In this case $d\geq c\geq a\geq d$ so $c=d$, and we
have a $\stack{d}{ }\stack{d}{\bar{d}}$ configuration, contradicting
our assumption that $T\in B(s\Lt)$.

{\bf Case 5:} $a$, $b$, $c$ all barred: Similarly to case 4, suppose
$c\geq a$ and let $d$ be the rightmost barred letter weakly to the
left of $a$ which does not have its $+$ bracketed by the letter
immediately above it.  (If none exists, we have a $\stack{\bar{d} }{
d}$ case, see below.)  It must be bracketed by a
$\overline{\bar{d}+1}$ in the top row to the right of $c$.  We then
have $d\leq a\leq c\leq\overline{\bar{d}+1}$; contradiction.

If we have a $\stack{\bar{d}}{d}$ column (note that $d$ is barred),
the $+$ from the $\bar{d}$ must be bracketed by a $d$ in the top row
to the right of $c$.  This implies that $d\geq c\geq a\geq d$, so
$a=d$ and we have a $\stack{\bar{d}}{d}\stack{}{d}$ configuration,
again contradicting our assumption that $T\in B(s\Lt)$.

\end{proof}

\begin{lemma} \label{lem:unmixed}
Let $T\in (\tbs\cap B(s\Lt))_{\min}$ such that $T$ does not contain
any column $\stack{a}{\bar{b}}$ for $1\le a,b\le n$ except possibly
$\stack{n-1}{\bar{n}}$, $\stack{n}{\bar{n}}$, or
$\stack{n}{\overline{n-1}}$.  Then $T\in \T(s)$.
\end{lemma}

\begin{proof}
    In the notation of section~\ref{sec:inj}, we have
$w_{2}=w_{4}=0$, and if $w_{3}=1$, $T_{3}=\stack{n}{\bar{n}}$ or
$T_{3}=\stack{\bar{n}}{n}$.  Lemma \ref{lem:plus minus} thus tells us
that $w_{1}=w_{5}$.

Next we show that a column $\stack{j}{i}$ must be of the form
$\stack{i-1}{i}$ for $T$ to be in $\tbs_{\min}$.  For $i=2$ we have
$j=1$ by column-strictness.  Now suppose that $\stack{j}{i}$ is the
leftmost column such that $j<i-1$.  Then $j$ contributes a $\La_j$ to
$\varphi(T)$ and hence $\lan c,\varphi(T)\ran \ge 2w_{1}+w_{3}+1=s+1$,
so that $T$ is not minimal.  By a similar argument $\lan
c,\varepsilon(T)\ran>s$ unless all columns of the form
$\stack{\bar{i}}{\bar{j}}$ must obey $j=i-1$.

A column $\stack{\overline{i}}{\overline{i-1}}$ for $i>2$ (resp.
$\stack{n}{\overline{n-1}}$) contributes a $-$ to the
$(i-2)$-signature (resp.  $(n-2)$-signature) of $T$.  This $-$ can
only be compensated by a $+$ in the $(i-2)$-signature (resp.
$(n-2)$-signature) from a column $\stack{i-1}{i}$ (resp.
$\stack{n-1}{\overline{n}}$).  Hence for $T$ to be minimal the number
of columns of the form $\stack{i-1}{i}$ (resp.
$\stack{n-1}{\overline{n}}$) needs to be the same as the number of
columns of the form $\stack{\overline{i}}{\overline{i-1}}$ (resp.
$\stack{n}{\overline{n-1}}$).  This proves that $T\in \T(s)$.
\end{proof}

\subsection{Uniqueness} \label{sec:unique}
Theorem~\ref{thm:b2s} follows as a corollary from the next
Proposition.

\begin{proposition}\label{prop:uniqueness}
$\tbs$ is the only affine finite-dimensional crystal satisfying the
properties of Conjecture~\ref{conj:Brs}.
\end{proposition}

For the proof of Proposition~\ref{prop:uniqueness} we must show that
our definition of $\sigma$ is the only crystal automorphism satisfying
the properties of Conjecture~\ref{conj:Brs}.  Recall from the
beginning of section~\ref{sec:KR:sigma} the relationship between $\sigma$
and $\hsig$.  Let $T\in\tbs$.  We know that
\begin{equation}
    \label{eq:sigma prop}
    \wt(T)=\sum_{i=0}^{n}k_{i}\Lambda_{i}\Leftrightarrow 
    \wt(\sigma(T))=k_{1}\Lz+k_{0}\Lo+\sum_{i=2}^{n}k_{i}\Lambda_{i}.
\end{equation}

Please note that in this section we often use the phrase ``the tableau
$b$ is in the branching component vertex $v$'' to mean $b\in
B(v)$.

\begin{lemma} \label{lem:middle}
    Let $v\in\Sup(k\Lt)$ be a branching component vertex in stratum
    $0$ associated with a rectangular partition, and let $\ell$ be
    minimal such that $v\in\iotil_{\ell}^{k}(\Sup(\ell\Lt))$.  Then
    the hypothesis that $\tbs$ is perfect of level $s$ can only be
    satisfied if $\hsig(v)$ is the vertex associated with the same
    shape as $v$ in stratum $0$ in $\Sup((s+\ell-k)\Lt)$.
\end{lemma}

\begin{proof} 
    We have already shown that $\hsig(v)$ has the same shape as $v$
    and is in stratum $0$, so it only remains to show that
    $\hsig(v)\in\Sup((s+\ell-k)\Lt)$.
    
    First, observe that $v$ must contain a minimal tableau as
    constructed in section \ref{sec:sur}, according to the following 
    table.
    
    \vspace{.1in}    
    \begin{center}
    $
    \begin{array}{c|c}
	\textrm{shape associated with }v & \textrm{weight of tableau
	in }v \\
	\hline
	(2m) & m\Lambda_2+(k-2m)\Lo+(s-k)\Lz \\
	(m,m)  & m\Lambda_{3}+(k-2m)\Lo+(s-k)\Lz \\
    \end{array}
    $
    \end{center}
    \vspace{.1in}

    Recall from \ref{sec:KR:sigma} that
    $\mathrm{cw}(b)=\mathrm{wt}(b)-(\varphi_{0}(b)-\vare_{0}(b))\Lambda_0$.
    Let $T$ be the tableau constructed by this prescription, so that
    $\lan c,\cw(T)\ran=k$.  The criterion that $\lan
    c,\varphi(T)\ran\geq s$ forces us to have
    $\varphi_{0}(T)=\varphi_{1}(\sigma(T))\geq s-k.$ We denote
    $T^{(i)}=\iota_{k}^{i}(T)$, and thus have
    $\varphi_{1}(T^{(i)})=i-\ell$.
    
    We show inductively that $\sigma(T^{(i)})=T^{(s+\ell-i)}$ for
    $\ell\leq i\leq s$.  As a base case, we see that $\lan
    c,\cw(T^{(\ell)})\ran=\ell$, so we must have
    $\varphi_{1}(\sigma(T^{(\ell)}))\geq s-\ell$.  The only $T^{(i)}$ for
    which this inequality holds is $T^{(s)}$, where we have
    $\varphi_{1}(T^{(s)})=s-\ell$.
    
    For the induction step, assume that $\sigma$ sends
    $T^{(\ell)},T^{(\ell+1)},\ldots,T^{(k-1)}$ to\\
    $T^{(s)},\ldots, T^{(s+\ell-k+1)}$, respectively.  By the above
    inequality this implies $\varphi_{1}(\sigma(T^{(k)}))\geq s-k$,
    which specifies that $\sigma(T^{(k)})=T^{(s+\ell-k)}$.
\end{proof}

\begin{definition}\label{def:ups check}
    Recall the map $\Upsilon_{s'}^{s}$ from
    Definition~\ref{def:upsilon}.  We define
    $\check{\Upsilon}_{s'}^{s}:\Sup(\tilde{B}^{2,s'})\hookrightarrow
    \Sup(\tilde{B}^{2,s})$ by $\check{\Upsilon}_{s'}^{s}(v)=v'$ if for
    some $T\in B(v)$, we have $\Upsilon_{s'}^{s}(T)\in B(v')$.
\end{definition}

\begin{lemma} \label{lem:down}
    Let $v\in\Sup(k\Lt)$ be a branching component vertex in stratum
    $p$ with $1\leq p\leq s-1$ associated to the shape
    $(\lambda_{1},\lambda_{2})$.  Suppose that for the branching
    component vertex $w\in\Sup(k\Lt)$ in stratum $p-1$ with shape
    $(\lambda_{1}-1,\lambda_{2})$, $\tbs$ has the correct energy
    function and is perfect only if $\hsig(w)$ is as described in
    section \ref{sec:hsig}.  Then $\tbs$ has the correct energy
    function only if $\hsig(v)$ is as described in section
    \ref{sec:hsig}.
\end{lemma}

\begin{proof}
    First, recall that the partitions associated to vertices in
    stratum $p$ in $\Sup(\tbs)$ are produced by adding or removing one
    box from the partitions associated to vertices of stratum $p+1$.
    Since the vertex in stratum $s$ is associated with a rectangle of
    shape $(s)$, the only stratum for which we can have a two-row
    rectangle is stratum $0$.  It follows that removing a box from the
    first row of $v$ results in a partition apppearing in stratum
    $p-1$, so there is in fact a vertex $w$ as described in the
    statement of the lemma.
        
    Let $\ell$ be minimal such that
    $v\in\iotil_{\ell}^{k}(\Sup(\ell\Lt))$.  We may assume
    $v\notin\Sup(\ell\Lt)$, and let $\hsig(v)$ be determined by the
    involutive property of $\hsig$ in the case $v\in\Sup(\ell\Lt)$.
    Specifically, we will show that the vertex $v'\in\Sup(s\Lt)$ with
    the same shape and stratum as $v$ has the property that
    $\hsig(v')$ is the complementary vertex of $v$, and therefore
    $\hsig(v)$ is the complementary vertex of $v'$.  (Recall the
    notion of complementary vertex from page \pageref{comp_vert}.)
    
    The top row of $w$ is one box shorter than the top row of $v$, so
    that $w\in\iotil_{\ell-1}^{k}(\Sup((\ell-1)\Lt))$, and $\ell-1$ is
    minimal with this property.  By our hypothesis, $\hsig(w)$ is the
    vertex with shape $(\lambda_{1}-1,\lambda_{2})$ of stratum $-p+1$
    in $\Sup((s+(\ell-1)-k)\Lt)$.
    
    We now use induction on $s$.  Suppose that the only choice of
    $\hsig$ for which $\tilde{B}^{2,s-1}$ is perfect and has an energy
    function is that of section~\ref{sec:hsig}.  Part 3 of
    Conjecture~\ref{conj:Brs} states that in $\tbs$, the energy on the
    component $B((s-k)\Lt)$ is $-k$, and so the difference in energy
    between $B((s-k)\Lt)$ and $B((s-j)\Lt)$ is $j-k$.  In order for
    this to be true for all $1\leq k,j\leq s-1$, the action of
    $\ftil{0}$ and $\etil{0}$ on $\tbs$ must agree with the action on
    $\tilde{B}^{2,s-1}$.  More precisely, if $v$ and $w$ are in
    different classical components of $\tilde{B}^{2,s-1}$ and
    $\ftil{0}(v)=w$ in $\tilde{B}^{2,s-1}$, then
    $\ftil{0}(\Upsilon_{s-1}^{s}(v))=\Upsilon_{s-1}^{s}(w)$; this
    statement extends naturally to $\Sup(\tilde{B}^{2,s-1})$ and
    $\Sup(\tbs)$.
    
    Let $v^{\dagger}$ denote the vertex with shape
    $(\lambda_{1},\lambda_{2})$ in stratum $-p$ in
    $\Sup((s+\ell-k)\Lt)$.  Since we assumed $k\neq\ell$, we know that
    $v^{\dagger}\notin\Sup(s\Lt)$, and therefore $v^{\dagger}$ has a
    preimage under $\check{\Upsilon}_{s-1}^{s}$.  From our
    construction of $\hsig$ we know that in $\Sup(\tilde{B}^{2,s-1})$,
    $(\check{\Upsilon}_{s-1}^{s})^{-1}(v^{\dagger})$ has a $0$ arrow
    to $(\check{\Upsilon}_{s-1}^{s})^{-1}(\hsig(w))$.  Our induction
    argument tells us that in $\Sup(\tbs)$, $v^{\dagger}$ has a $0$
    arrow to $\hsig(w)$.  Since $v$ has a $1$ arrow to $w$ and we must
    have $\ftil{0}=\sigma\ftil{1}\sigma$, we conclude that in fact
    $\hsig(v)=v^{\dagger}$.
\end{proof}

\begin{lemma} \label{lem:up}
    Let $v\in\Sup(k\Lt)$ be a branching component vertex in stratum
    $0$ associated to the non-rectangular shape
    $(\lambda_{1},\lambda_{2})$, and suppose that for the branching
    component vertex $w\in\Sup(k\Lt)$ in stratum $1$ with shape
    $(\lambda_{1},\lambda_{2}+1)$, $\tbs$ has the correct energy
    function and is perfect only if $\hsig(w)$ is as described in
    section \ref{sec:hsig}.  Then $\tbs$ has the correct energy
    function only if $\hsig(v)$ is as described in section
    \ref{sec:hsig}.
\end{lemma}

\begin{proof}
    (This proof is very similar to the proof of Lemma~\ref{lem:down}.)
    
    Let $\ell$ be minimal such that
    $v\in\iotil_{\ell}^{k}(\Sup(\ell\Lt))$, assuming $k\neq\ell$.
    Note that $\ell$ is also minimal for
    $w\in\iotil_{\ell}^{k}(\Sup(\ell\Lt))$, since the shapes for $v$
    and $w$ have the same number of boxes in the first row.  By our
    hypothesis, $\hsig(w)$ is the vertex with shape
    $(\lambda_{1},\lambda_{2}+1)$ in stratum $1$ in
    $\Sup((s+\ell-k)\Lt)$.  Let $v^{\dagger}$ be the vertex with shape
    $(\lambda_{1},\lambda_{2})$ in stratum $0$ in
    $\Sup((s+\ell-k)\Lt)$.  From our construction of $\hsig$, we know
    that in $\Sup(\tilde{B}^{2,s-1})$,
    $(\check{\Upsilon}_{s-1}^{s})^{-1}(w)$ has a $0$ arrow to
    $(\check{\Upsilon}_{s-1}^{s})^{-1}(\hsig(v^{\dagger}))$.  It
    follows from our induction argument that in $\Sup(\tbs)$, $w$ has
    a $0$ arrow to $\hsig(v^{\dagger})$.  Since $w$ has a $1$ arrow to
    $v$ and we must have $\ftil{0}=\sigma\ftil{1}\sigma$, we conclude
    that in fact $\hsig(v)=v^{\dagger}$.
\end{proof}

\begin{corollary} Corollary~\ref{prop:flip} and
Lemmas~\ref{lem:middle}, \ref{lem:down} and \ref{lem:up} determine
$\hsig$ on $\Sup(\tbs)$ uniquely.
\end{corollary}

\begin{proof}
For any vertex $v$ associated with shape $(\lambda_{1},\lambda_{2})$
in stratum $p\geq 0$, $\hsig(v)$ is fixed by the image under $\hsig$
of a vertex with shape $(\lambda_{1}-s+p,\lambda_{2})$ and in stratum
$0$ by Lemma~\ref{lem:down}.  If $\lambda_{1}-s+p\neq\lambda_{2}$,
Lemmas~\ref{lem:down} and~\ref{lem:up} may be used together to reduce
determining $\hsig(v)$ to determining the action of $\hsig$ on a
rectangular vertex in stratum $0$, which is given by
Lemma~\ref{lem:middle}.
\end{proof}
        
    %
    %

    \newchap{Local Properties of Crystals for Modules over Doubly 
    Laced Algebras}
    \label{sec:local}
    
        \section{Local Characterization of Simply Laced Crystals}
        \label{sec:local:simply}
        In \cite{Stem}, a local charactarization of crystals for highest
weight integrable representations of simply-laced algebras is given.
Furthermore, Proposition 2.4.4 of \cite{KKMMNN:1992a} states that a
crystal with a unique maximal vertex comes from a representation if
and only if it decomposes as a disjoint union of crystals of
representations relative to the rank 2 subalgebras corresponding to
each pair of edge colors.  It therefore suffices to address the
problem of locally characterizing crystal graphs for rank 2 algebras.
The results of \cite{Stem} apply to the algebras $A_{1}\times A_{1}$
and $A_{2}$; the obvious next case to consider is $C_{2}\simeq B_{2}$.
In the sequel, a ``doubly laced algebra'' is a Kac-Moody algebra all
of whose regular rank $2$ subalgebras are $A_{1}\times A_{1}$,
$A_{2}$, or $C_{2}$.

We now recall the results of \cite{Stem} on the local structure of
crystal graphs.  While the main theorem of \cite{Stem} pertains only
to simply laced algebras, several statements that apply to crystals
for an arbitrary Kac-Moody algebra are also proved there.

We begin by choosing an index set $I$ and a Cartan matrix $A = 
[a_{ij}]_{i,j\in I}$ for a Kac-Moody algebra, and consider a directed 
graph $X$ whose edges are colored by the elements of $I$.  We assume that
\begin{itemize}
    \item[(P1)] All monochromatic directed paths in $X$ have finite length

    \item[(P2)] For any vertex $v$ of $X$ and every $i\in I$, there is at 
    most one $i$ edge into $v$ and one $i$ edge out of $v$.
\end{itemize}
These assumptions (the first two axioms of \cite{Stem}) allow us to
define the $i$-string through $v$ to be the maximal path in $X$ of the
form
$$
f_{i}^{-d}v \to \cdots \to f_{i}^{-1}v \to v \to f_{i}v \to \cdots 
\to f_{i}^{r}v.
$$
These familiar statistics $d = \vare_{i}(v)$ and $r = \varphi_{i}(v)$
are here called the $i$-depth and $i$-rise of $v$, respectively.

\begin{remark}
    Our conventions differ from those of \cite{Stem} by the sign of
    the $i$-depth and the greek letters denoting these quantities. 
    This is more consistent with the existing literature on crystals. 
\end{remark}

Stembridge introduces the following local statistics for a crystal
vertex $v$  such that $\vare_{i}(v)>0$ and $\vare_{j}(v)>0$:
$$
\Delta_{i}\vare_{j}(v) = \vare_{j}(v) - \vare_{j}(e_{i}v)
\qquad
\mathrm{and}
\qquad
\Delta_{i}\varphi_{j}(v) = \varphi_{j}(v) - \varphi_{j}(e_{i}v),
$$
and complementarily,
$$
\nabla_{i}\vare_{j}(v) = \vare_{j}(v) - \vare_{j}(f_{i}v)
\qquad
\mathrm{and}
\qquad
\nabla_{i}\varphi_{j}(v) = \varphi_{j}(v) - \varphi_{j}(f_{i}v).
$$
He shows that for a crystal of a module over any Kac-Moody algebra, we
have the following
\begin{enumerate}
    \item[(P3)] $\Delta_{i}\vare_{j}(v) + \Delta_{i}\varphi_{j}(v) =
    a_{ij}$

    \item[(P4)] $\Delta_{i}\varphi_{j}(v) \leq 0$ and
    $\Delta_{i}\vare_{j}(v) \leq 0$

    \item[(P5)] If $\Delta_{i}\vare_{j}(v) = 0$, then $e_{i}e_{j}v =
    e_{j}e_{i}v$ and $\nabla_{j}\varphi_{i}(e_{i}e_{j}v) = 0$

    \item[(P6)] If $\Delta_{i}\vare_{j}(v) = \Delta_{j}\vare_{i}(v) =
    -1$, then $e_{i}e_{j}^{2}e_{i}v = e_{j}e_{i}^{2}e_{j}v$, and
    $\nabla_{i}\varphi_{j}(e_{i}e_{j}^{2}e_{i}v) =
    \nabla_{i}\varphi_{j}(e_{i}e_{j}^{2}e_{i}v) = -1$.

    \item[(P5')] If $\nabla_{i}\vare_{j}(v) = 0$, then $e_{i}e_{j}v =
    e_{j}e_{i}v$ and $\Delta_{j}\varphi_{i}(e_{i}e_{j}v) = 0$

    \item[(P6')] If $\nabla_{i}\vare_{j}(v) = \nabla_{j}\vare_{i}(v) =
    -1$, then $e_{i}e_{j}^{2}e_{i}v = e_{j}e_{i}^{2}e_{j}v$, and
    $\Delta_{i}\varphi_{j}(e_{i}e_{j}^{2}e_{i}v) =
    \Delta_{i}\varphi_{j}(e_{i}e_{j}^{2}e_{i}v) = -1$.
\end{enumerate}

Stembridge proved that these axioms hold for all crystals of highest
weight modules of Kac-Moody algebras, and that they precisely
characterize those over simply laced algebras.  In the work that
follows, we show that when generalizing to doubly laced algebras,
there are only two other possible local behaviors,
$$
e_{i}e_{j}e_{i}e_{j}e_{i}v = e_{j}e_{i}^{3}e_{j}v
\qquad
\textrm{and}
\qquad
e_{i}e_{j}^{2}e_{i}^{3}e_{j}v = e_{j}e_{i}^{3}e_{j}^{2}e_{i}v.
$$

Note that Theorem \ref{proposition:int} does not provide a local
characterization of crystals coming from representations of doubly
laced algebras.  In order to have such a characterization, it would be
necessary to provide additional axioms to those stated above and to
show that any graph satisfying those axioms is in fact a crystal.
Here, we only show the other half of the characterization; we are
assured that any graph with a relation not explicitly described in the
above theorem is in fact not a crystal over one of these algebras.

Another application of the results of \cite{Stem} appears in
\cite{DKK}, which gives a novel construction of $A_{2}$-crystals based
on half-grids.  Results such as those in \cite{Stem} and \cite{DKK},
as well as those discussed below, impart an increased understanding of
the graph/poset theoretic aspects of crystal graphs, revealing a new
perspective on the representation theory of Lie algebras.

	\section{Realization of Crystals for Modules over Type $C_{2}$
	Algebras}
        \label{sec:local:B}

Recall that a crystal is a colored directed graph in which we
interpret an $i$-colored edge from the vertex $x$ to the vertex $y$ to
mean that $f_{i}x = y$ and $e_{i}y = x$, where $e_{i}$ and $f_{i}$ are
Kashiwara crystal operators.  The vector representation of $C_{2}$ has
the following crystal: \vspace{.2in}
\begin{center}
    \setlength{\unitlength}{1in}
    \begin{picture}(.5,0)
	\put(-1.05,0){\framebox{$1$}}
	\put(-.25,0){\framebox{$2$}}
	\put(.55,0){\framebox{$\bar{2}$}}
	\put(1.3,0){\framebox{$\bar{1}$}}
	\put(-.8,0.05){\vector(1,0){.4}}
	\put(-.6,.16){\footnotesize{1}}
	\put(.05,0.05){\vector(1,0){.4}}
	\put(.27,.16){\footnotesize{2}}
	\put(.8,0.05){\vector(1,0){.4}}
	\put(.95,.16){\footnotesize{1}}
    \end{picture}
    
\end{center}
We realize crystals using tableaux filled with the letters
$\{1,2,\bar{2},\bar{1}\}$ with the total ordering $1 < 2 < \bar{2} <
\bar{1}$.
The following definition is
adapted from that in \cite{KN}.

\begin{definition}
    \label{def:tableau}
	
    A $C_{2}$ Young diagram is a partition with no more than two
    parts; we draw it as a left-justified two-row arrangement of boxes
    such that the second row is no longer than the first.
    
    A $C_{2}$ tableau is a filling of a $C_{2}$ Young diagram by the
    letters of the above alphabet with the following properties:
    \begin{enumerate}
	\item each row is weakly increasing by the ordering in the
	vector representation;

	\item each column is strictly increasing by the ordering in
	the vector representation;

	\item no column may contain $1$ and $\bar{1}$ simultaneously;
	
	\item the configuration 	     
	$\begin{array}{|c|c|}
		\hline
		2 & *  \\
		\hline
		* & \raisebox{-1pt}{$\bar{2}$} \\
		\hline
		 
	\end{array}
	$ does not appear in $T$.

    \end{enumerate}
\end{definition}

\begin{example}
    \label{ex:tableau}
    The following is an example of a $B_{2}$ tableau.
	$$
	\begin{array}{|c|c|c|c|c|c|c|c}
	    \cline{1-7}
	    \raisebox{-2pt}{$1$} & \raisebox{-2pt}{$1$} & 
	    \raisebox{-2pt}{$2$} & \raisebox{-2pt}{$2$} & 
	    \raisebox{-2pt}{$\bar{2}$} & \raisebox{-2pt}{$\bar{2}$} &
	    \raisebox{-2pt}{$\bar{1}$} & \\
	    \cline{1-7}        
	    \raisebox{-2pt}{$2$} & \raisebox{-2pt}{$2$} &
	    \raisebox{-2pt}{$\bar{2}$} & \raisebox{-2pt}{$\bar{1}$} &
	    \multicolumn{4}{c}{} \\
	    \cline{1-4} 
	\end{array}
	$$
\end{example}

\begin{definition}
    \label{def:colword}
    
    Let $T$ be a type $C_{2}$ tableau.  The column word $W$ of $T$ is
    the word on the alphabet $\{ 1, 2, \bar{2}, \bar{1}\}$ consisting
    of $cd$ for each column $
    \begin{array}{|c|} 
	\hline
        d  \\
        \hline
        c  \\
        \hline
         
    \end{array}
    $ in $T$, reading left to right, then followed by each entry
    appearing in a one-row column in $T$, again reading left to right.
    
\end{definition}
(Compare to Definition \ref{def:col_word})

\begin{example}
    \label{ex:colword}
    The column word of the tableau in Example \ref{ex:tableau} is \\
    $2121\bar{2}2\bar{1}2\bar{2}\bar{2}\bar{1}$.
\end{example}

We now present a definition of the $1$-signature and $2$-signature of
the column word of a type $C_{2}$ tableau, which is easily seen to be
equivalent to the conventional definitions (e.g. \cite{HK}).  Our
definition differs by using the extra symbol $*$ to keep track of
vacant spaces in the signatures.

\begin{definition}
    \label{def:sigs}
    
    Let $a$ be in the alphabet $\{ 1, 2, \bar{2}, \bar{1}\}$.  Then
    the $1$-signature of $a$ is
    \begin{itemize}
	\item $-$, if $a$ is $\bar{2}$ or $1$;
	
	\item $+$, if $a$ is $\bar{1}$ or $2$;
    
    \end{itemize}
    The $2$-signature of $a$ is
        \begin{itemize}

	\item $-$, if $a$ is $2$;

	\item $+$, if $a$ is $\bar{2}$;

        \item $*$, if $a$ is $1$ or $\bar{1}$.
    \end{itemize}
    
    Let $W$ be the column word of a type $C_{2}$ tableau $T$.  Then
    for $i\in\{1,2\}$ the $i$-signature of $T$ is the word on the
    alphabet $\{ +, -, *\}$ that results from concatenating the
    $i$-signatures of the entries of $W$.

\end{definition}

\begin{example}
    \label{ex:sigs}
    The $1$- and $2$-signatures of the tableau in example 
    \ref{ex:tableau} are 
    $$
    +-+--+++--+
    $$
    and 
    $$
    -*-*+-*-++*,
    $$
    respectively.
\end{example}

\begin{definition}
    \label{def:redsig}
    
    Let $S = s_{1}s_{2}\cdots s_{\ell}$ be a signature in the sense of
    Definition \ref{def:sigs}.  The reduced form of $S$ is the word on
    the alphabet $\{+,-,*\}$ that results from iteratively replacing
    every occurance of $+\underbrace{*\cdots*}_{k}-$ in $S$ with
    $\underbrace{*\cdots*}_{k+2}$ until there are no occurances of
    $+\underbrace{*\cdots*}_{k}-$ in $S$.
\end{definition}

\begin{example}
    \label{ex:redsigs}
    The reduced forms of the signatures in Example \ref{ex:sigs} are
    $$
    ****-+****+
    $$
    and 
    $$
    -*-****-++*,
    $$    
    respectively.
\end{example}

The result of applying the Kashiwara operator $e_{i}$ to a tableau $T$
breaks into several cases.  If there are no $+$'s in the reduced form
of the $i$-signature of $T$, we say that $e_{i}T = 0$, where $0$ is a
formal symbol.  Otherwise, let $a$ be the entry corresponding to the
leftmost $+$ in the reduced form of the $i$-signature of $T$.  Then
$e_{i}T$ is the tableau that results from changing $a$ to $e_{i}a$ in
$T$, as defined by the vector representation.

Similarly for $f_{i}$, if there are no $-$'s in the reduced form of
the $i$-signature of $T$, we say that $f_{i}T = 0$.  Otherwise, let
$a$ be the entry corresponding to the rightmost $-$ in the reduced
form of the $i$-signature of $T$.  Then $f_{i}T$ is the tableau that
results from changing $a$ to $f_{i}a$ in $T$.

\begin{example}
    \label{ex:kash}
    Let 
    $$
    T=
    \begin{array}{|c|c|c|c|c|c|c|c}
	\cline{1-7}
	\raisebox{-2pt}{$1$} & \raisebox{-2pt}{$1$} & 
	\raisebox{-2pt}{$2$} & \raisebox{-2pt}{$2$} & 
	\raisebox{-2pt}{$\bar{2}$} & \raisebox{-2pt}{$\bar{2}$} &
	\raisebox{-2pt}{$\bar{1}$} & \\
	\cline{1-7}        
	\raisebox{-2pt}{$2$} & \raisebox{-2pt}{$2$} &
	\raisebox{-2pt}{$\bar{2}$} & \raisebox{-2pt}{$\bar{1}$} &
	\multicolumn{4}{c}{} \\
	\cline{1-4} 
    \end{array}
    $$ 
    be the tableau of Example \ref{ex:tableau}.  We then have
    $$
    f_{1}T=
    \begin{array}{|c|c|c|c|c|c|c|c}
	\cline{1-7}
	\raisebox{-2pt}{$1$} & \raisebox{-2pt}{$1$} & 
	\raisebox{-2pt}{$2$} & \raisebox{-2pt}{$2$} & 
	\raisebox{-2pt}{$\bar{2}$} & \raisebox{-2pt}{$\bar{2}$} &
	\raisebox{-2pt}{$\bar{1}$} & \\
	\cline{1-7}        
	\raisebox{-2pt}{$2$} & \raisebox{-2pt}{$2$} &
	\raisebox{-2pt}{$\bar{1}$} & \raisebox{-2pt}{$\bar{1}$} &
	\multicolumn{4}{c}{} \\
	\cline{1-4} 
    \end{array},
    $$
    $$
    f_{2}T=
    \begin{array}{|c|c|c|c|c|c|c|c}
	\cline{1-7}
	\raisebox{-2pt}{$1$} & \raisebox{-2pt}{$1$} & 
	\raisebox{-2pt}{$2$} & \raisebox{-2pt}{$\bar{2}$} & 
	\raisebox{-2pt}{$\bar{2}$} & \raisebox{-2pt}{$\bar{2}$} &
	\raisebox{-2pt}{$\bar{1}$} & \\
	\cline{1-7}        
	\raisebox{-2pt}{$2$} & \raisebox{-2pt}{${2}$} &
	\raisebox{-2pt}{$\bar{2}$} & \raisebox{-2pt}{$\bar{1}$} &
	\multicolumn{4}{c}{} \\
	\cline{1-4} 
    \end{array},
    $$
    $$
    e_{1}T=
    \begin{array}{|c|c|c|c|c|c|c|c}
	\cline{1-7}
	\raisebox{-2pt}{$1$} & \raisebox{-2pt}{$1$} & 
	\raisebox{-2pt}{$1$} & \raisebox{-2pt}{$2$} & 
	\raisebox{-2pt}{$\bar{2}$} & \raisebox{-2pt}{$\bar{2}$} &
	\raisebox{-2pt}{$\bar{2}$} & \\
	\cline{1-7}        
	\raisebox{-2pt}{$2$} & \raisebox{-2pt}{$2$} &
	\raisebox{-2pt}{$\bar{2}$} & \raisebox{-2pt}{$\bar{1}$} &
	\multicolumn{4}{c}{} \\
	\cline{1-4} 
    \end{array},
    $$
    and
    $$
    e_{2}T=
    \begin{array}{|c|c|c|c|c|c|c|c}
	\cline{1-7}
	\raisebox{-2pt}{$1$} & \raisebox{-2pt}{$1$} & 
	\raisebox{-2pt}{$2$} & \raisebox{-2pt}{$2$} & 
	\raisebox{-2pt}{${2}$} & \raisebox{-2pt}{$\bar{2}$} &
	\raisebox{-2pt}{$\bar{1}$} & \\
	\cline{1-7}        
	\raisebox{-2pt}{$2$} & \raisebox{-2pt}{$2$} &
	\raisebox{-2pt}{$\bar{2}$} & \raisebox{-2pt}{$\bar{1}$} &
	\multicolumn{4}{c}{} \\
	\cline{1-4} 
    \end{array}.
    $$
\end{example}

        \section{Analysis of Generic $C_{2}$ tableaux}
        \label{sec:local:genB}
        A generic $C_{2}$ tableau is of the form
\begin{equation}
    \label{eq:gentab}
    \begin{array}{|c|c|c|c|c|c|c|c|c|c|c|c|c|c|c|c|}
	\hline
	\raisebox{-1pt}{$1$} & \cdots & \raisebox{-1pt}{$1$} &
	\raisebox{-1pt}{$1$} & \cdots & \raisebox{-1pt}{$1$}  &
	\raisebox{-1pt}{$2$} & \raisebox{-1pt}{$2$} & \cdots &
	\raisebox{-1pt}{$2$} & 
	\raisebox{-1pt}{$\bar{2}$} & \cdots &
	\raisebox{-1pt}{$\bar{2}$} & \raisebox{-1pt}{$\bar{1}$} & \cdots &
	\raisebox{-1pt}{$\bar{1}$}  \\
	\cline{1-16}
	\raisebox{-1pt}{$2$} & \cdots & \raisebox{-1pt}{$2$} &
	\raisebox{-1pt}{$\bar{2}$} & \cdots &
	\raisebox{-1pt}{$\bar{2}$} &
	\raisebox{-1pt}{$\bar{2}$} & \raisebox{-1pt}{$\bar{1}$} & \cdots &
	\raisebox{-1pt}{$\bar{1}$} & \raisebox{-1pt}{$\bar{1}$} & \cdots &
	\raisebox{-1pt}{$\bar{1}$} 
	& \multicolumn{3}{}{c}  \\
	\cline{1-13}
    \end{array}
\end{equation}
where 
\begin{itemize}
    \item any column may be omitted;

    \item  any of the columns other than
    $\begin{array}{|c|} \hline
	2 \\
	\hline
	\raisebox{-1pt}{$\bar{2}$}  \\
	\hline
	 
    \end{array}
    $ may be repeated an arbitrary number of times;

    \item the bottom row
    may be truncated at any point.

\end{itemize}

We are interested in how the Kashiwara operators $e_{1}$ and $e_{2}$
act on this tableaux, so we must determine where the left-most $+$
appears in the reduced form of the signatures of the tableau.  The
relevant $+$'s in the signatures of a generic tableau naturally fall
into two groups as described by definition \ref{def:blocks}.

\begin{definition}
    \label{def:blocks}
    Let $T$ be a $C_{2}$ tableau.  
    \begin{itemize}
	\item We define the left block of $+$'s in the $1$-signature
	of $T$ to be those $+$'s from $2$'s in the top row and
	$\bar{1}$'s in the bottom row.  If no such entries appear in
	$T$, we say that the left block of $+$'s in the $1$-signature
	of $T$ has size $0$ and its left edge is located on the
	immediate left of symbol coming from the leftmost $\bar{2}$ or
	$\bar{1}$ in the top row of $T$.  If there is furthermore no
	such entry, its left edge is located at the right end of the
	$1$-signature of $T$.
    
	\item We define the right block of $+$'s in the $1$-signature
	of $T$ to be those $+$'s from $\bar{1}$ in the top row of $T$.
	If no such entry appears in $T$, we say that the right block
	of $+$'s in the $1$-signature of $T$ has size $0$ and its left
	edge is located at the right end of the $1$-signature of $T$.
    
	\item We define the left block of $+$'s in the $2$-signature
	of $T$ to be those $+$'s from $\bar{2}$ in the bottom row.  If
	such an entry does not appear in $T$, we say that the left
	block of $+$'s in the $2$-signature of $T$ has size $0$ and
	its left edge is located on the immediate left of the $*$ in
	the $2$-signature coming from the leftmost $\bar{1}$ in the
	bottom row of $T$.  If $T$ has no $\bar{1}$'s in the bottom
	row, we say that its left edge is located at the right end of
	the $2$-signature of $T$.
    
	\item We define the right block of $+$'s in the $2$-signature
	of $T$ to be those $+$'s from $\bar{2}$ in the top row of $T$.
	If such an entry does not appear in $T$, we say that the right
	block of $+$'s in the $2$-signature of $T$ has size $0$ and
	its left edge is located on the immediate left of the $*$ in
	the $2$-signature coming from the leftmost $\bar{1}$ in the
	top row of $T$.  If there are furthermore no $\bar{1}$'s in
	the top row of $T$, we say that its left edge is located at
	the right end of the $2$-signature of $T$.
	
    \end{itemize}
    
    In the above cases when a block of $+$'s has positive size we say
    that its left edge is on the immediate left of its leftmost $+$.
    
\end{definition}

Motivated by this definition, we define the following statistics on a
$C_{2}$ tableaux $T$.
\begin{itemize}
    \item $A(T)$ is the number of $\bar{2}$'s in the top row of $T$,

    \item $B(T)$ is the number of $2$'s in the top row of $T$ plus the
    number of $\bar{1}$'s in the bottom row of $T$,

    \item $C(T)$ is the number of $2$'s in the top row of $T$,

    \item $D(T)$ is the number of $\bar{2}$'s in the bottom row of
    $T$.
\end{itemize}

\begin{example}
    Let 
    $$
    T = 
    \begin{array}{|c|c|c|c|c|c|c}
        
	\cline{1-6}
        \raisebox{-1pt}{1} & \raisebox{-1pt}{1} & \raisebox{-1pt}{2} & \raisebox{-1pt}{$\bar{2}$} & \raisebox{-1pt}{$\bar{2}$} & \raisebox{-1pt}{$\bar{1}$} &   \\
        \cline{1-6}
        \raisebox{-1pt}{2} & \raisebox{-1pt}{$\bar{2}$} & \raisebox{-1pt}{$\bar{2}$} & \raisebox{-1pt}{$\bar{1}$} & \multicolumn{3}{c}{}    \\
        
        \cline{1-4} 
    \end{array}
    .
    $$
    Then $A(T) = 2$, $B(T) = 2$, $C(T)  = 1$, and $D(T) = 2$.
\end{example}

\begin{claim} \label{claim:lr}
    $\begin{array}{c}
        
          \\

    \end{array}
    $
    \begin{itemize}
	\item If $e_{1}$ acts on a tableau $T$ with $A(T) < B(T)$, the
	entry on which $e_{1}$ acts corresponds to a symbol in the
	left block of $+$'s in the $1$-signature of $T$;
    
	\item If $e_{1}$ acts on a tableau $T$ with $A(T) \geq B(T)$,
	the entry on which $e_{1}$ acts corresponds to a symbol in
	the right block of $+$'s in the $1$-signature of $T$;
    
	\item If $e_{2}$ acts on a tableau $T$ with $C(T) < D(T)$, the
	entry on which $e_{2}$ acts corresponds to a symbol in the
	left block of $+$'s in the $2$-signature of $T$;
    
	\item If $e_{2}$ acts on a tableau $T$ with $C(T) \geq D(T)$,
	the entry on which $e_{2}$ acts corresponds to a symbol in
	the right block of $+$'s in the $2$-signature of $T$.
    \end{itemize}
    
\end{claim}

\begin{example}
    \label{ex:1lr}
    Let 
    $$
    T_{1}=
        \begin{array}{|c|c|c|c|c|c|c}
	\cline{1-6}
	\raisebox{-2pt}{$1$} & \raisebox{-2pt}{$2$} & 
	\raisebox{-2pt}{$2$} & \raisebox{-2pt}{$\bar{2}$} & 
	\raisebox{-2pt}{$\bar{2}$} & \raisebox{-2pt}{$\bar{1}$} & \\
	\cline{1-6}        
	\raisebox{-2pt}{$2$} & \raisebox{-2pt}{$\bar{2}$} &
	\raisebox{-2pt}{$\bar{1}$} & 
	\multicolumn{4}{c}{} \\
	\cline{1-3} 
    \end{array}
    $$
    and
    $$
    T_{2}=
        \begin{array}{|c|c|c|c|c|c|c}
	\cline{1-6}
	\raisebox{-2pt}{$1$} & \raisebox{-2pt}{$1$} & 
	\raisebox{-2pt}{$2$} & \raisebox{-2pt}{$\bar{2}$} & 
	\raisebox{-2pt}{$\bar{2}$} & \raisebox{-2pt}{$\bar{1}$} & \\
	\cline{1-6}        
	\raisebox{-2pt}{$2$} & \raisebox{-2pt}{$2$} &
	\raisebox{-2pt}{$\bar{1}$} & 
	\multicolumn{4}{c}{} \\
	\cline{1-3} 
    \end{array}.
    $$
    In the case of $T_1$, we have $A = 2 < 3 = B$, so $e_{1}$ will act
    on the left block of $+$'s in the $1$-signature; in this case, it
    will act on the leftmost $2$ in the top row of $T_{1}$.  By
    constrast, in $T_{2}$ we have $A = 2 = B$, so $e_{1}$ will act on
    the right block of $+$'s in the $1$-signature; this corresponds to
    the $\bar{1}$ in the top row.

    \label{ex:2lr}
    Let 
    $$
    T_{3}=
        \begin{array}{|c|c|c|c|c|c|c}
	\cline{1-6}
	\raisebox{-2pt}{$1$} & \raisebox{-2pt}{$1$} & 
	\raisebox{-2pt}{$2$} & \raisebox{-2pt}{$\bar{2}$} & 
	\raisebox{-2pt}{$\bar{2}$} & \raisebox{-2pt}{$\bar{1}$} & \\
	\cline{1-6}        
	\raisebox{-2pt}{$\bar{2}$} & \raisebox{-2pt}{$\bar{2}$} &
	\raisebox{-2pt}{$\bar{1}$} & 
	\multicolumn{4}{c}{} \\
	\cline{1-3} 
    \end{array}
    $$
    and    
    $$
    T_{4}=
        \begin{array}{|c|c|c|c|c|c|c}
	\cline{1-6}
	\raisebox{-2pt}{$1$} & \raisebox{-2pt}{$2$} & 
	\raisebox{-2pt}{$2$} & \raisebox{-2pt}{$\bar{2}$} & 
	\raisebox{-2pt}{$\bar{2}$} & \raisebox{-2pt}{$\bar{1}$} & \\
	\cline{1-6}        
	\raisebox{-2pt}{$\bar{2}$} & \raisebox{-2pt}{$\bar{2}$} &
	\raisebox{-2pt}{$\bar{1}$} & 
	\multicolumn{4}{c}{} \\
	\cline{1-3} 
    \end{array}.
    $$
    In the case of $T_{3}$, we have $C = 1 < 2 = D$, so $e_{2}$ will 
    act on the left block of $+$'s in the $2$-signature; in this 
    case, it will act on the leftmost $\bar{2}$ in the bottom row of 
    $T_{3}$.  By contrast, in $T_{4}$ we have $C = 2 = D$, so $e_{2}$ 
    will act on the right block of $+$'s in the $2$-signature; this 
    will be the leftmost $\bar{2}$ in the top row.
\end{example}

We now show that the entries in $T$ on which a sequence
$e_{i_{1}}\cdots e_{i_{\ell}}$ acts are determined by which blocks of
$+$'s correspond to those entries.  To achieve this, we verify that
the left edge of a block of $+$'s can be changed only by acting on
that block of $+$'s.  This goal motivates the following notation.

We write $e_{1}^{{\bf \ell}}$ to indicate the Kashiwara operator
$e_{1}$ when applied to a tableau $T$ such that $A(T) < B(T)$ and
$e_{1}^{{\bf r}}$ to indicate the Kashiwara operator $e_{1}$ when
applied to a tableau $T$ such that $A(T) \geq B(T)$.  Similarly, we
write $e_{2}^{{\bf \ell}}$ to indicate the Kashiwara operator $e_{2}$
when applied to a tableau $T$ such that $C(T) < D(T)$ and $e_{2}^{{\bf
r}}$ to indicate the Kashiwara operator $e_{2}$ when applied to a
tableau $T$ such that $C(T) \geq D(T)$.  Note that these are not new
operators: we simply use the superscript notation to record additional
information about how the operators act on specific tableaux.

\begin{claim} \label{prop:leftedge}
    Let $T$ be a tableau such that $e_{1}^{{\bf r}}T \neq 0$.  Then
    the left edges of both the left and right blocks of $+$'s in the
    $2$-signature of $e_{1}T$ and the left edge of the left block of
    $+$'s in the $1$-signature of $e_{1}T$ are in the same place as
    they are in the signatures of $T$, and the left edge of the right
    block of $+$'s in the $1$-signature of $e_{1}T$ is one position to
    the right of that in $T$.  
    
    Symmetrically, if $e_{1}^{\ell}T \neq 0$, the left edges of the
    blocks in the $2$-signature and the right block in the
    $1$-signature are unchanged and the left edge of the left block in
    the $1$-signature moves one position to the right; if $e_{2}^{{\bf
    r}}T \neq 0$, the left edges of the blocks in the $1$-signature
    and the left block in the $2$-signature are unchanged and the left
    edge of the right block in the $2$-signature moves one position to
    the right; and if $e_{2}^{\ell}T \neq 0$, the left edges of the
    blocks in the $1$-signature and the right block in the
    $2$-signature are unchanged and the left edge of the left block in
    the $2$-signature moves one position to the right.
\end{claim}

\begin{proof}
    Let $i = 1$ and $j = 2$ or vice versa, and let ${\bf x} = {\bf
    \ell}$ and ${\bf y} = {\bf r}$ or vice versa.  It is clear from
    the combinatorially defined action of $e_{i}$ on $C_{2}$ tableaux
    that if $e_{i}^{{\bf x}}T \neq 0$, the left edge of the {\bf y}
    block of $+$'s in the $i$-signature of $e_{i}T$ is in the same
    place as in the $i$-signature of $T$, and that the left edge of
    the {\bf x} block of $+$'s in the $i$-signature of $e_{i}T$ is one
    space to the right of its position in the $i$-signature of $T$.
    We may therefore devote our attention to the $j$-signature in each
    of the four cases of concern.

    First, consider the case of $e_{1}^{{\bf \ell}}T\neq 0$.  By Claim
    \ref{claim:lr}, we know that this operator changes a $2$ in the
    top row to a $1$ or a $\bar{1}$ in the bottom row to a $\bar{2}$.
    In the first case, no $+$'s are added to the $2$-signature, and no
    change is made to those entries of interest to the location of a
    block of $+$'s of size $0$ in the $2$-signature.  In the other
    case, one $+$ is added to the left block of $+$'s in the
    $2$-signature; this addition is to the right of the left edge of
    this block.  Finally, the right block of $+$'s in the
    $2$-signature of $e_{1}T$ is the same as in the $2$-signature of
    $T$ in any case.

    Next, consider the case of $e_{1}^{{\bf r}}T\neq 0$.  By Claim
    \ref{claim:lr}, we know that this operator changes a $\bar{1}$ in the
    top row of $T$ into a $\bar{2}$.  This adds one $+$ to the right block
    of $+$'s in the $2$-signature to the right of its left edge and makes
    no change to the left block of $+$'s.

    Now, consider the case of $e_{2}^{{\bf \ell}}T\neq 0$.  By Claim
    \ref{claim:lr}, we know that this operator changes a $\bar{2}$ in the
    bottom row to a $2$.  This does not contribute a $+$ to either the
    left or right blocks of the $1$-signature of $e_{2}T$, nor does it
    pertain to the location of a block of $+$'s of size $0$ in the
    $1$-signature, so the left edges in this signature are the same as in
    the $1$-signature of $T$.

    Finally, we consider the case of $e_{2}^{{\bf r}}T\neq 0$.  By Claim
    \ref{claim:lr}, we know that this operator changes a $\bar{2}$ in the
    top row to a $2$.  This has the effect of adding a $+$ to the left
    block of $+$'s in the $1$-signature to the right of the left edge of
    this block.  Finally, the right block of $+$'s in the $1$-signature of
    $e_{2}T$ is the same as in the $1$-signature of $T$.

\end{proof}

\begin{corollary}
    Let $T$ be a tableau such that $e_{1}^{\ell}T \neq 0$, and let $E$
    be a sequence of operators from the set $\{ e_{1}^{{\bf r}},
    e_{2}^{\ell}, e_{2}^{{\bf r}} \}$ such that $ET \neq 0$.  Then
    $e_{1}^{\ell}$ acts on the same entry in $T$ as it does in $ET$.
    The symmetric statements corresponding to the cases of Claim
    \ref{prop:leftedge} hold as well.
\end{corollary}

\begin{example}
    Let 
    $$
    T=
    \begin{array}{|c|c|c}
	\cline{1-2}
	\raisebox{-2pt}{$\bar{2}$} & \raisebox{-2pt}{$\bar{1}$} & \\
	\cline{1-2}        
	\raisebox{-2pt}{$\bar{1}$} & 	\multicolumn{2}{c}{} \\
	\cline{1-1} 
    \end{array},
    $$
    so 
    $$
    e_{1}^{{\bf r}}T = 
    \begin{array}{|c|c|c}
	\cline{1-2}
	\raisebox{-2pt}{$\bar{2}$} & \raisebox{-2pt}{$\bar{2}$} & \\
	\cline{1-2}        
	\raisebox{-2pt}{$\bar{1}$} & 	\multicolumn{2}{c}{} \\
	\cline{1-1} 
    \end{array}
    $$	
    is the result of acting on the third tensor factor in the column 
    word of $T$.  We find that 
    $$
    e_{2}^{\ell}e_{1}^{\ell}e_{2}^{{\bf r}}T = 
    \begin{array}{|c|c|c}
	\cline{1-2}
	\raisebox{-2pt}{$1$} & \raisebox{-2pt}{$\bar{1}$} & \\
	\cline{1-2}        
	\raisebox{-2pt}{$\bar{2}$} & 	\multicolumn{2}{c}{} \\
	\cline{1-1} 
    \end{array},
    $$
    so that indeed we have
    $$
    e_{1}^{{\bf r}}e_{2}^{\ell}e_{1}^{\ell}e_{2}^{{\bf r}}T = 
    \begin{array}{|c|c|c}
	\cline{1-2}
	\raisebox{-2pt}{$1$} & \raisebox{-2pt}{$\bar{2}$} & \\
	\cline{1-2}        
	\raisebox{-2pt}{$\bar{2}$} & 	\multicolumn{2}{c}{} \\
	\cline{1-1} 
    \end{array},
    $$
    which differs from $e_{2}^{\ell}e_{1}^{\ell}e_{2}^{{\bf r}}T$ in 
    the third tensor factor.
\end{example}

The following four Sublemmas state that the relative values of $A(T),
B(T), C(T),$ and $D(T)$ not only determine where $e_{i}$ acts within a
tableau, but also what the values of $A(e_{i}T), B(e_{i}T),
C(e_{i}T),$ and $D(e_{i}T)$ are.  This will be an invaluable tool for
our analysis in section \ref{sec:local:proof}.

\begin{sublemma} \label{sublemma:alb}
    Suppose $T$ is a tableau such that $e_{1}$ acts on the left block
    of $+$'s in the $1$-signature of $T$ (i.e., such that $A(T) <
    B(T)$).  Then $A(e_{1}T) = A(T)$, $B(e_{1}T) = B(T) - 1$, and
    $C(e_{1}T) - D(e_{1}T) = C(T) - D(T) - 1$.
\end{sublemma}

\begin{proof}
    
We have two cases to consider; $e_{1}$ may act by changing a $2$ to a
$1$ in the top row or a $\bar{1}$ to a $\bar{2}$ in the bottom row.
In both of these cases, it is easy to see that the number of
$\bar{2}$'s in the top row is unchanged and the number of $2$'s in the
top row plus the number of $\bar{1}$'s in the bottom row is diminished
by one; hence $A(e_{1}T) = A(T)$ and $B(e_{1}T) = B(T) - 1$.

Observe that in the case of a $\bar{1}$ changing into a $\bar{2}$ in
the bottom row, the content of the top row is unchanged, but the
number of $\bar{2}$'s in the bottom row is increased by $1$.  In the
case of a $2$ changing into a $1$ in the top row, the bottom row is
unchanged, but the number of $2$'s in the top row is decreased by $1$.
In both of these cases, we find that $C(e_{1}T) - D(e_{1}T) = C(T) -
D(T) - 1$.

\end{proof}

\begin{sublemma} \label{sublemma:agb}
    Suppose $T$ is a tableau such that $e_{1}$ acts on the right block
    of $+$'s in the $1$-signature of $T$ (i.e., such that $A(T) \geq
    B(T)$).  Then $A(e_{1}T) = A(T) + 1$, $B(e_{1}T) = B(T)$,
    $C(e_{1}T) = C(T)$, and $D(e_{1}T) = D(T)$.
\end{sublemma}

\begin{proof}

Since the right block of $+$'s in the $1$-signature comes entirely
from $\bar{1}$'s in the top row of $T$, it follows that acting by
$e_{1}$ changes one of these $\bar{1}$'s into a $\bar{2}$.  We
immediately see that the number of $\bar{2}$'s in the top row
increases by $1$, and that the number of $2$'s in the top row and
$\bar{2}$'s and $\bar{1}$'s in the bottom row are all unchanged.

\end{proof}

\begin{sublemma} \label{sublemma:cld}
    Suppose $T$ is a tableau such that $e_{2}$ acts on the left block
    of $+$'s in the $1$-signature of $T$ (i.e., such that $C(T) <
    D(T)$).  Then $A(e_{2}T) = A(T)$, $B(e_{2}T) = B(T)$, $C(e_{2}T) =
    C(T)$, and $D(e_{2}T) = D(T) - 1$.
\end{sublemma}

\begin{proof}

The entry on which $e_{2}$ acts is a $\bar{2}$ the bottom row, which
will be changed into a $2$.  We immediately see that the number of
$\bar{2}$'s in the bottom row decreases by $1$, and that the number of
$2$'s and $\bar{2}$'s in the top row and $\bar{1}$'s in the bottom row
are all unchanged.

\end{proof}

\begin{sublemma} \label{sublemma:cgd}
    Suppose $T$ is a tableau such that $e_{2}$ acts on the right block
    of $+$'s in the $1$-signature of $T$ (i.e., such that $C(T) \geq
    D(T)$).  Then $A(e_{2}T) = A(T) - 1$, $B(e_{2}T) = B(T) + 1$,
    $C(e_{2}T) = C(T) + 1$, and $D(e_{2}T) = D(T)$.
\end{sublemma}

\begin{proof}

In this case $e_{2}$ will change a $\bar{2}$ to a $2$ in the top row.
It is easy to see that the number of $2$'s in the top row increases by
$1$ and the number of $\bar{2}$'s in the bottom row is unchanged;
hence $C(e_{2}T) = C(T) + 1$ and $D(e_{2}T) = D(T)$.
  
Likewise, since the number of $\bar{2}$'s in the top row is decreased
by $1$ and the number of $2$'s in the top row is increased by $1$, we
find that $A(e_{2}T) = A(T) - 1$ and $B(e_{2}T) = B(T) + 1$.

\end{proof}

        \section{Proof of Theorem \ref{proposition:int}}
        \label{sec:local:proof}

We are now equipped to address Theorem \ref{proposition:int}.  It is
proved as a consequence of Lemmas \ref{lemma:deg2:1} through
\ref{lemma:deg7}, each of which deals with a certain case of the
relative values of $A(T), B(T), C(T),$ and $D(T)$.  To see that these
cases are exhaustive, refer to Table \ref{relstats}.

\begin{table}
    $$
    \begin{array}{c|c|c|c|c|}
	 & A<B & A=B & A=B+1 & A>B+1  \\
	\hline
	 C<D & 2 & 2 & 2 & 2  \\
	\hline
	 C=D & 4 & 7 & 4 & 2  \\
	\hline
	 C>D & 2 & 5 & 4 & 2  \\
	\hline
	 
    \end{array}
    $$

    \caption{Degree of relation over $T$, given its $ABCD$ statistics} 
    \protect\label{relstats}
\end{table}


\begin{lemma}
    \label{lemma:deg2:1}
    Suppose $T$ is a tableau such that $C(T) < D(T)$, $e_{1}T\neq 0$, and
    $e_{2}T\neq 0$.  Then $T$ has a degree 2 relation above it.
\end{lemma}

\begin{proof}

    From Claim \ref{claim:lr}, we know that $e_{2}$ acts on the left
    block of $+$'s, and by Sublemma \ref{sublemma:cld}, we know that
    $A(e_{2}T) = A(T)$ and $B(e_{2}T) = B(T)$; it follows that
    $e_{1}e_{2}T \neq 0$, and that $e_{1}$ acts on the same entry in
    $e_{2}T$ as it does in $T$.  Furthermore, by Sublemmas
    \ref{sublemma:alb} and \ref{sublemma:agb}, we know that either
    $C(e_{1}T) - D(e_{1}T) = C(T) - D(T) - 1$ or $C(e_{1}T) = C(T)$
    and $D(e_{1}T) = D(T)$; in either case, we still find that
    $C(e_{1}T) < D(e_{1}T)$.  Since $C(e_{1}T) \geq 0$, we are assured
    that $D(e_{1}T) \geq 1$, and thus $e_{2}e_{1}T \neq 0$.  We
    conclude that $e_{2}$ acts on the same entry in $e_{1}T$ as it
    does in $T$.
    
\end{proof}
\psset{unit=1cm}

\begin{example}
    
    In Figure \ref{fig:lem1}, we have an interval of a crystal in 
    which the bottom tableau $T$ has the statistics $C(T) = 1 < 2 = 
    D(T)$, illustrating Lemma \ref{lemma:deg2:1}.
    
    \begin{figure}[!hbp]
	
	\begin{pspicture}(-7.75,-5.25cm)
	    \rput(0,-.5){    
	    $
	    \begin{array}{|c|c|c|c|c}
		\cline{1-4}
		\raisebox{-1pt}{$1$} & \raisebox{-1pt}{$1$} & 
		\raisebox{-1pt}{$\bar{1}$} & \raisebox{-1pt}{$\bar{1}$} & \\
		\cline{1-4}        
		\raisebox{-1pt}{$2$} & \raisebox{-1pt}{$\bar{2}$} & 
		\multicolumn{3}{c}{} \\
		\cline{1-2} 
	    \end{array}
	    $
	    }
	    
	    \rput(-3,-2.75){    
	    $
	    \begin{array}{|c|c|c|c|c}
		\cline{1-4}
		\raisebox{-1pt}{$1$} & \raisebox{-1pt}{$1$} & 
		\raisebox{-1pt}{$\bar{1}$} & \raisebox{-1pt}{$\bar{1}$} & \\
		\cline{1-4}        
		\raisebox{-1pt}{$\bar{2}$} & \raisebox{-1pt}{$\bar{2}$} & 
		\multicolumn{3}{c}{} \\
		\cline{1-2} 
	    \end{array}
	    $
	    }
	    
	    \rput(3,-2.75){    
	    $
	    \begin{array}{|c|c|c|c|c}
		\cline{1-4}
		\raisebox{-1pt}{$1$} & \raisebox{-1pt}{$2$} & 
		\raisebox{-1pt}{$\bar{1}$} & \raisebox{-1pt}{$\bar{1}$} & \\
		\cline{1-4}        
		\raisebox{-1pt}{$2$} & \raisebox{-1pt}{$\bar{2}$} & 
		\multicolumn{3}{c}{} \\
		\cline{1-2} 
	    \end{array}
	    $
	    }	
	    
	    \rput(0,-5){    
	    $
	    \begin{array}{|c|c|c|c|c}
		\cline{1-4}
		\raisebox{-1pt}{$1$} & \raisebox{-1pt}{$2$} & 
		\raisebox{-1pt}{$\bar{1}$} & \raisebox{-1pt}{$\bar{1}$} & \\
		\cline{1-4}        
		\raisebox{-1pt}{$\bar{2}$} & \raisebox{-1pt}{$\bar{2}$} & 
		\multicolumn{3}{c}{} \\
		\cline{1-2} 
	    \end{array}
	    $
	    }
	    
	    \psline{->}(-1.5,-.9)(-3,-2)
	    \rput[Br](-2.25,-1.45){$2\;$}
	    
	    \psline{->}(1.5,-.9)(3,-2)
	    \rput[Bl](2.25,-1.45){$\;1$}
	    
	    \psline{->}(-3,-3.15)(-1.5,-4.25)
	    \rput[Bl](-2.25,-3.7){$\;1$}
	    
	    \psline{->}(2.8,-3.35)(1.3,-4.45)
	    \rput[Br](2.05,-3.9){$2\;$}	
	    
	\end{pspicture}	
	
        \caption{}
        \protect\label{fig:lem1}
    \end{figure}

\end{example}


\begin{lemma}
    \label{lemma:deg2:2}
    Suppose $T$ is a tableau such that $A(T) > B(T) + 1$, $e_{1}T\neq
    0$, and $e_{2}T\neq 0$.  Then $T$ has a degree 2 relation above
    it.
\end{lemma}

\begin{proof}
    
    From Claim \ref{claim:lr}, we know that $e_{1}$ acts on the right
    block of $+$'s, and by Sublemma \ref{sublemma:agb}, we know that
    $C(e_{1}T) = C(T)$ and $D(e_{1}T) = D(T)$; thus $e_{2}e_{1}T \neq
    0$, and $e_{2}$ acts on the same entry in $e_{1}T$ as it does in
    $T$.  Furthermore, by Sublemmas \ref{sublemma:cld} and
    \ref{sublemma:cgd}, we know that either $A(e_{2}T) = A(T)$ and
    $B(e_{2}T) = B(T)$ or $A(e_{2}T) = A(T) - 1$ and $B(e_{2}T) = B(T)
    + 1$; in either case, we find that $A(e_{2}T) > B(e_{2}T)$ and the
    size of the right block of $+$'s in the $1$-signature is not
    diminished.  We therefore conclude that $e_{1}e_{2}T \neq 0$ and
    that $e_{1}$ acts on the same entry in $e_{2}T$ as it does in
    $T$.

\end{proof}

\begin{example}
    
    In Figure \ref{fig:lem2}, we have an interval of a crystal in 
    which the bottom tableau $T$ has the statistics $A(T) = 2 > 1 = 
    B(T)+1$ illustrating Lemma \ref{lemma:deg2:2}.
    
    \begin{figure}[!hbp]
	
	\begin{pspicture}(-7.5,-5.25cm)
	    \rput(0,-.5){    
	    $
	    \begin{array}{|c|c|c|c|c|}
		\cline{1-5}
		\raisebox{-1pt}{$1$} & \raisebox{-1pt}{$1$} & 
		\raisebox{-1pt}{$2$} & \raisebox{-1pt}{$\bar{2}$} & 
		\raisebox{-1pt}{$\bar{2}$} \\
		\cline{1-5}        
	    \end{array}
	    $
	    }
	    
	    \rput(-3,-2.75){    
	    $
	    \begin{array}{|c|c|c|c|c|}
		\cline{1-5}
		\raisebox{-1pt}{$1$} & \raisebox{-1pt}{$1$} & 
		\raisebox{-1pt}{${2}$} & \raisebox{-1pt}{$\bar{2}$} & 
		\raisebox{-1pt}{$\bar{2}$} \\
		\cline{1-5}        
	    \end{array}
	    $
	    }
	    
	    \rput(3,-2.75){    
	    $
	    \begin{array}{|c|c|c|c|c|}
		\cline{1-5}
		\raisebox{-1pt}{$1$} & \raisebox{-1pt}{$1$} & 
		\raisebox{-1pt}{$2$} & \raisebox{-1pt}{$\bar{2}$} & 
		\raisebox{-1pt}{$\bar{1}$} \\
		\cline{1-5}        
	    \end{array}
	    $
	    }	
	    
	    \rput(0,-5){    
	    $
	    \begin{array}{|c|c|c|c|c|}
		\cline{1-5}
		\raisebox{-1pt}{$1$} & \raisebox{-1pt}{$1$} & 
		\raisebox{-1pt}{$\bar{2}$} & \raisebox{-1pt}{$\bar{2}$} & 
		\raisebox{-1pt}{$\bar{1}$} \\
		\cline{1-5}        
	    \end{array}
	    $
	    }
	    
	    \psline{->}(-1.5,-.9)(-3,-2.2)
	    \rput[Br](-2.25,-1.55){$2\;$}
	    
	    \psline{->}(1.5,-.9)(3,-2.2)
	    \rput[Bl](2.25,-1.55){$\;1$}
	    
	    \psline{->}(-3,-3.15)(-1.5,-4.55)
	    \rput[Bl](-2.25,-3.85){$\;1$}
	    
	    \psline{->}(3,-3.15)(1.5,-4.55)
	    \rput[Br](2.25,-3.85){$2\;$}	
	    
	\end{pspicture}

        \caption{}
        \protect\label{fig:lem2}
    \end{figure}

\end{example}


\begin{lemma}\label{lemma:deg2:3}
    Suppose $T$ is a tableau such that $A(T) < B(T)$, $C(T) > D(T)$,
    $e_{1}T\neq 0$, and $e_{2}T\neq 0$.  Then $T$ has a degree 2
    relation above it.
\end{lemma}

\begin{proof}
    
    By Claim \ref{claim:lr}, we know that $e_{1}$ acts on the left
    block of $+$'s in $T$ and $e_{2}$ acts on the right block of $+$'s
    in $T$.  By Sublemma \ref{sublemma:cgd}, we know that $A(e_{2}T) =
    A(T) - 1$ and $B(e_{2}T) = B(T) + 1$.  It follows that $A(e_{2}T)
    < B(e_{2}T)$, and since $A(e_{2}T) \geq 0$, this ensures that
    $B(e_{2}T) \geq 1$, and thus $e_{1}e_{2}T \neq 0$.  We conclude
    that $e_{1}$ acts on the same entry in $e_{2}T$ as it does in
    $T$.  Furthermore, by Sublemma \ref{sublemma:alb}, we know that
    $C(e_{1}T) - D(e_{1}T) = C(T) - D(T) - 1$; it follows that
    $C(e_{1}T) \geq D(e_{1}T)$.  Since we also know that the size of
    the right block of $+$'s in the $2$-signature of $e_{1}T$ is at
    least as large as that of $T$, it is the case that $e_{2}e_{1}T
    \neq 0$, and so $e_{2}$ acts on the same entry in $e_{1}T$ as it
    does in $T$.

\end{proof}

\begin{example}

    In Figure \ref{fig:lem3}, we have an interval of a crystal in
    which the bottom tableau $T$ has the statistics $A(T) = 1 < 2 =
    B(T)$ and $C(T) = 1 > 0 = D(T)$, illustrating Lemma
    \ref{lemma:deg2:3}.
    
    \begin{figure}[!hbp]
	\begin{pspicture}(-7.75,-5.25cm)
	    \rput(0,-.5){    
	    $
	    \begin{array}{|c|c|c}
		\cline{1-2}
		\raisebox{-1pt}{$2$} & \raisebox{-1pt}{$2$} & 
		\\
		\cline{1-2}        
		\raisebox{-1pt}{$\bar{2}$} & 
		\multicolumn{2}{c}{} \\
		\cline{1-1} 
	    \end{array}
	    $
	    }
	    
	    \rput(-3,-2.75){    
	    $
	    \begin{array}{|c|c|c}
		\cline{1-2}
		\raisebox{-1pt}{$2$} & \raisebox{-1pt}{$\bar{2}$} & 
		\\
		\cline{1-2}        
		\raisebox{-1pt}{$\bar{2}$} & 
		\multicolumn{2}{c}{} \\
		\cline{1-1} 
	    \end{array}
	    $
	    }
	    
	    \rput(3,-2.75){    
	    $
	    \begin{array}{|c|c|c}
		\cline{1-2}
		\raisebox{-1pt}{$2$} & \raisebox{-1pt}{$2$} & 
		\\
		\cline{1-2}        
		\raisebox{-1pt}{$\bar{1}$} & 
		\multicolumn{2}{c}{} \\
		\cline{1-1} 
	    \end{array}
	    $
	    }	
	    
	    \rput(0,-5){    
	    $
	    \begin{array}{|c|c|c}
		\cline{1-2}
		\raisebox{-1pt}{$2$} & \raisebox{-1pt}{$\bar{2}$} & 
		\\
		\cline{1-2}        
		\raisebox{-1pt}{$\bar{1}$} & 
		\multicolumn{2}{c}{} \\
		\cline{1-1} 
	    \end{array}
	    $
	    }
	    
	    \psline{->}(-1.1,-.9)(-3,-2)
	    \rput[Br](-2.05,-1.45){$2\;$}
	    
	    \psline{->}(.8,-.9)(3,-2)
	    \rput[Bl](1.9,-1.45){$\;1$}
	    
	    \psline{->}(-3,-3.15)(-1.1,-4.25)
	    \rput[Bl](-2.05,-3.7){$\;1$}
	    
	    \psline{->}(2.8,-3.35)(.8,-4.45)
	    \rput[Br](1.8,-3.9){$2\;$}	
	    
	\end{pspicture}
        \caption{}
        \protect\label{fig:lem3}
    \end{figure}
\end{example}


To prove Lemmas \ref{lemma:deg4:1} through \ref{lemma:deg7}, we must
not only show that the given sequences of operators act on the same
entries, but also that no pair of homogeneous sequences of operators
(i.e., a pair $(P_{1},P_{2})$ such that $P_{1}$ and $P_{2}$ have the
same number of instances of $e_{1}$ and $e_{2}$) with shorter or equal
length act on the same entries.  To assist in our illustration of this
fact, we will refer to figures that encode the generic behavior of all
sequences of operators on a tableau with content as specified by the
hypothesis of each lemma.  Table \ref{piclegend} is a legend for the
figures used to prove Lemmas \ref{lemma:deg4:1} through
\ref{lemma:deg5}.  In the picture used to prove Lemma
\ref{lemma:deg7}, we instead use an edge pointing down to indicate
acting by $e_{1}$ and an edge pointing up to indicate acting by
$e_{2}$; otherwise the legend is the same.

To assist in proving that the sequences in question do not kill our
tableaux, we have the following Sublemma.

\begin{sublemma}
    \label{sublemma:dashed}
    Let $E$ be a dashed edge from $v$ up to $w$; i.e., an operator
    $e_{i}^{\ell}$ acts on $v$ to produce $w$.  Then $e_{i}v \neq 0$.
\end{sublemma}

\begin{proof}
    
    The Kashiwara operator $e_{i}$ acts on the left block of $+$'s of
    a tableau $T$ precisely when $A(T) < B(T)$ or $C(T) < D(T)$ in the
    cases of $i = 1$ or $i = 2$, respectively.  Since these numbers
    are all non-negative integers, we conclude that $B(T) > 0$ or
    $D(T) > 0$.  Since these statistics indicate the number of $+$'s
    in the left block of their respective signatures, we are assured
    that there is an entry on which $e_{i}$ can act.
    
\end{proof}

Thus it suffices to prove that the solid edges in the paths of concern
do not produce $0$.

\begin{table}
    
    $
    \begin{array}{|c|c|}
	\hline
	\textrm{edge pointing to the right} & \textrm{acting by }e_{1}  \\
	\hline
	\textrm{edge pointing to the left} & \textrm{acting by }e_{2}  \\
	\hline
	\textrm{solid edge} & \textrm{acting on the right block of }+\textrm{'s}  \\
	\hline
	\textrm{dashed edge} & \textrm{acting on the left block of }+\textrm{'s}  \\
	\hline
	\textrm{vertex labeled by }(d_{1},d_{2}) & \textrm{tableau
	}T\textrm{ with statistics such that }\\&A(T) = B(T) +
	d_{1}\textrm{ and }C(T) = D(T) + d_{2} \\
	\hline
	 
    \end{array}
    $

    \caption{Legend for Figures \ref{piclem4} through \ref{piclem7}}
    \protect\label{piclegend}
\end{table}


\begin{lemma}
    \label{lemma:deg4:1}
    Suppose $T$ is a tableau such that $A(T) = B(T) + 1$, $C(T) \geq
    D(T)$, $e_{1}T\neq 0$, and $e_{2}T\neq 0$.  Then $T$ has a degree
    4 relation above it.
\end{lemma}

\begin{figure}
    \centerline{\includegraphics[height=2.5in]{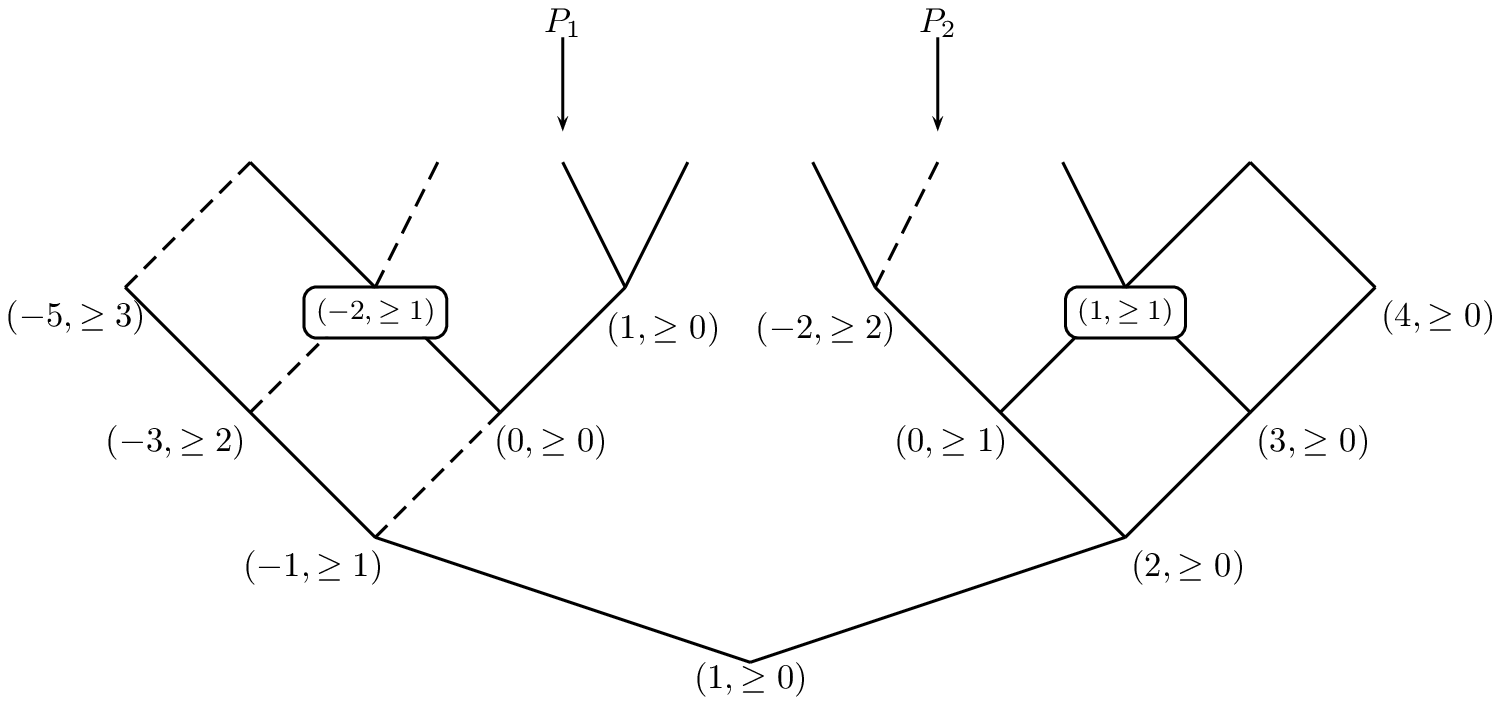}}
    \caption{Picture for Lemma \ref{lemma:deg4:1}}
    \protect\label{piclem4}
\end{figure}

\begin{proof}

    We must first confirm that the sequences $e_{1}e_{2}^{2}e_{1}$ and
    $e_{2}e_{1}^{2}e_{2}$ do not produce $0$ when applied to $T$.
    First, observe that since $e_{1}$ acts on the right block of $+$'s
    in $T$, it changes a $\bar{1}$ to a $\bar{2}$ in the top row.
    This adds a $+$ to the reduced $2$-signature of the tableaux, so
    we know that $e_{2}^{2}e_{1}T \neq 0$.  By Sublemma
    \ref{sublemma:dashed}, we know that $e_{1}e_{2}^{2}e_{1}T \neq 0$.
    On the other hand, we know that $e_{2}$ acts on $T$ by changing a
    $\bar{2}$ to a $2$ in the top row; this means that the reduced
    $1$-signature of $e_{2}T$ has a single $+$ in the left block and
    its right block has at least one $+$, as did the $1$-signature of
    $T$.  We conclude that $e_{1}^{2}e_{2}T \neq 0$.  We know that in
    $e_{1}e_{2}T$, $e_{1}$ changes a $\bar{1}$ to a $\bar{2}$ in the
    top row; since the $+$'s in the $2$-signature from this entry
    cannot be paired with any $-$'s, we conclude that
    $e_{2}e_{1}^{2}e_{2}T \neq 0$.

    Now that we know that neither of these sequences produces $0$ when
    applied to $T$, it is clear that we have $e_{1}e_{2}^{2}e_{1}T =
    e_{2}e_{1}^{2}e_{2}T$, as the paths $P_{1}$ and $P_{2}$ leading
    from the base of the graph in Figure \ref{piclem4} to the
    indicated leaves both have one solid right edge, one dashed right
    edge, and two solid left edges.  We must now confirm that among
    all pairs $(Q_{1},Q_{2})$ of increasing paths from the base in
    these graphs such that $Q_{1}$ begins by following the left edge
    and $Q_{2}$ begins by following the right edge, $(P_{1}, P_{2})$
    is the only pair with the same number of each type of edge.

    Since the right edge from the base of the graph is solid, our
    candidate for $Q_{1}$ must have a solid right edge.  Inspecting
    the graph tells us that this path must begin with the path
    corresponding to $e_{1}^{2}e_{2}$.  This path has a dashed left
    edge, and the only candidate for $Q_{2}$ with this feature is in
    fact $P_{2}$, which has two solid right edges.  The only way to
    extend $e_{1}^{2}e_{2}$ to have the same edge content as $P_{2}$
    is by extending it to $P_{1}$.

\end{proof}

\begin{example}
    
    In Figure \ref{fig:lem4}, we have an interval of a crystal in
    which the bottom tableau $T$ has the statistics $A(T) = 2 = B(T) +
    1$ and $C(T) = 1 > 0 = D(T)$, illustrating Lemma
    \ref{lemma:deg4:1}.
    
    \begin{figure}[!hbp]
        \centerline{\includegraphics[height=2.95in]{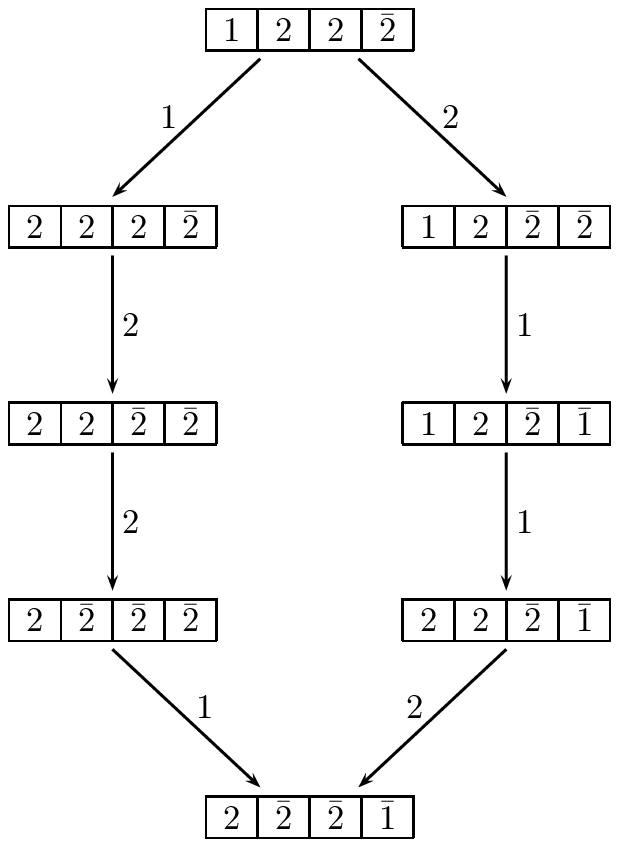}}
        \caption{}
        \protect\label{fig:lem4}
    \end{figure}
\end{example}


\begin{lemma}
    \label{lemma:deg4:2}
    Suppose $T$ is a tableau such that $A(T) < B(T)$, $C(T) = D(T)$,
    $e_{1}T\neq 0$, and $e_{2}T\neq 0$.  Then $T$ has a degree 4
    relation above it.
\end{lemma}

\begin{figure}
    \centerline{\includegraphics[height=2.5in]{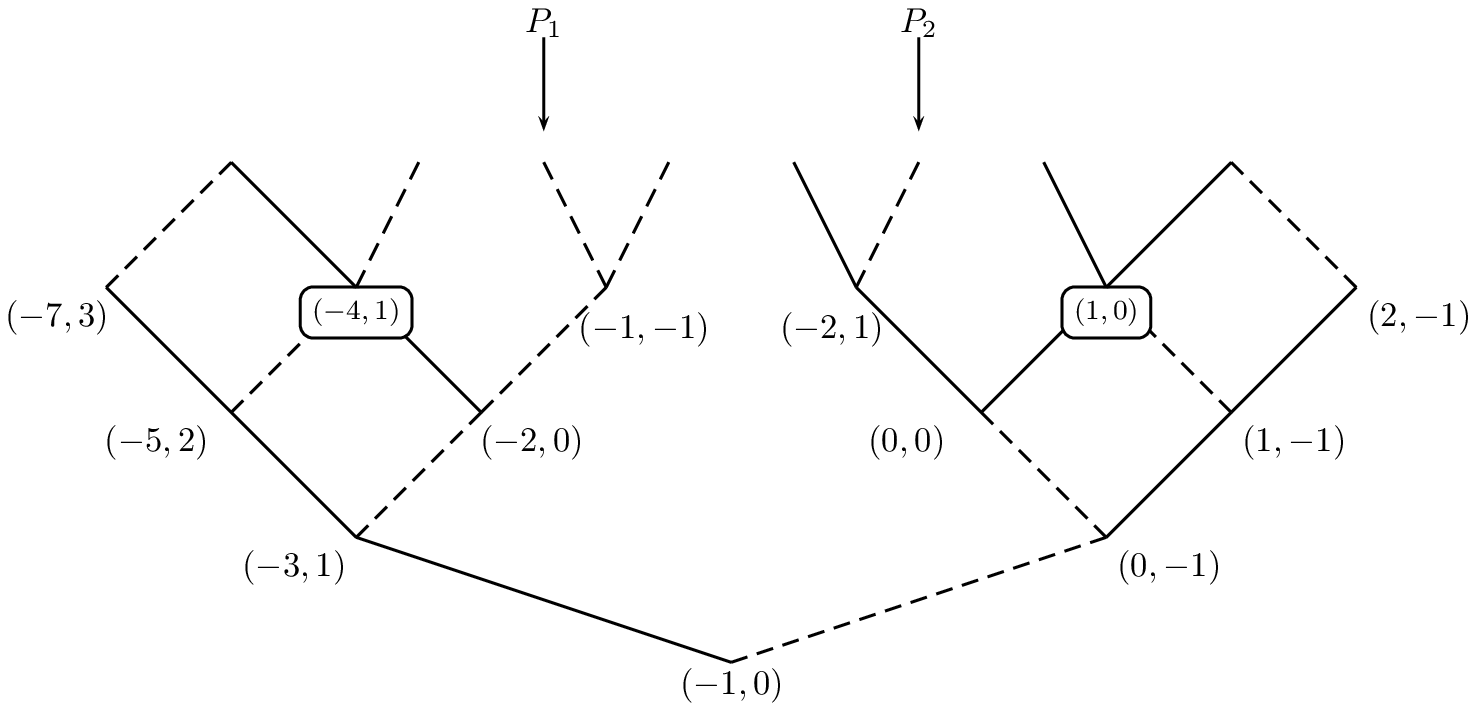}}
    \caption{Picture 1 for Lemma \ref{lemma:deg4:2}}
    \protect\label{pic1lem5}
\end{figure}

\begin{figure}
    \centerline{\includegraphics[height=2.5in]{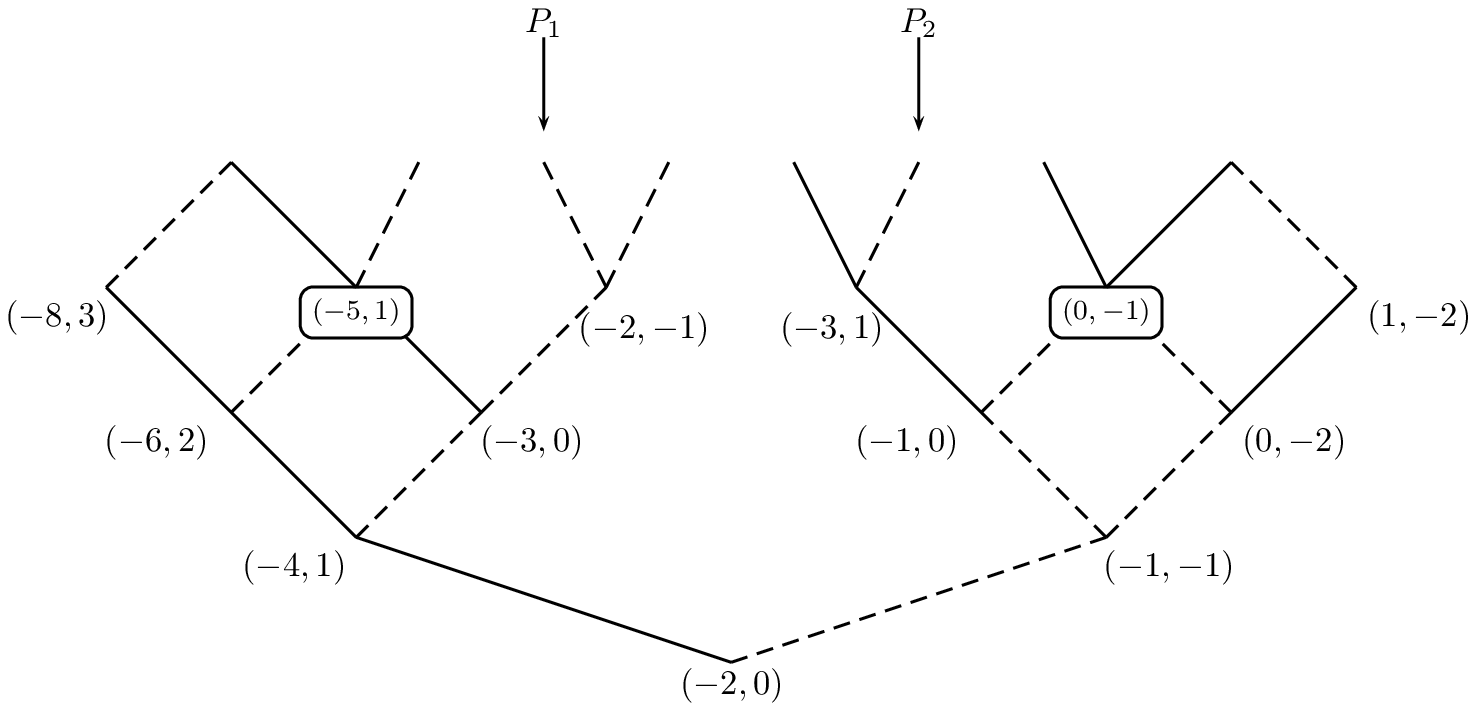}}
    \caption{Picture 2 for Lemma \ref{lemma:deg4:2}}
    \protect\label{pic2lem5}
\end{figure}

\begin{figure}
    \centerline{\includegraphics[height=2.5in]{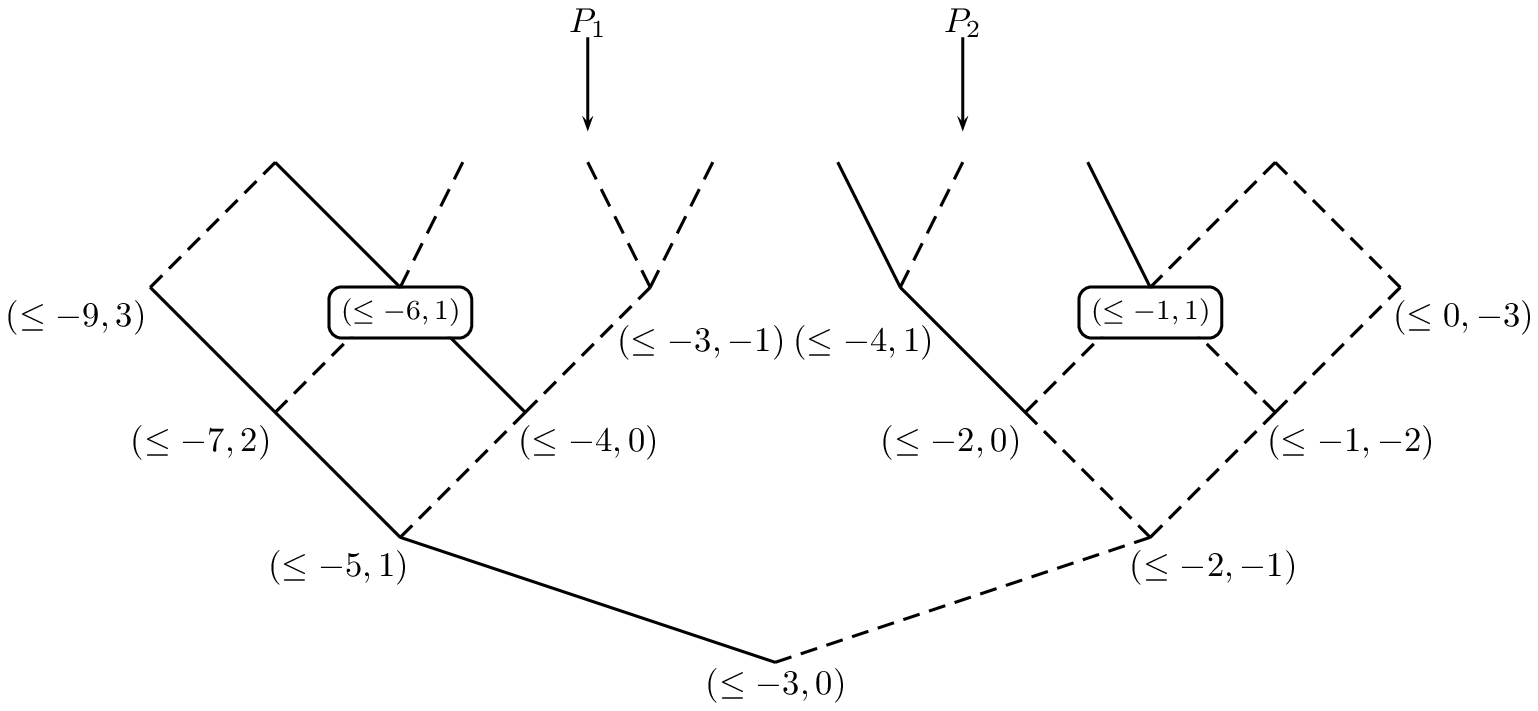}}
    \caption{Picture 3 for Lemma \ref{lemma:deg4:2}}
    \protect\label{pic3lem5}
\end{figure}

\begin{proof}

    We must first confirm that the sequences $e_{1}e_{2}^{2}e_{1}$ and
    $e_{2}e_{1}^{2}e_{2}$ do not produce $0$ when applied to $T$.  By
    Sublemma \ref{sublemma:dashed} it suffices to show that
    $e_{2}^{2}e_{1}T \neq 0$, since $e_{2}T \neq 0$ by assumption and
    all other edges in $P_{1}$ and $P_{2}$ are dashed.  To see this,
    simply observe that there is at least one $+$ in the right block
    of the reduced $2$-signature of $T$; the sequence $e_{2}e_{1}$
    acts on the left blocks of $+$'s, so the corresponding entry
    remains available for $e_{2}$ to act on.

    Now that we know that neither of these sequences produces $0$ when
    applied to $T$, it is clear that we have $e_{1}e_{2}^{2}e_{1}T =
    e_{2}e_{1}^{2}e_{2}T$, as the paths $P_{1}$ and $P_{2}$ leading
    from the base of the graph in Figures \ref{pic1lem5} through
    \ref{pic3lem5} to the indicated leaves have two dashed right
    edges, one solid left edge, and one dashed left edge.  We must now
    confirm that among all pairs $(Q_{1},Q_{2})$ of increasing paths
    from the base in these graphs such that $Q_{1}$ begins by
    following the left edge and $Q_{2}$ begins by following the right
    edge, $(P_{1}, P_{2})$ is the only pair with the same number of
    each type of edge.

    This is easy to see by the following argument.  Every candidate
    for $Q_{2}$ (i.e., every path in the right half of the graphs in
    Figures \ref{pic1lem5} through \ref{pic3lem5}) has at least one
    dashed left edge.  The only candidate for $Q_{1}$ (i.e., the only
    path in the left half of the graphs in Figures \ref{pic1lem5}
    through \ref{pic3lem5}) with a dashed left edge is $P_{1}$.  By
    inspecting Figures \ref{pic1lem5} through \ref{pic3lem5}, $P_{2}$
    is the only candidate for $Q_{2}$ with two dashed right edges and
    one solid left edge.

\end{proof}

\begin{example}
    
    In Figure \ref{fig:lem5}, we have an interval of a crystal in
    which the bottom tableau $T$ has the statistics $A(T) = 1 < 2 = B(T)$
    and $C(T) = 1 = D(T)$, illustrating Lemma \ref{lemma:deg4:2}.
    
    \begin{figure}[!hbp]
        \centerline{\includegraphics[height=2.95in]{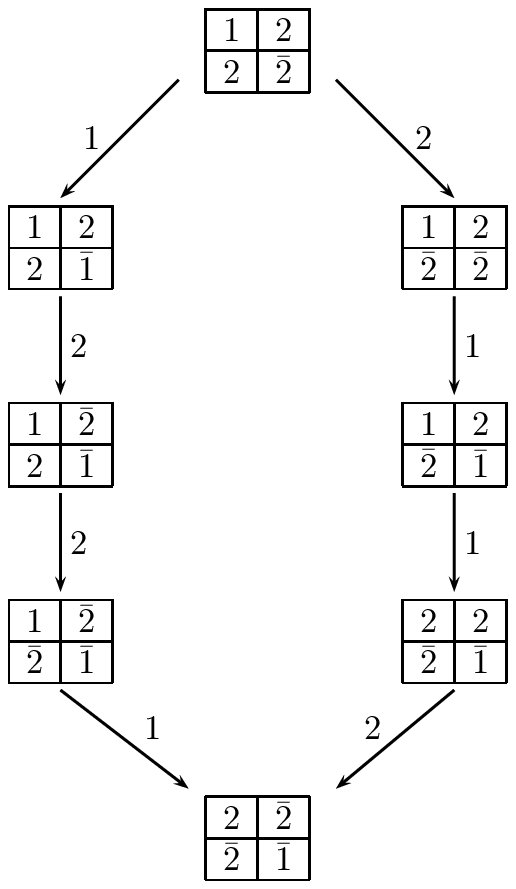}}
        \caption{}
        \protect\label{fig:lem5}
    \end{figure}
\end{example}


\begin{lemma}
    \label{lemma:deg5}
    Suppose $T$ is a tableau such that $A(T) = B(T)$, $C(T) \geq D(T)
    + 1$, $e_{1}T\neq 0$, and $e_{2}T\neq 0$.  Then $T$ has a degree 5
    relation above it.
\end{lemma}

\begin{figure}
    \centerline{\includegraphics[height=2.5in]{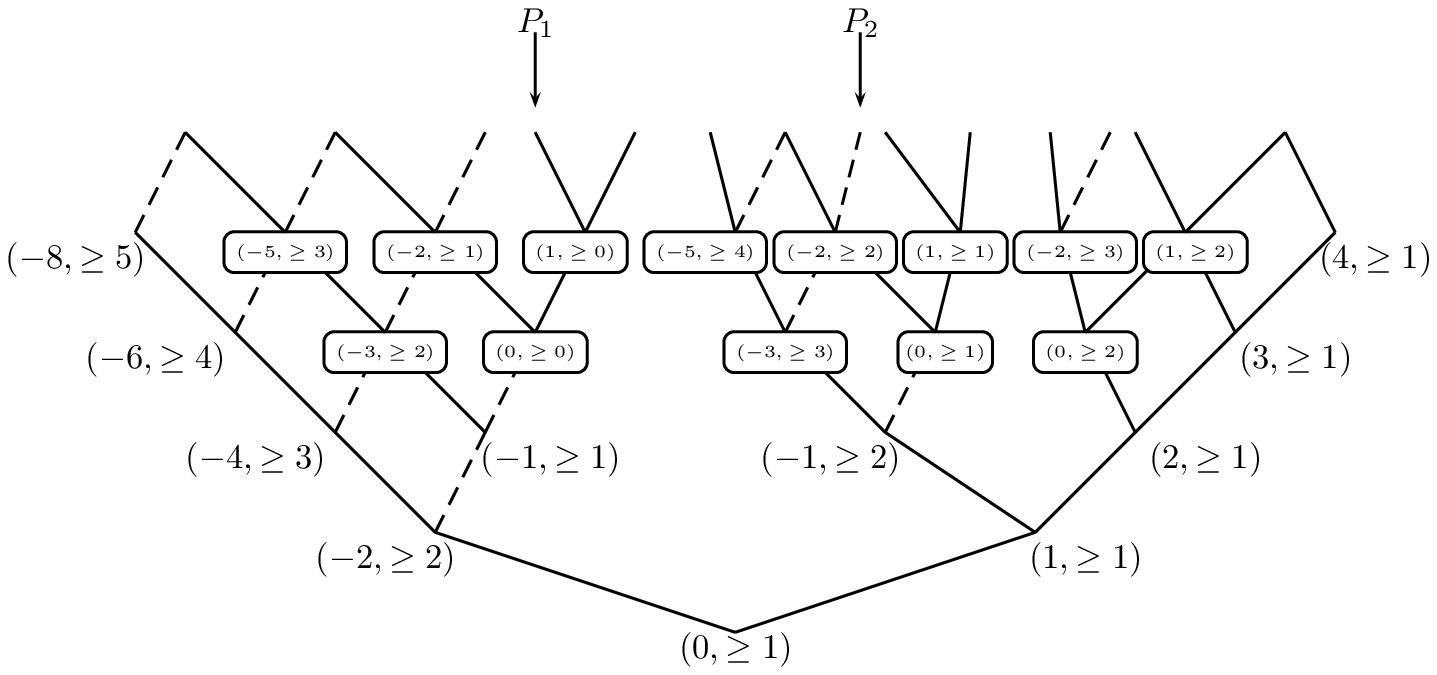}}
    \caption{Picture for Lemma \ref{lemma:deg5}}
    \protect\label{pic1lem6}
\end{figure}

\begin{proof}

    We must first confirm that the sequences $e_{2}e_{1}^{3}e_{2}$ and
    $e_{1}e_{2}e_{1}e_{2}e_{1}$ do not produce $0$ when applied to
    $T$.  First, note that there is at least one $+$ in the right
    block of the reduced $1$-signature of $T$.  Since $e_{2}$ acts on
    the right block of $+$'s in the $2$-signature of $T$, there are as
    many $+$'s in the right block of the $1$-signature of $e_{2}T$ as
    in that of $e_{1}T$.  Observe that $A(e_{2}T) = B(e_{2}T) - 2$, so
    we know that there are additionally two $+$'s in the left block of
    the reduced $1$-signature of $e_{2}T$.  This implies that
    $e_{1}^{3}e_{2}T \neq 0$.  The third of these applications of
    $e_{1}$ changes a $\bar{1}$ to a $\bar{2}$ in the top row; the $+$
    in the $2$-signature of $e_{1}^{3}e_{2}T$ entry cannot be
    bracketed, so we know that $e_{2}e_{1}^{3}e_{2}T \neq 0$.  On the
    other hand, we know that the right block of the reduced
    $2$-signature of $T$ has at least one $+$.  Since $e_{1}$ changes
    a $\bar{1}$ to a $\bar{2}$ in the top row of $T$, we know that
    $e_{2}$ will change this entry to a $2$ so that the right block
    of the reduced signature of $e_{1}T$ has at least two $+$'s that
    cannot be bracketed by $-$'s.  The leftmost of these entries will
    be acted upon by $e_{2}$, so $e_{2}e_{1}T \neq 0$.  Furthermore,
    since $A(e_{2}e_{1}T) = B(e_{2}e_{1}T) - 1$, we know that
    $e_{1}e_{2}e_{1}T \neq 0$.  At least one $+$ remains in the right
    block of the reduced $2$-signature of $e_{1}e_{2}e_{1}T$, so
    $e_{2}e_{1}e_{2}e_{1}T \neq 0$.  Finally, since
    $A(e_{2}e_{1}e_{2}e_{1}T) = B(e_{2}e_{1}e_{2}e_{1}T) - 2$, we know
    that $e_{1}e_{2}e_{1}e_{2}e_{1}T \neq 0$.
    
    Now that we know that neither of these sequences produces $0$ when
    applied to $T$, it is clear that we have $e_{2}e_{1}^{3}e_{2}T =
    e_{1}e_{2}e_{1}e_{2}e_{1}T$, since the paths $P_{1}$ and $P_{2}$
    leading from the base of the graph in Figure \ref{pic1lem6} to the
    indicated leaves have no solid left edges, two dashed left edges,
    one solid right edge, and two dashed right edges.  Note that these
    paths are equivalent to $e_{1}^{2}e_{2}^{2}e_{1}T$, due to the
    degree 2 relation above $e_{2}e_{1}T$; we may denote this
    alternative path by $P'_{2}$.  We must now confirm that among all
    pairs $(Q_{1},Q_{2})$ of increasing paths from the base in these
    graphs such that $Q_{1}$ begins by following the left edge and
    $Q_{2}$ begins by following the right edge, $(P_{1}, P_{2})$ and
    $(P_{1},P'_{2})$ are the only pairs with the same number of each
    of the above types of edges.

    Since the right edge from the base of the graph is solid, our
    candidate for $Q_{1}$ must have a solid right edge.  Observe that
    all paths in the left half of this graph with at least one solid
    right edge have two dashed right edges.  The only candidates for
    $Q_{2}$ with two dashed edges are $P_{2}$ and $P'_{2}$, both of
    which have two dashed left edges.  The only remaining candidate
    for $Q_{1}$ with two dashed left edges is in fact $P_{1}$.

\end{proof}

\begin{example}
    
    In Figure \ref{fig:lem6}, we have an interval of a crystal in
    which the bottom tableau $T$ has the statistics $A(T) = 2 = B(T)$
    and $C(T) = 1 \geq 1 = D(T) + 1$, illustrating Lemma \ref{lemma:deg5}.
    
    \begin{figure}[!hbp]
        \centerline{\includegraphics[height=3.55in]{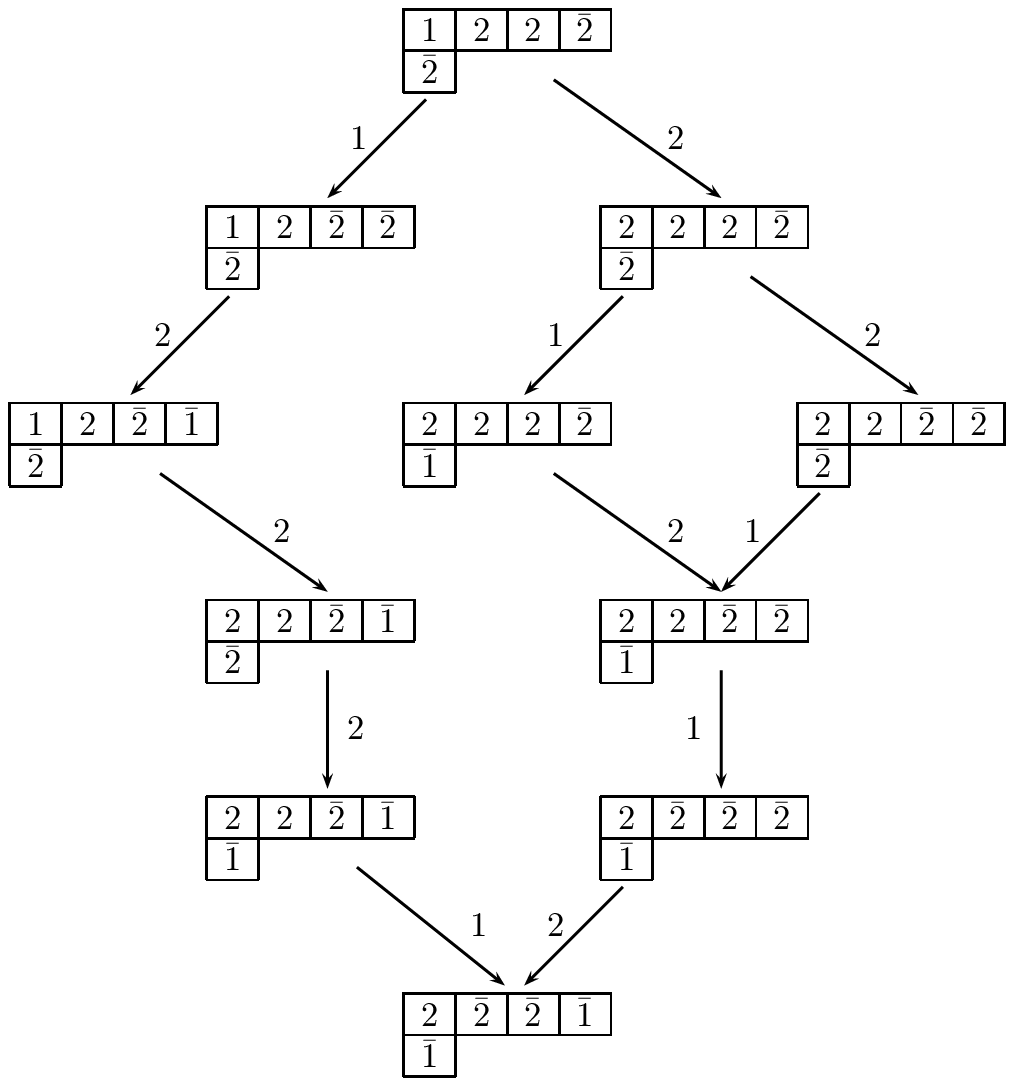}}
        \caption{}
        \protect\label{fig:lem6}
    \end{figure}
\end{example}


\begin{lemma}
    \label{lemma:deg7}
    Suppose $T$ is a tableau such that $A(T) = B(T)$, $C(T) = D(T)$,
    $e_{1}T\neq 0$, and $e_{2}T\neq 0$.  Then $T$ has a degree 7
    relation above it.
\end{lemma}

\begin{figure}[tbp]
    \centerline{\includegraphics[height=7in]{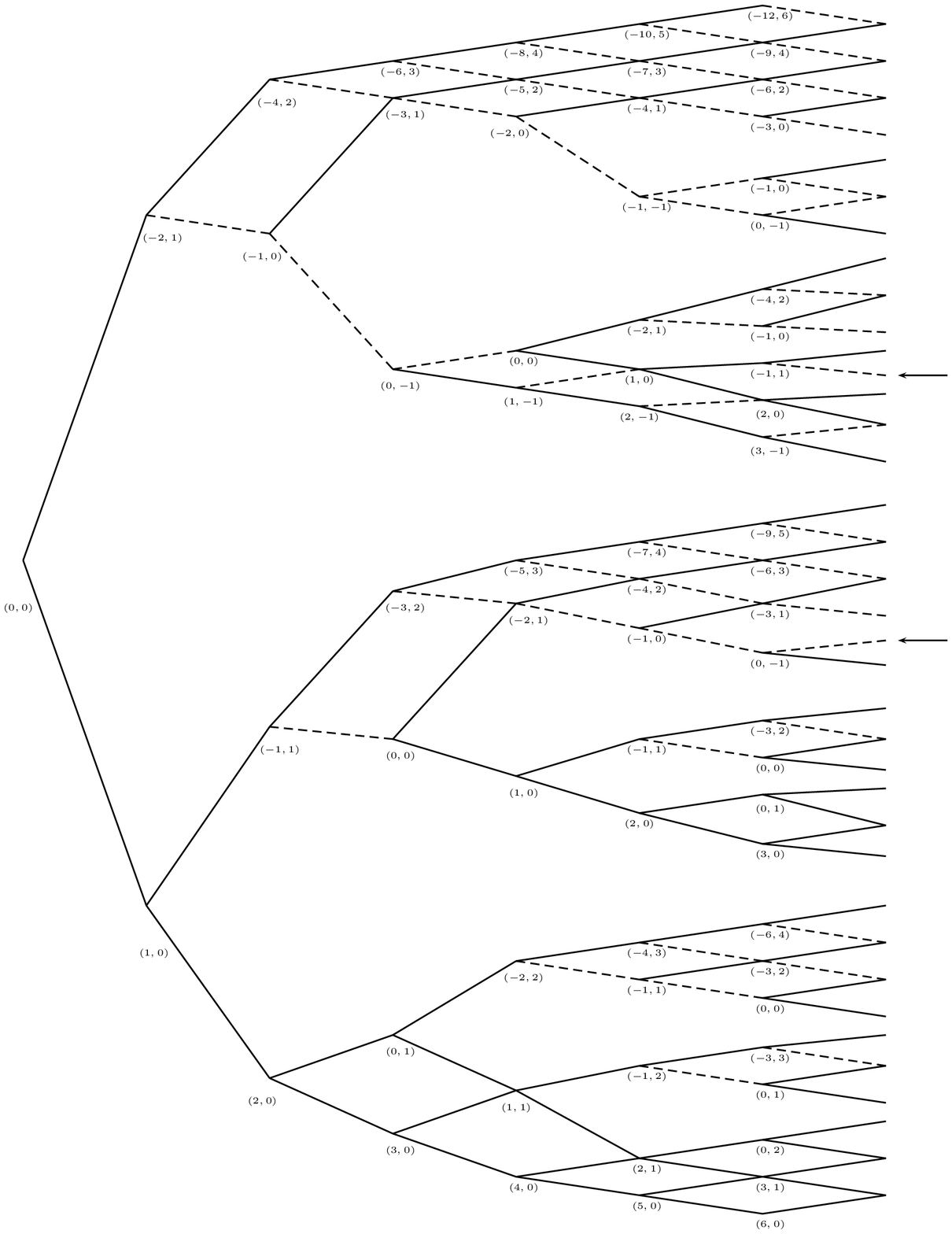}}
    \caption{Picture for Lemma \ref{lemma:deg7}}
    \protect\label{piclem7}
\end{figure}

\begin{proof}

    Note that in order to increase the readability of the graph in
    Figure \ref{piclem7}, it has been oriented to grow to the right
    rather than up.  We therefore take a down edge to indicate acting
    by $e_{1}$ and an up edge to indicate acting by $e_{2}$.
    Otherwise, the legend from Table \ref{piclegend} applies.

    We must first confirm that the sequences
    $e_{2}e_{1}^{2}e_{2}^{3}e_{1}$ and $e_{1}e_{2}^{3}e_{1}^{2}e_{2}$
    do not produce $0$ when applied to $T$.  By Sublemma
    \ref{sublemma:dashed}, we need only show that $e_{1}^{3}e_{2}T$,
    $e_{2}^{2}e_{1}^{3}e_{2}T$, and $e_{2}^{2}e_{1}T$ are not $0$.
    First note that there is at least one $\bar{1}$ in the top row of
    $T$, and the application of $e_{1}^{2}$ to $e_{2}T$ acts on
    entries corresponding to the left block of $+$'s.  It follows that
    the $\bar{1}$'s in the top row of $T$ are also present in
    $e_{1}^{2}e_{2}T$, so $e_{1}^{3}e_{2}T \neq 0$.  This final
    application of $e_{1}$ changes a $\bar{1}$ to a $\bar{2}$.  Since
    $e_{2}$ acts on the left block of $+$'s in $e_{1}^{3}e_{2}T$, it
    leaves this $\bar{2}$ alone, and it can be acted on by the next
    application of $e_{2}$, so $e_{2}^{2}e_{1}^{3}e_{2}T \neq 0$.
    Finally, note that there is a $\bar{2}$ in the top row of $T$ and
    $e_{1}$ changes a $\bar{1}$ to a $\bar{2}$ in the top row of $T$.
    Thus, there are at least two $\bar{2}$'s in the top row of
    $e_{1}T$, and $e_{2}^{2}e_{1}T \neq 0$.

    Now that we know that neither of these sequences produces $0$ when
    applied to $T$, it is clear that we have
    $e_{2}e_{1}^{2}e_{2}^{3}e_{1}T = e_{1}e_{2}^{3}e_{1}^{2}e_{2}T$,
    since the paths corresponding to these sequences leading from the
    base of the graph in Figure \ref{piclem7} to the leaves marked by
    arrows have one solid down edge, three dashed down edges, two
    solid up edges, and one dashed up edge.  Note that these paths are
    equivalent to $e_{2}e_{1}e_{2}e_{1}e_{2}^{2}e_{1}T =
    e_{1}e_{2}^{2}e_{1}e_{2}e_{1}e_{2}T$, due to the degree 2
    relations above $e_{2}^{2}e_{1}T$ and $e_{1}e_{2}T$; we denote
    these alternative paths by $P'_{1}$ and $P'_{2}$ respectively.  We
    must now confirm that among all pairs $(Q_{1},Q_{2})$ of
    increasing paths from the base in these graphs such that $Q_{1}$
    begins by following the up edge and $Q_{2}$ begins by following
    the down edge, $(P_{1}, P_{2})$, $(P_{1},P'_{2})$,
    $(P'_{1},P_{2})$ and $(P'_{1},P'_{2})$, are the only pairs with
    the same number of each of the above types of edges.

    We first address pairs of paths of length no greater than 5.  For
    a path to be a candidate for $Q_{1}$, it must have at least one
    solid down edge.  The only such paths are those beginning with
    $e_{1}^{2}e_{2}T$.  As these paths have two dashed down edges,
    their only possible $Q_{2}$ mate is $e_{1}^{2}e_{2}^{2}e_{1}T$,
    but none of our $Q_{1}$ candidates have the same edge content as
    this path.
    
    We now consider paths of length 6.  As in the preceding paragraph,
    our only candidates for $Q_{1}$ are those paths that contain a
    solid down edge and begin with $e_{1}^{2}e_{2}T$; all such paths
    have exactly two dashed down edges.  Up to degree 2 relations,
    there are three candidates for $Q_{2}$:
    $e_{1}^{2}e_{2}^{3}e_{1}T$, $e_{1}e_{2}e_{1}^{2}e_{2}e_{1}T$, and 
    $e_{1}^{2}e_{2}^{2}e_{1}^{2}T$.  None of these paths contain a
    dashed up edge, which leaves only $e_{1}^{5}e_{2}T$ as our only
    candidate for $Q_{1}$; this cannot be paired with any of our three
    potential $Q_{2}$ paths.

    Finally, we restrict our attention to paths of length 7.  There
    are six paths (again, up to degree 2 relations) in the top half of
    the graph with solid down edges: $e_{1}^{5}e_{2}^{2}T$ and those
    paths beginning with $e_{1}^{3}e_{2}T$.  The former has four
    dashed down edges, a feature lacking from all paths in the bottom 
    half of the graph.  We may also exclude from our consideration
    $e_{1}^{6}e_{2}T$, as all candidates for $Q_{2}$ with only one up 
    edge have at most one dashed down edge.
    
    The remaining three paths that might be $Q_{1}$ all have a dashed 
    up edge; the only $Q_{2}$ candidates with this feature are $P_{2}$
    and $P'_{2}$.  The only paths in the top half of the graph with
    the same edge content as these are $P_{1}$ and $P'_{1}$.

\end{proof}

\begin{example}
    
    In Figure \ref{fig:lem7}, we have a crystal in which the bottom
    tableau $T$ has the statistics $A(T) = B(T) = 1$ and $C(T) =
    D(T) = 0$, illustrating Lemma \ref{lemma:deg7}.
    
    \begin{figure}[!t]
        \centerline{\includegraphics[height=4.55in]{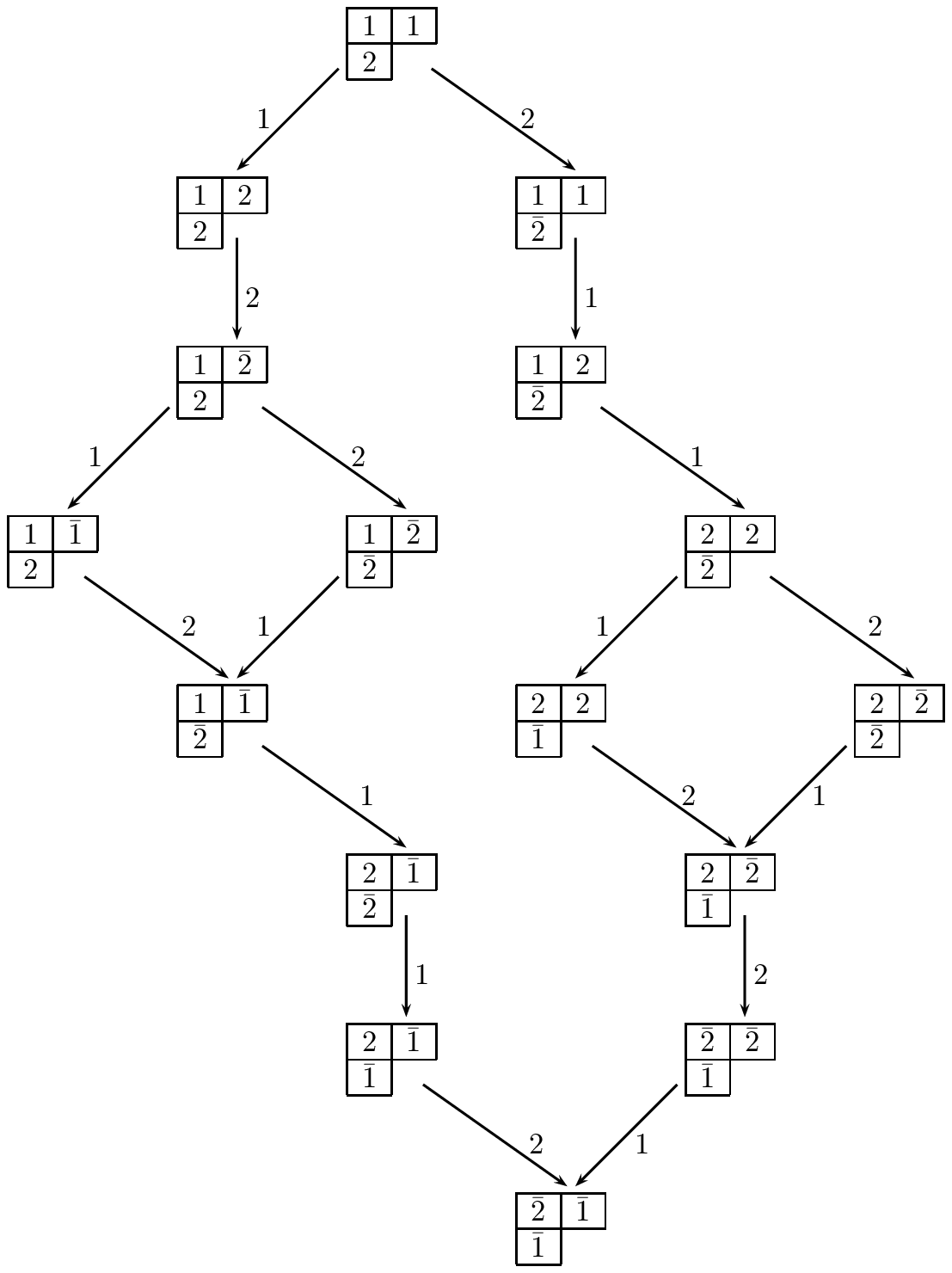}}
        \caption{}
        \protect\label{fig:lem7}
    \end{figure}
\end{example}

    %
    %

    \newchap{A User's Manual for CrystalView}
    \label{sec:CV}
    
        \section{Overview}

CrystalView is a software package for visualizing crystals for
irreducible highest weight modules of classical Lie algebras.  Based
on user input, the program will produce an image file with the
requested crystal graph.  The program will automatically produce an
.epsf file, which can be included directly in a \LaTeX~file.  Using
the web interface, the user can see the tableau associated to a vertex
of the crystal by moving the mouse pointer over it.

When running this program locally (as opposed to via the web), the
program produces a list of all tableaux in the specified crystal.
This can be used to carry out calculations on crystals, including
searches for tableaux with specified properties.  Additionaly,
Kashiwara operators may be applied to the tableaux.

All calculations on tableaux are carried out using python.  Image
files are automatically generated PostScript.  The web interface uses 
html, javascript, and css.

\section{Requirements}

CrystalView may be run using a web interface or from source code.  The
web interface can, in principle, be used from any browser that
supports form input (any ``modern'' browser).  However, some aspects
of the interface will only work with a browser that complies with
standard html, javascript, and css.  In particular, Internet Explorer
is known to have issues with some dynamic aspects of the web
interface.  Firefox is a recommended alternative.

To run the software from source code, the user must have Python 2.4
installed on their local machine.  Other versions of Python (both
older and newer) may not run CrystalView properly.  See {\tt
http://www.python.org} for further information regarding python.
Python is available free of charge for all major operating systems,
and is included pre-installed on many modern computers, including
almost all distributions of Unix/Linux.

If used to generate image files locally, the user is advised that for
large crystals, these images can get quite large.  See section
\ref{sec:legal}.

\section{User input}

The user may specify the following Lie theoretic data:
\begin{itemize}
    \item symmetry type of the algebra being represented;

    \item rank of the algebra being represented;

    \item highest weight of the representation.
\end{itemize}

Additionally, the user may specify how the edges of the crystal graph
will be drawn.  There are two default settings, color and grayscale,
as well as an option for custom settings.  If the custom option is
selected, the user may specify the following:
\begin{itemize}
    \item red/green/blue values for each edge color on a scale from
    $0$ to $1$ (Default = $1$),

    \item line width on a scale from $1$ to $5$ (Default = $5$),

    \item dash pattern; none, short, long (Default = none).
\end{itemize}
There are numerous resources on the web and preinstalled on many
computers to assist the user in finding red/green/blue values for
their desired colors.

The user may also choose to have the output converted to .pdf, .jpg,
.gif, and/or .tiff formats.  These conversions are carried out by
ImageMagick.

\section{Limitations}

Currently, only types $A$ and $B$ (i.e., $\mathfrak{sl}_{n}$ and
$\mathfrak{so}_{2n+1}$) are supported.  Furthermore, in the case of
type $B$, only even multiples of the highest fundamental weight
(corresponding to the short root) may be specified.  These limitations
are due to the current stage of the development cycle; future versions
of the software will add support for types $C$ and $D$ and all
dominant weights.

The rank of the algebra is currently restricted by the web interface
to be no larger than $5$.  This is an artificial limitation; any rank 
of algebra may be specified when running CrystalView from source.  

The web interface only allows crystals with as many as 4,000 vertices
to be calculated to prevent excessive strain on the server.
Considering the resolution at which these images can be
viewed/printed, it is unlikely that producing crystals larger than
this would be useful to most users.  However, this limitation is
artificial; when running the program on a local machine, the user is
limited only by their own patience and hard disk space.

\section{How it works}

The web interface for CrystalView is written in dhtml; i.e., html
enhanced by javascript and css.  The Weyl dimension formula is used to
calculate the number of vertices in the currently specified crystal.

The tableaux are produced by generating the list of all column
tableaux for the column lengths appearing in the shape specified by
the dominant weight.  To determine what constitutes a legal tableaux,
the criteria of \cite{HK} are used.  The columns are then compared
pairwise in order of decreasing length to build the set of all legal
tableaux.  In the case of type $B$ crystals, the ``split form'' 
criterion of \cite{L} is used to determine which columns can be 
adjacent in a legal tableau.  

The graph is ensured to have a reasonable number of edge crossings by
ordering the vertices of the graph as follows, starting from the top
and going down the rows, and proceeding through each row from left to
right.  First, the tableaux are collected into rows according to
content.  There is only one tableau in the first row (the highest
weight tableau), so the first row is in order.  Now, given that row
$n$ is in order, row $n+1$ is put in the following order.  The
leftmost vertices in row $n+1$ will be the non-zero images of the
Kashiwara operators $f_{i}$ on the leftmost vertex in row $n$, taken
in order from $f_{1}$ up through $f_{r}$, where $r$ is the rank of the
algebra.  The next leftmost vertices in row $n+1$ are those tableaux
that result from applying $f_{i}$ to the next vertex from row $n$,
excluding those that have been placed to the left already.  This
process continues until all vertices have been placed in their final 
position.

\section{How to get the software}

This software can be accessed at the following url:

{\tt http://www.math.ucdavis.edu/$\sim$sternberg/crystalview/}\\
The source code is available upon request from the author.

\section{Examples of output}

\begin{figure}[!hbp]
    \centerline{\includegraphics[height=2.95in]{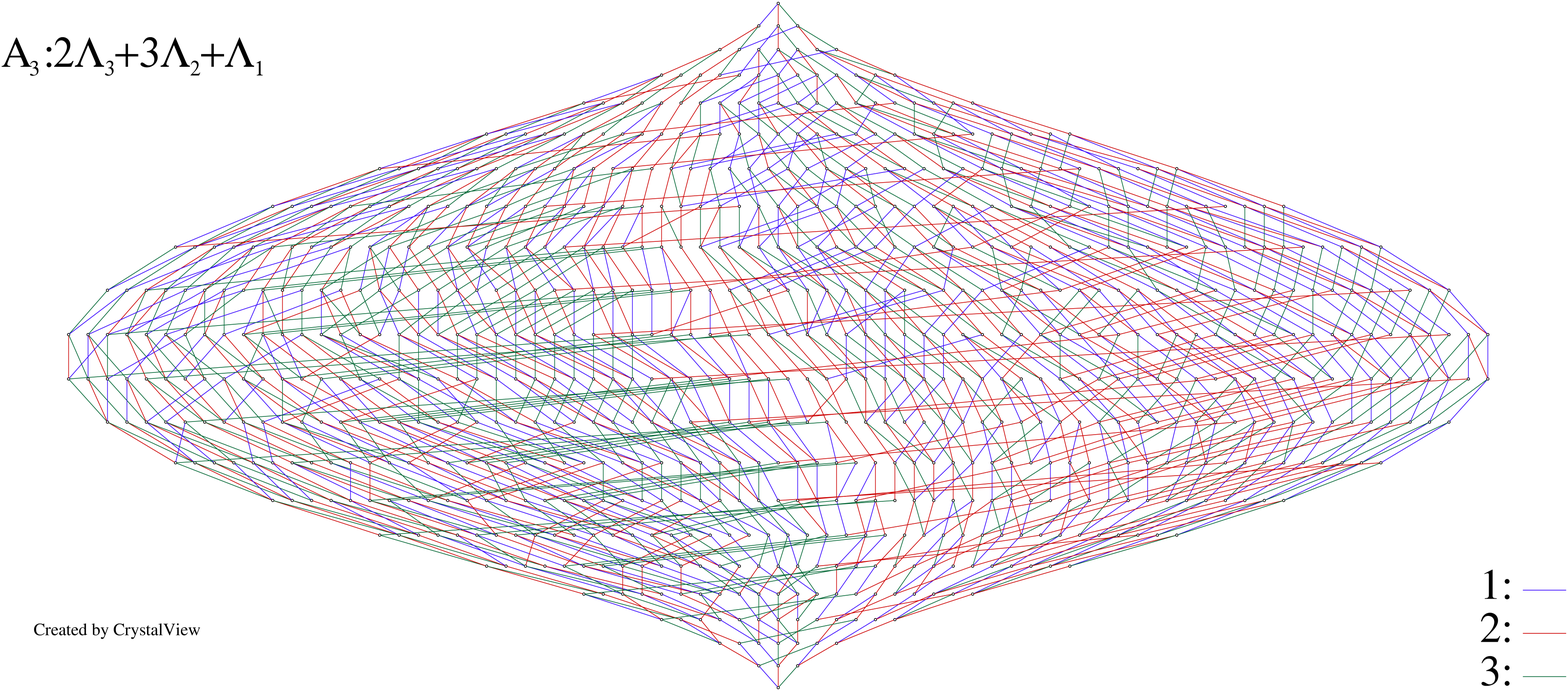}}
    \caption{}
    \protect\label{fig:cv1}
\end{figure}

\begin{figure}[!hbp]
    \centerline{\includegraphics[height=5.5in]{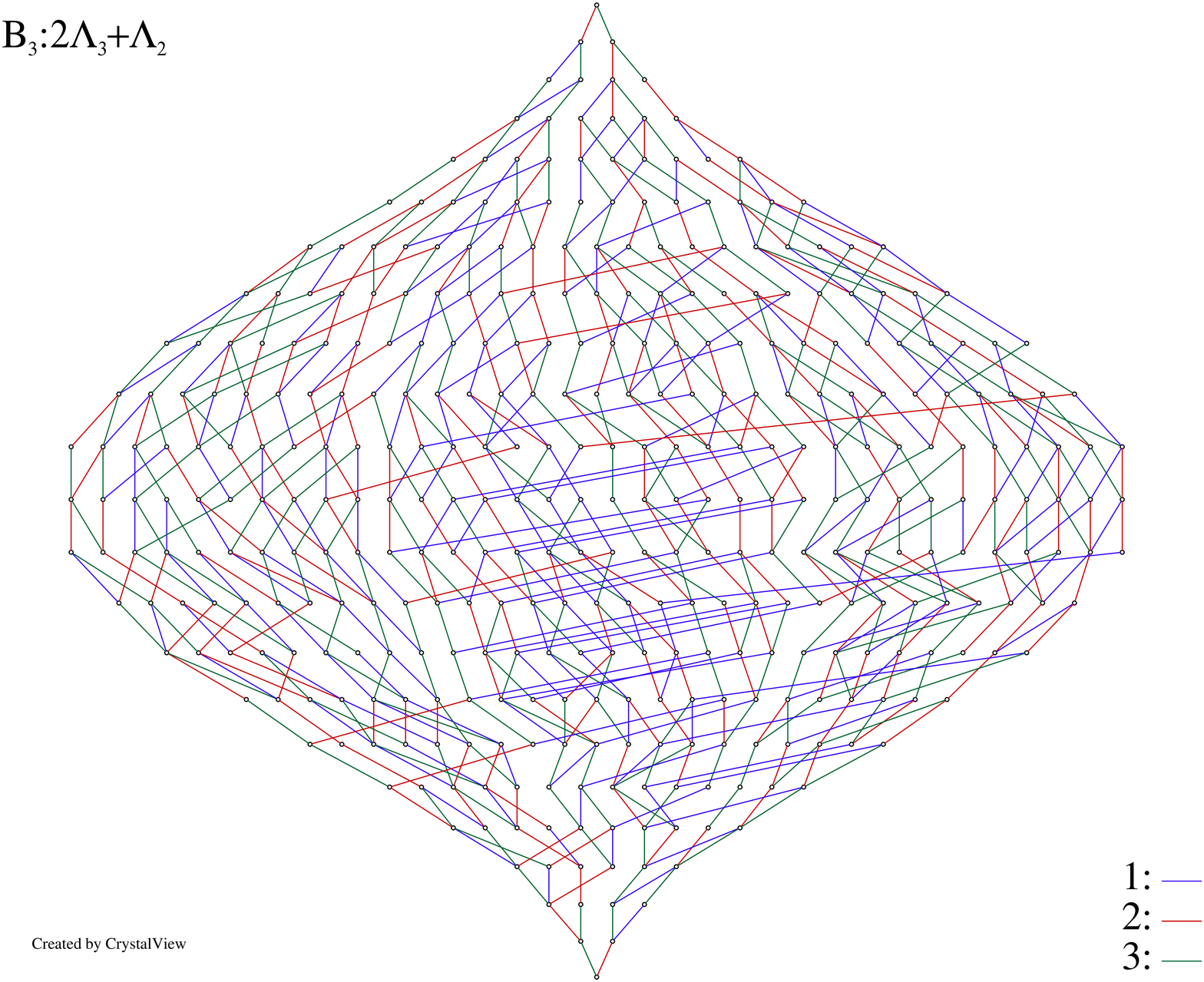}}
    \caption{}
    \protect\label{fig:cv2}
\end{figure}

Figures \ref{fig:cv1} and \ref{fig:cv2} feature examples of the
graphical output of CrystalView.
    
\section{Legal disclaimer}\label{sec:legal}

This software is provided without warranty; use is at the user's sole 
risk.  Under no circumstances will any user hold the author of this 
software liable for any damages that may result from the use of this 
software.  Use of the software implies that the user agrees to these 
terms.

    %
    %

\end{document}